\documentclass[nobib,a4paper,symmetric,justified]{tufte-handout}

\hypersetup{%
  citecolor = Black,
  linkcolor = Black,
  urlcolor = Black,
}%

\usepackage{graphicx} 
  \setkeys{Gin}{width=\linewidth,totalheight=\textheight,keepaspectratio}
  \graphicspath{{graphics/}} 
  
\usepackage{array}
\usepackage{mdframed}
\usepackage[fleqn]{amsmath} 
\usepackage{amsfonts,amssymb,amsthm,thmtools}
\DeclareSymbolFontAlphabet{\mathbb}{AMSb}
\usepackage[most]{tcolorbox} 
\usepackage{bbding} 
\usepackage{stmaryrd} 
\usepackage{booktabs} 
\usepackage{units}    
\usepackage{multicol} 
\usepackage{lipsum}   
\usepackage{fancyvrb} 
  \fvset{fontsize=\normalsize}
\usepackage{scalerel}
\usepackage{enumitem}
\setlist[enumerate]{nosep, label=\arabic*), leftmargin=1.44em, labelsep=4pt}
\setlist[description]{nosep}
\let\marginnote\relax

\usepackage{marginnote}

\newcommand{\mathnote}[1]{\marginnote{\raggedright #1}[-1.137\baselineskip]}
\newcommand{\proofnote}[1]{\marginnote{\raggedright #1}[-1.39\baselineskip]}

\newcommand{\s}[1][1]{\mkern #1mu}
\newcommand\w{\s[2]}%

\newcommand\mquad{\s[-18]}
\newcommand\os{\s[2]\text{--}\s[4]}
\newcommand{\vl}{\w\vert\w}
\newcommand{\cs}{,\s[1]}
\newcommand{\Ref}[1]{\scriptstyle{\s[-5]\ref{#1}\s[5]}}%

\newcommand{\Deq}{\color{LightGray}\blacksquare\,}
\newcommand\Qed{$\color{lightgray}\qedsymbol$}%

\newcommand{\DEQ}{$\blacksquare$\nopunct}
\newcommand{\QED}{\qedhere}
\newcommand{\Item}[1]{\s[8]{\scriptstyle \text{#1)}}\s[8]}

\newcommand{\Point}{\text{\rm Point}}
\newcommand{\Line}{\text{\rm Line}}
\newcommand{\Length}{\text{\rm Length}}
\newcommand{\Angle}{\text{\rm Angle}}
\newcommand{\Measure}{\text{\rm Measure}}
\newcommand{\Figure}{\text{\rm Figure}}
\newcommand{\Ray}{\text{\rm Ray}}
\newcommand{\Flag}{\text{\rm Flag}}
\newcommand{\Postulate}[1]{\text{P$_{\s {#1}}$}}
\newcommand{\Axiom}[1]{\text{A$_{\s {#1}}$}}

\newcommand{\Solution}[1]{\text{S$_{\s {#1}}$}}
\newcommand{\CC}{\blacktriangleright}

\definecolor{MyRed}{RGB}{200,100,50} 
\definecolor{MyBlue}{RGB}{40,100,160} 
\definecolor{MyGreen}{RGB}{10,100,10}
\definecolor{MyGray}{RGB}{220, 220, 220}
\definecolor{BlockGray}{rgb}{1, 1, 1}
\newcommand{\Def}{\color{MyRed}} 
\newcommand{\Prop}{\color{MyBlue}} 
\newcommand{\Set}{\color{MyGreen}} 
\newcommand{\Mrg}{\color{Black}} 

\newcommand{\refcolor}[1]{\hypersetup{linkcolor=#1}} 

\newtheorem{definition}{\Def Definition}
\newtheorem{axiom}{\Prop Axiom}
\newtheorem{postulate}{\Set Postulate}
\newtheorem{problem}{\Set Problem}
\newtheorem{theorem}{\Prop Theorem}

\newenvironment{axiombis}
  {%
   \addtocounter{axiom}{-1}%
   \begin{axiom}}
  {\end{axiom}}

\tcolorboxenvironment{definition}{
  oversize,
  if odd page={leftrule=4pt}{leftrule=0.2pt},
  if odd page={rightrule=0.2pt}{rightrule=4pt},
  toprule=0.2pt,
  bottomrule=0.2pt,
  colback=MyRed!0!BlockGray,
  colframe=MyRed,
  arc=0pt,
  outer arc=0pt,
  titlerule=0pt,
  bottomtitle=0pt,
  top=6pt,
  bottom=0pt,
  left=6pt,
  right=6pt,
}
\tcolorboxenvironment{axiom}{
  oversize,
  if odd page={leftrule=4pt}{leftrule=0.2pt},
  if odd page={rightrule=0.2pt}{rightrule=4pt},
  toprule=0.2pt,
  bottomrule=0.2pt,
  colback=MyBlue!0!BlockGray,
  colframe=MyBlue,
  arc=0pt,
  outer arc=0pt,
  titlerule=0pt,
  bottomtitle=0pt,
  top=6pt,
  bottom=0pt,
  left=6pt,
  right=6pt,
}
\tcolorboxenvironment{postulate}{
  oversize,
  if odd page={leftrule=4pt}{leftrule=0.2pt},
  if odd page={rightrule=0.2pt}{rightrule=4pt},
  toprule=0.2pt,
  bottomrule=0.2pt,
  colback=MyGreen!0!BlockGray,
  colframe=MyGreen,
  arc=0pt,
  outer arc=0pt,
  titlerule=0pt,
  bottomtitle=0pt,
  top=6pt,
  bottom=0pt,
  left=6pt,
  right=6pt,
}

\makeatletter
\newsavebox{\@brx}
\newcommand{\llangle}[1][]{\savebox{\@brx}{\(\m@th{#1\langle}\)}%
\mathopen{\copy\@brx\kern-0.675\wd\@brx\usebox{\@brx}}}
\newcommand{\rrangle}[1][]{\savebox{\@brx}{\(\m@th{#1\rangle}\)}%
\mathclose{\copy\@brx\kern-0.675\wd\@brx\usebox{\@brx}}}

\newcommand{\llpar}[1][]{\savebox{\@brx}{\(\m@th{#1(}\)}%
\mathopen{\copy\@brx\kern-0.675\wd\@brx\usebox{\@brx}}}
\newcommand{\rrpar}[1][]{\savebox{\@brx}{\(\m@th{#1)}\)}%
\mathclose{\copy\@brx\kern-0.675\wd\@brx\usebox{\@brx}}}

\newcommand{\llcurb}[1][]{\savebox{\@brx}{\(\m@th{#1\{}\)}%
\mathopen{\copy\@brx\kern-0.725\wd\@brx\usebox{\@brx}}}
\newcommand{\rrcurb}[1][]{\savebox{\@brx}{\(\m@th{#1\}}\)}%
\mathclose{\copy\@brx\kern-0.725\wd\@brx\usebox{\@brx}}}

\newcommand\bbs{\mathpalette\bbs@{.5}}%
\newcommand\bss{\mathpalette\bbs@{.6}}%
\newcommand\bbs@[2]{\mathbin{\vcenter{\hbox{\scalebox{#2}{$\m@th#1\text{\FiveStar}$}}}}}%
\newcommand\wbs{\mathpalette\wbs@{.5}}%
\newcommand\wbs@[2]{\mathbin{\vcenter{\hbox{\scalebox{#2}{$\m@th#1\text{\FiveStarOpen}$}}}}}%
\newcommand{\bs}{\protect\bbs}%
\newcommand{\ws}{\protect\wbs}%

\newcommand\trg{\mathpalette\trg@{.65}}%
\newcommand\trg@[2]{\mathbin{{\hbox{\scalebox{#2}{$\m@th#1\triangle$}}}}}%
\newcommand{\tr}{\protect\trg}%

\newcommand{\Obj}[2]{#1 : #2}%
\newcommand{\Forall}[3]{\forall\w\Obj{#1}{#2}, \s #3 \s}%
\newcommand{\Exists}[3]{\exists\w\Obj{#1}{#2}, \s #3 \s}%
\newcommand{\NotExists}[3]{\nexists\w\Obj{#1}{#2}, \s #3 \s}%
\newcommand{\Build}[3]{\{\s\Obj{#1}{#2} \vl \s #3\s\}}%
\newcommand{\BuildDouble}[5]{\{(\s \Obj{#1}{#2}\s)(\s \Obj{#3}{#4}\s) \vl \s #5\s\}}%
\newcommand{\Or}[2]{ \{\s #1 \s\} \s[-2]+\s[-2] \{\s #2 \s\}}%
\newcommand{\ObjDouble}[5]{(\s \Obj{#1}{#2}\s)(\s \Obj{#3}{#4}\s), \s #5 \s}%
\newcommand{\ForallDouble}[5]{\forall\w(\s \Obj{#1}{#2}\s)(\s \Obj{#3}{#4}\s), \s #5 \s}%
\newcommand{\In}[2]{(\s #1 \in #2 \s)}%
\newcommand{\NotIn}[2]{(\s #1 \notin #2 \s)}%

\newcommand{\Apart}[2]{[ \w #1 \bs #2 \w ]}%
\newcommand{\ColTwo}[2]{[ \w #1 \s #2 \w ]}%
\newcommand{\ColThree}[3]{[ \w #1 \s #2 \s #3 \w ]}%
\newcommand{\Bet}[3]{[ \w #1 \bs #2 \bs #3 \w ]}%
\newcommand{\BetFour}[4]{[ \w #1 \bs #2 \bs #3 \bs #4 \w ]}%
\newcommand{\BetOX}[3]{[ \w #1 \ws #2 \bs #3 \w ]}%
\newcommand{\BetXO}[3]{[ \w #1 \bs #2 \ws #3 \w ]}%
\newcommand{\BetOO}[3]{[ \w #1 \ws #2 \ws #3 \w ]}%
\newcommand{\SameRay}[3]{[ \w #1 \bs #2 \s #3 \w ]}%
\newcommand{\SameRayO}[3]{[ \w #1 \ws #2 \s #3 \w ]}%
\newcommand{\OppositeSide}[3]{[ \w #1 \vl #2 \vl #3 \w ]}%
\newcommand{\SameSide}[3]{[ \w #1 \vl #2 \s #3 \w ]}%
\newcommand{\SameHalfplane}[4]{[ \w #1 \bs #2 \vl #3 \s #4 \w ]}%
\newcommand{\Middle}[3]{[ \w #1 \os #2 \os #3 \w ]}%
\newcommand{\Seg}[2]{#1 #2}%
\newcommand{\Circ}[2]{#1 #2}%
\newcommand{\Lin}[2]{#1 #2}%
\newcommand{\RayObj}[2]{\s \overrightarrow{#1 #2} \s }%
\newcommand{\Len}[2]{\llbracket \s #1 #2 \s \rrbracket}%
\newcommand{\IntAngPs}[3]{{\llangle \s #1 #2 #3 \s \rrangle}}%
\newcommand{\AngRs}[2]{{\llangle \s #1 \cdot #2 \s \rrangle}}%
\newcommand{\IntAngRs}[2]{{\llangle \s #1 \bs #2 \s \rrangle}}%
\newcommand{\LC}[3]{\Len{#1}{#2} \s[-3]=\s[-3] #3}%
\newcommand{\AC}[3]{\AngRs{#1}{#2} \s[-3]=\s[-3] #3}%
\newcommand{\Lcong}[4]{\Len{#1}{#2} \s[-3]=\s[-3] \Len{#3}{#4}}%
\newcommand{\LessPlus}[6]{\Len{#1}{#2} \s[-3]<\s[-3] \Len{#3}{#4} \s[-3]+\s[-3] \Len{#5}{#6}}%
\newcommand{\LcongThree}[6]{\Len{#1}{#2} \s[-3]=\s[-3] \Len{#3}{#4}\s[-3]=\s[-3] \Len{#5}{#6}}%
\newcommand{\IntAngCongPs}[6]{\IntAngPs{#1}{#2}{#3} \s[-3]=\s[-3] \IntAngPs{#4}{#5}{#6}}%
\newcommand{\AngCongRs}[4]{\AngRs{#1}{#2} \s[-3]=\s[-3] \AngRs{#3}{#4}}%
\newcommand{\Triangle}[3]{{\tr\s[1]#1 #2 #3}}%
\newcommand{\CongTriangles}[6]{\Triangle{#1}{#2}{#3} \cong \Triangle{#4}{#5}{#6}}%
\newcommand{\Prime}[1]{#1'\s[-0]}%
\newcommand{\Op}[1]{{#1}^{*}}%
\newcommand{\Flip}[1]{\widehat{#1}}%
\newcommand{\Iff}{\s[4]\xleftrightarrow{\hspace*{2.25mm}}\s[4]}%

\renewenvironment{proof}[1][\proofname]{\par
  \pushQED{\qed}%
  \normalfont  \small
  \topsep=0pt 
  \partopsep=0.5\baselineskip
  \trivlist
  \item[\hskip\labelsep
        \itshape
    #1\@addpunct{.}]\ignorespaces
}{%
  \popQED\endtrivlist\@endpefalse
\addvspace{-0.5\baselineskip} 
}

\newcommand*\bigcdot{\mathpalette\bigcdot@{.5}}
\newcommand*\bigcdot@[2]{\mathbin{\vcenter{\hbox{\scalebox{#2}{$\m@th#1\bullet$}}}}}

\newsavebox\myboxA
\newsavebox\myboxB
\newlength\mylenA

\newcommand*\xoverline[2][0.75]{%
    \sbox{\myboxA}{$\m@th#2$}%
    \setbox\myboxB\null
    \ht\myboxB=\ht\myboxA%
    \dp\myboxB=\dp\myboxA%
    \wd\myboxB=#1\wd\myboxA
    \sbox\myboxB{$\m@th\overline{\copy\myboxB}$}
    \setlength\mylenA{\the\wd\myboxA}
    \addtolength\mylenA{-\the\wd\myboxB}%
    \ifdim\wd\myboxB<\wd\myboxA%
       \rlap{\hskip 0.5\mylenA\usebox\myboxB}{\usebox\myboxA}%
    \else
        \hskip -0.5\mylenA\rlap{\usebox\myboxA}{\hskip 0.5\mylenA\usebox\myboxB}%
    \fi}

\newcommand\xleftrightarrow[2][]{%
  \ext@arrow 9999{\longleftrightarrowfill@}{#1}{#2}}
\newcommand\longleftrightarrowfill@{%
  \arrowfill@\leftarrow\relbar\rightarrow}

\makeatother

\title{\sc On Constructive-Deductive Method\\For Plane Euclidean Geometry\\}

\author{Evgeny V. Ivashkevich}
\date{{\footnotesize\rm\vskip-2mm Moscow Institute of Physics and Technology, Russia}\marginnote[0.001mm]{\rm E-mail: ivashkev@yandex.ru}} 

\begin{document}

\maketitle

\begin{abstract}\noindent
Constructive-deductive method for plane Euclidean geometry is proposed and formalized within Coq Proof Assistant.
This method includes both {\em postulates} that describe elementary constructions by idealized geometric tools (pencil, straightedge and compass),
and {\em axioms} that describes properties of basic geometric figures (points, lines, circles and triangles).
The proposed system of postulates and axioms can be considered as a constructive version of the Hilbert's formalization of plane Euclidean geometry.
\end{abstract}

\section{0. Introduction}

\subsection{0.1. Logic of Euclid}

\marginnote{\raggedright
{\Set\sc\normalsize constructive method}\\[1.0\baselineskip]
{\Set\bf Postulate} --- prescribes how to construct an elementary object with the required properties.\\[1.0\baselineskip]
{\Set\bf Problem} --- is a specification that formalizes required properties for the object to be constructed.\\[1.0\baselineskip]
{\Set\bf Construction} --- is a method of acting (using set of tools) by which problems are solved from postulates.\\[1.0\baselineskip]
{\Set\bf Solution} --- is a sufficient evidence for the decidability of a problem.}
Euclidean geometry \cite{Euclid} studies not only the properties of geometric figures, but also the rules of their construction. Propositions, that we encounter in Euclid's {\em Elements},  can be divided onto two types: problems and theorems. 

To solve {\em a problem} one needs to find a sequence of elementary actions ({\em postulates}) that leads to the {\em construction} of a figure with the required properties. 
Once the construction algorithm is found, it gives us a {\em solution} to the problem.
The set of allowed elementary actions is completely determined by chosen geometric tools (a pencil, a straightedge and a compass) and provides us with the instructions on how to handle these instruments. 

Pencil helps us to draw some points on a plane. Straightedge helps us to draw a straight lines passing through two distinct points. Compass helps us to construct a circle through a given point and the center of the circle. 

\marginnote{\raggedright
{\Prop\sc\normalsize deductive method}\\[1.0\baselineskip]
{\Prop\bf Axiom} --- asserts some inherent property of an elementary objects.\\[1.0\baselineskip]
{\Prop\bf Theorem} --- is a proposition that formalizes some derived properties of  already constructed objects.\\[1.0\baselineskip]
{\Prop\bf Deduction} --- is a method of reasoning (using set of rules) by which theorems are proved from axioms.\\[1.0\baselineskip]
{\Prop\bf Proof} --- is a sufficient evidence for the truth of a theorem.}
To prove {\em a theorem} one needs to find a sequence of elementary inferences ({\em axioms}) that leads to the {\em deduction} of a certain property in the given figure. 
Once the deductive inference is found, it gives us a {\em proof} to the theorem.
The set of allowed elementary inferences describes some basic properties of figures (points, lines, circles and triangles) that were constructed using the geometric tools that we have chosen. 

Besides axioms we also need some rules of logic (e.g. {\em Modus Ponens}) to prove our theorems. Nevertheless, it is the axioms, not the rules of logic, that fill our theory with geometric meaning.

A distinctive feature of Euclid's logic, is that it is not purely deductive \cite{Rodin2014, Rodin2018, Melikhov2015}.
Postulates and problems represent an important constructive part of Euclid's logic which cannot be completely separated from the deductive one. Indeed, proofs of Euclid's theorems often require auxiliary constructions, while problem solutions often require intermediate proofs to ensure applicability of postulates or previously solved problems. So, the inherent logic of Euclid's {\em Elements} can be considered as constructive-deductive. 

Later interpreters of Euclid, however, sought to reduce all his problems into existential theorems. In other words, they treated all objects of the theory as idealized and pre-existing and considered every problem not as instruction on how to construct the object with the required properties, but simply as propositional proof of object's existence.

This forced reduction of constructive-deductive logic of Euclid to purely propositional logic was not without its difficulties. Indeed, in the logic of Euclid, every object under consideration must first be constructed. While in purely propositional logic the mere existence of an object with the required properties can be proved indirectly. One first assumes that such an object does not exist, and then prove that this assumption leads to a contradiction. Such an indirect proof by contradiction relies upon so-called "the law of excluded middle" first formulated by Aristotle. This law states that either an assertion or its negation must be true.

Intuitionistic logic \cite{Heyting1980, vonPlato1995, Beeson2015b} was proposed prohibiting all indirect proofs of existence. The main idea of the method was to reduce further not only the constructive-deductive logic of Euclid, but also the classical propositional logic. Namely, {\em intuitionists} suggested to exclude "the law of excluded middle" from the inference rules of propositional logic. As a results, from the existential proofs of the intuitionistic logic one can reproduce original Euclid's constructions and vice versa.

Despite of the apparent success, this reduction of classical logic cannot be considered satisfactory.
Indeed, let us recall that Euclid builds his theory using a limited set of geometric instruments --- a pencil, a straightedge and a compass.
This means, that translating postulates and problems into existential axioms and theorems, and then limiting ourselves by the rules of intuitionistic logic, we shall be able to prove existence of only those geometric figures that can be constructed with this limited set of geometric tools. However, it is well known that not all geometric figures can be constructed by these simple tools. For example, it is impossible to do a trisection of an angle. Obviously, this does not mean that the rays that divide an angle onto three equal parts do not exist. To get around this new difficulty, intuitionists would have to change the original Euclidean postulates, thereby introducing some idealized tools that would give us the ability to construct more complex geometric figures.

\subsection{0.2. Hilbert's system of axioms}

The standards of mathematical rigor have changed dramatically since Euclid wrote his {\em Elements}. Subsequently, many mathematicians added missing but implicitly used statements to the original list of Euclid's postulates and axioms. The culmination of this development was a small book by german mathematician David Hilbert entitled {\em Foundations of geometry} \cite{Hilbert1988}.  Although the overall approach of Hilbert is fairly close in spirit to that of Euclid, yet in some respects it deviates substantially from it. Let us discuss these differences.

First of all, we have to note that Hilbert does not adhere to the constructive-deductive logic of Euclid. In fact, Hilbert was one of the main proponents of the existential approach and hence reduced all the postulates and problems of Euclid's {\em Elements} to purely existential statements within classical logic.

Another striking feature of the Hilbert's axioms is the complete absence of circles. Instead of Euclid's tools (pencil, straightedge and compass), Hilbert introduces some imaginary instruments for existential transfer of both segments and angles. Surprisingly, those same authors who criticized Euclid for his method of superposition (which is nothing more than an imaginary tool for existential transfer of triangles) seem to accept Hilbert's tool for similar angle transfer without any questions. 

Since without compass it is not possible to reproduce many propositions from Euclid's {\em Elements}, Hilbert followers \cite{Greenberg1993, Hartshorne2013} have to introduce some additional "principles": such as line-circle intersection and circle-circle intersection. These "principles" are nothing but additional postulates that describe simple constructions by straightedge and compass. However, those authors do not add these "principles" to the list of Hilbert's axioms. Instead, they treat them as being derived from the Hilbert continuity axiom. Needless to say, that if such a derivation is actually provided by those authors, it goes far beyond Euclid's original approach.

As for the Hilbert's treatment of angles, it cannot be considered completely satisfactory either. The thing is that, Hilbert defines only interior angles (greater than null angle and less than straight angle). As a result of this truncated definition, it is impossible to universally define the addition of two arbitrary interior angles. Indeed, the sum of the two interior angles in Hilbert's theory is defined only iff the resulting angle is also interior. Thus, the addition of angles in Hilbert's theory does not form a group. 

In the last edition of his book, to define area of an arbitrary polygon Hilbert introduces the idea of triangle orientation. However, Hilbert did not provide the reader with any clear definition of the triangle orientation and has not studied the properties of the concept within his system of axioms. 

Hilbert introduces a number of separate congruence axioms for segments and angles. In Euclid's approach, on the contrary, the lengths of segments , angular measures and areas were all considered as {\em quantities} that obey the same axioms: "two quantities equal to the third are equal to each other", etc. In other words, for Euclid, "quantities" are more general concepts than segments and angles themselves. As such, Euclid considers "quantities" as classes o congruence. 

\subsection{0.3. Motivation and results}

In this paper, we return to Euclid's original constructive-deductive logic. This means that we are going to use classical logic for propositions and intuitionistic logic for constructions. In other words, we shall clearly distinguish between existential statements and real constructions. The same applies to disjunction. We shall consider both propositional disjunctions which is a logical assertion that one of the two alternatives is true without specifying which one, and constructive disjunction, which is a practical algorithm that for every configuration of objects decides which one of the two alternatives actually takes place.

We will show that the {\em Calculus of Constructions} (CoC), as it was implemented in the Coq Proof Assistant, suits ideally for the constructive-deductive logic of Euclid. In the CoC we have two sorts: {\tt \Prop Prop} and {\tt \Set Set}. The first one is inhabited by propositions and their proofs, while the second by constructed objects and their specifications (required properties). Without contradiction we may accept classical logic for proofs in the {\tt \Prop Prop} sort, and at the same time limit ourselves to intuitionistic logic for constructions in the {\tt \Set Set} sort. As a result, our constructions will be limited to only those that can be performed with our simple tools --- a pencil, a straightedge and a compass, while we shall be able to prove propositional existence of more complex geometric figures.

We abandon Hilbert's axioms for existential segment and angle transfer ($\text{III}_1$ and  $\text{III}_4$) and replace them by simple versions of line-circle and circle-circle intersection postulates (Postulate \ref{pst:DrawIntersectionOfLineAndCircle} and \ref{pst:DrawIntersectionOfTwoCircles} below). Detailed comparison of Hilbert's and proposed constructive systems of axioms is given at the end of the paper.

Constructive-deductive logic of Euclid also gives us the opportunity to formalize Euclid's "method of superposition". To prove triangle congruence theorems, Euclid conjectured that every triangle could be moved and placed over another triangle so that their respective sides and angles could be compared. In this interpretation, the method of superposition seems to be quite a strong postulate, which makes it possible to move, rotate and flip triangles on the plane. That is why it is not very popular nowadays, and it was questioned even in the times of Euclid. In our interpretation the method of superposition turns into an axiom about the existence of an equal triangle anywhere on the plane. The actual construction of equal triangle can be done later only with the help of pencil, straightedge and compass.

We define angles as a pairs of rays originating from the same point and prove that addition of angles can be defined uniformly for all angles (modulo the full angle). Then we define the orientation of angles on the plane and use this concept to define a subclass of convex angles, thus returning to the Hilbert's definition of angle. As a results, all Hilbert axioms can be derived form our formalization of Euclid's constructive-deductive method. 

We replaced separate congruence axioms for segments and angles by the classes of congruences and consider the classes as definitions of Euclid's quantities, that is segment lengths and angle measures. Such an approach simplifies our formalization of Euclid's constructive-deductive method in the Coq Proof Assistant. 

All our development has been formalized in the Coq Proof Assistant \cite{Ivashkevich2019} and we demonstrate that from our set of axioms and postulates it is possible to reproduce all the Von Plato (for constructive incidence geometry) and Hilbert axioms exactly in the form proposed in the respectively formalizations by Gilles Kahn \cite{Kahn2008} and Julien Narboux et al \cite{Narboux2015}. 

\section{1. Points}

\subsection{1.1. Conceptualization of point}

The simplest geometric figures on the plane are points.
Euclid described point as such a figure that can not be divided into parts.
In other words, according to the ancient geometers, a point on the plane has no dimensions: neither length nor width.
In fact, when we draw an ordinary pencil point on a sheet of paper, the resulting graphite spot always has some finite size.
To more accurately determine the position of a point, we can sharpen a pencil and draw a smaller graphite spot.
If this is not enough, we can arm ourselves with thin needles and microscopes.
The graphite spots will become smaller and smaller, they will represent the point more and more accurately, but they will all have some finite size.
Thus, the idea of the ancient Greeks about the point as an elementary indivisible object implies that this process of successively dividing a graphite spot into smaller parts can continue indefinitely.
The {\em ideal} point itself is thought of as the ultimate result of this infinite process.

\begin{marginfigure}[-3\baselineskip]\Mrg
\includegraphics[width=1\linewidth]{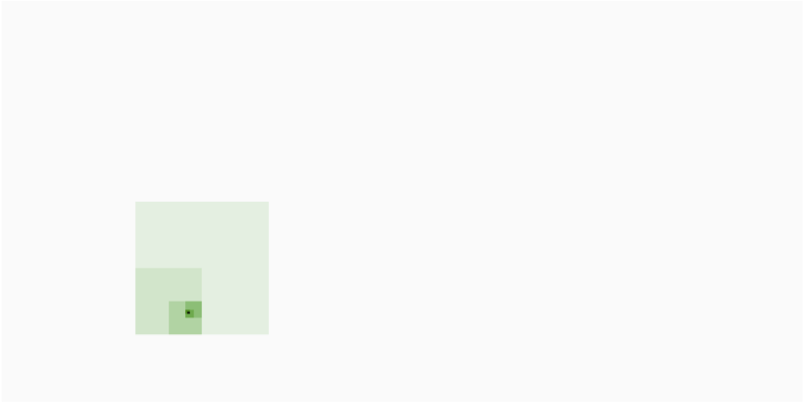}
\caption{Successive approximations of a point on a plane.}
\label{fig:Point}
\end{marginfigure}
Of course, this description of the process of infinitely increasing the observed scale and the successive refinement of the position of a point is not mathematically rigorous.
In the end, we have not defined what the "spot size" and "scale" are.
Any attempt to give these concepts clear definitions inevitably returns us to the concept of a point and to the idea of the distance between two points.

\mathnote{\Mrg
\begin{align}
&\Obj{A\w B\w C\ldots}{\Point}\nonumber\\
&\quad \CC\text{  the symbols $A, B, C, \ldots$ denote points}\nonumber\\[1.5\baselineskip]
&A = B \nonumber\\
&\quad \CC\text{  points $A$ and $B$ coincide with each other}\nonumber\\[0.5\baselineskip]
&A \neq B \nonumber\\
&\quad \CC\text{  points $A$ and $B$ are distinct from each other}\nonumber
\end{align}}
To break the vicious circle of definitions, we have to consider the notion of a point as {\em undefined}. 
Points are usually denoted by uppercase Latin letters: $A, B, C, \ldots$ 
These symbols are understood as a kind of labels of real objects.
Two distinct points cannot have the same label.
At the same time, two different labels may indicate either two distinct points, or one and the same point.
In the latter case, it is said that the points denoted by these different labels {\em coincide} with each other. 
The assertions of coincidence or, conversely, the distinction of points are logical propositions and can be true or false depending on the actual location of points on the plane. The coincidence of objects will be denoted by the equality sign"$=$".

\subsection{1.2. Drawing a point on the plane}

\mathnote{
\begin{align}
&\exists A : \Point \nonumber\\
&\quad \CC\text{  there exists a point $A$ on the plane}\nonumber
\end{align}}
We assume that some point on the plane {\em exist}.
This statement is, however, somewhat speculative.
Having made such an assertion, we have not presented any real figure on the plane. 
We call such weak assertions {\em existential} propositions. 
At the same time, we also assume that we have at our disposal an ideal tool --- a pencil, which allows us in certain situations to draw on the plane real points with the required properties.
Following the ancient Greeks, we call elementary constructions using such ideal tools as {\em postulates}. 
The postulates can also be understood as a kind of {\em commands} for constructing some simple figures with the required properties.

So, our first postulate states that with the help of an ideal pencil one can draw at least one point on the plane.
\vspace{0.5\baselineskip}
\begin{postulate}\label{pst:DrawPoint}\Mrg
\mathnote{\Set
\begin{align}
&\Build{A}{\Point}{\top} \nonumber\\
&\quad \CC \text{\rm  ~draw an arbitrary point $A$ on the plane}\nonumber
\end{align}}
Draw an arbitrary point on the plane.
\end{postulate}

Here we have to comment on the notation that we use.
First of all, in this paper we {\em do not} use the set-theoretic notation.
Thus, curly brackets in our case {\em do not} denote the set of all figures with the required property. 
Instead, we take here the notation of Coq, where the same curly brackets with a vertical line in the middle are understood as simply a pair of two elements. 
The first element of the pair is the constructed geometric figure itself, while the second element is the required property that the constructed figure should have.
Since in our first postulate we do not assume that the constructed point has any intrinsic default properties, we simply write the identically true proposition ($\top$) as the second element of the pair.

\subsection{1.3. Drawing a distinct point on the plane}

\begin{marginfigure}[13.5\baselineskip]\Mrg
\includegraphics[width=\linewidth]{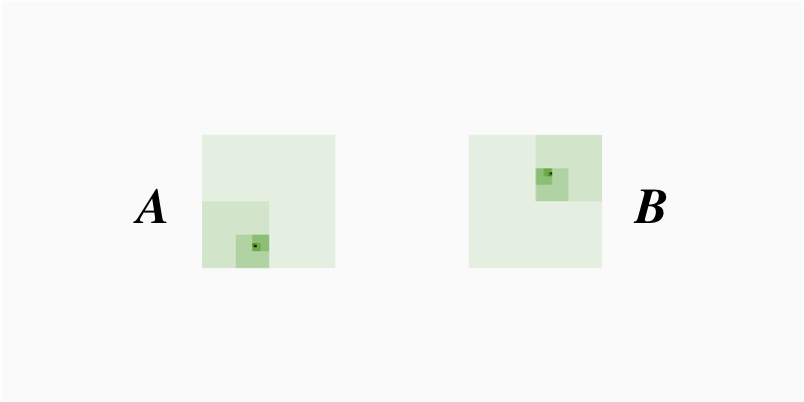}
\caption{Successive approximations of two distinct points on a plane}
\label{fig:TwoPoints}
\end{marginfigure}
In practice, however, there is a fundamental difference between the statements about coincidence and distinction between two points.
Indeed, let us again represent our points in the form of infinite sequences of ever decreasing graphite spots embedded in each other.
It is easy to see that if, on some scale, the graphite spots that represent our two points were distinct and did not touch each other, they would also be distinct at all larger scales (Fig. \ref{fig:TwoPoints}). 
Thus, in order to {\em empirically} verify the statement that the two given points are distinct, we only need to perform a {\em finite} number of steps to increase the scale.

\begin{marginfigure}[0.3\baselineskip]\Mrg
\includegraphics[width=\linewidth]{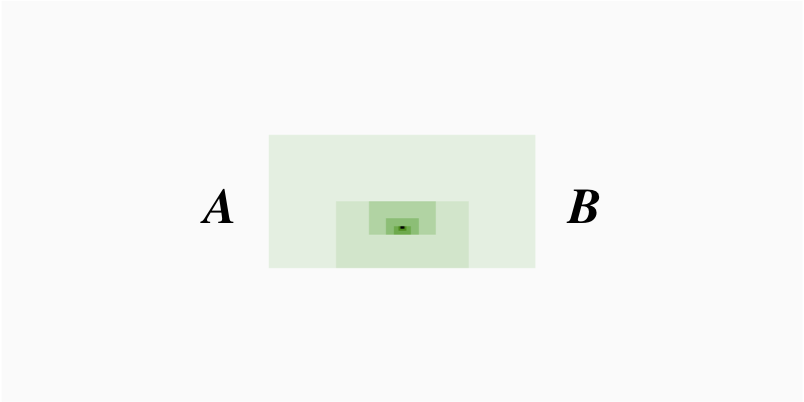}
\caption{Successive approximations of two coinciding points on the plane}
\label{fig:TwoEqualPoints}
\end{marginfigure}

However, is impossible to verify {\em empirically} the fact of coincidence of two arbitrary points. 
Indeed, in this case we need to make sure that the graphite spots that represent our points overlap or at least touch each other on {\it all} larger scales (Fig. \ref{fig:TwoEqualPoints}). 
Obviously, even if we were not able to distinguish from each other two arbitrary points with the strongest available microscope, this does not mean that our points really coincide. 
In the end, in the future we may find an even more powerful microscope, which could allow us to notice the distinction between these points on some larger scale.

Thus, two arbitrary points are distinct if on some {\it finite} scale the graphite spots representing them are clearly distinct and do not touch each other. 
On the contrary, these two points coincide if the graphite spots representing them overlap or touch each other on {\it all} scales. 

Whatever the point on the plane, with a pencil, we can always draw another point, {\it clearly} distinct from the given one. 
The following postulate reflects this our ability.
\vspace{0.5\baselineskip}
\begin{postulate}\Mrg
\mathnote{\Set
\begin{align}
&\Forall{A}{\Point}{}\nonumber\\
&\Build{B}{\Point}{A \neq B}\nonumber\\
&\quad \CC \text{\rm  ~draw point $B$ on the plane that is}\nonumber\\
&\quad \text{\rm  ~~~~~distinct from the point $A$}\nonumber
\end{align}}
Given an arbitrary point on the plane.
Draw another point distinct from the given one.
\label{pst:DrawDistinctPoint}
\end{postulate}

\noindent Here we again use the standard Coq's notations which should not be confused with set-theoretic ones. 
This postulate states that given a point on the plane, we can always construct and present a pair objects represented by curly brackets.
The first element of this pair is a newly constructed point, and its second element provides an additional "certificate" that confirms that this new point is distinct from the given one. In other words, curly brackets in this notation represent a pair --- a point $B$ and a certificate that this points is distinct form the previously constructed point $A$.

\subsection{1.4. How to decide whether points on the plane are distinct?}

Let us recall one of the classical laws of logic --- the law of excluded middle, --- which, as applied to the notion of coincidence of points, can be reformulated as follows:

\vspace{0.5\baselineskip}
\mathnote{\Prop
\begin{align}
&\Forall{A\w B}{\Point}{}\nonumber\\
& A = B\w \lor \s A\neq B \nonumber
\end{align}}
{\noindent{\Prop \bf The law of excluded middle.} \it Two arbitrary points either coincide with each other or are distinct from each other.}
\vspace{0.5\baselineskip}

\noindent This law, of course, gives us a logically true statement, which we can use in our further reasoning. 
But how to put it into practice? 
If we are given a sheet of paper on which two arbitrary points $A$ and $B$ are drawn, can we at any position of these points on the plane unambiguously decide whether they are distinct from each other or do they coincide with each other? 
Of course, we will immediately see the difference if these two points are significantly distant from each other.
If they are close to each other, but still distinct, we can empirically verify this fact, although perhaps we will need to arm ourselves with a magnifying glass or a sufficiently powerful microscope. 
The problem, however, arises when the points $A$ and $B$ actually coincide with each other, but we do not know about it in advance.
In this case, trying to discern the distinction between the points, we will use more and more powerful microscopes, but we will never be completely sure that we can reliably prove the fact that these two points coincide with each other. In the end, there will always be an opportunity to see the distinction between these two points with an even more powerful microscope. In other words, with a certain arrangement of points on the plane and with only a finite set of microscopes, we may not have the practical opportunity to decide which of the two alternatives in the law of the excluded middle is actually true.

This circumstance was noted by Arend Heiting, who proposed another statement, which in practice replaces the law of excluded middle. To formulate it, let us again draw two distinct points $A$ and $B$ on the plane. Then, let us also take an arbitrary third point $C$. This point cannot simultaneously coincide with both distinct points $A$ and $B$. Therefore, it must be distinct from at least one of them.

\begin{marginfigure}[-2\baselineskip]\Mrg
\includegraphics[width=\linewidth]{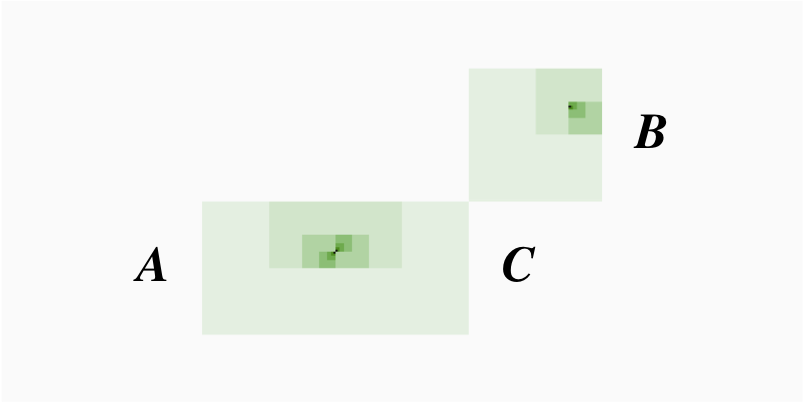}
\caption{Three points on the plane. The point $C$ is distinct from both the point $A$ and the point $B$.}
\label{fig:ThreePoints}
\end{marginfigure}In other words, since the points $A$ and $B$ are distinct, we can choose such a scale at which the graphite spots representing these points are not in contact with each other. Even if on this scale the point $C$ overlaps or touches both the points $A$ and $B$, as shown in the figure \ref {fig:ThreePoints}, then on the next, larger scale, we can certainly discern the distinction between the point $C$ and at least one of the two distinct points $A$ and $B$. Thus, we come to the next postulate.
\vspace{0.5\baselineskip}
\begin{postulate}[Arend Heyting]\label{pst:DecideDistinctPoints}\Mrg
\mathnote{\Set
\begin{align}\
&\Forall{A \w B \w C}{\Point}{}\nonumber\\
& A\neq B \rightarrow \Or{A \neq C}{C \neq B}\nonumber\\
&\quad \CC \text{\rm  ~decide if the point $C$ is distinct}\nonumber\\
&\quad \text{\rm  ~~~~~from the point $A$ or from the point $B$}\nonumber
\end{align}}
Given two distinct points $A$ and $B$ and an arbitrary point  $C$. 
Decide whether the point $C$ is distinct from the point $A$ or from the point $B$.
\end{postulate}
\noindent Here, again, we use standard Coq's notations where two curly brackets joined by plus symbol denote strong constructive disjunction, that is, a practical algorithm that allows us to decide which of the two alternatives is true at any particular position of points and lines on the plane. 

We see that this postulate is significantly different from previous ones. 
At first glance, it is even unclear what kind of ideal tool it describes, because as a result of its use no new objects appear on the plane. 
The answer is simple. 
The Hayting's postulate does not describe properties of any construction tool  (pencil, straightedge or compass), but the observation tool, that is, the "eye" of a person, possibly armed with a magnifying glass or microscope. 
The result of the application of this postulate is not the construction of some new figure on the plane, but the classification of previously constructed figures and the extraction of an additional evidence that these figures possess one or another property.
In other words, this postulate allows us to empirically extract an additional "certificate" of distinction of the third point from one of the first two.

\subsection{1.5. Propositions and algorithms}

Let us summarize our discussion of Euclid's  constructive-deductive logic in the following table, where we contrast weak logical propositions about the speculative existence of some objects with strong algorithms for their construction. 
Constructing a figure on a plane always entails the truth of a purely logical proposition about the existence of this figure, but the reverse is not true.
Even if we have proven the truth of the statement about the existence of a figure with a given property, we, nevertheless, may not have an algorithm for constructing it using our limited set of geometric tools.
\vspace{0.5\baselineskip}
\begin{table*}
\begin{tabular*}{\textwidth}{p{0.475\textwidth}p{0.475\textwidth}}
\text{\sc\Prop DEDUCTIVE METHOD} \em 
&\text{\sc \Set CONSTRUCTIVE METHOD} \em 
\\
\toprule\\[-4mm]
\text{\bf \Prop Existential proposition.} \em 
Given figure $A$ with property ${\mathcal P}_A$. There exists figure $B$ with the property ${\mathcal Q}_B$.
&\text{\bf \Set Construction algorithm.} \em 
Given figure $A$ with property ${\mathcal P}_A$. Draw figure $B$ with the property ${\mathcal Q}_B$.
\\[2mm]
\begin{tcolorbox}
[colback=MyBlue!5!white, frame hidden, enhanced, attach boxed title to top right={yshift=-5.2mm,xshift=-0.2mm}, 
boxed title style={size=small, frame hidden, colback=MyBlue},
title=\footnotesize\Axiom]
\footnotesize $\Forall{A}{\Figure}{{\mathcal P}_A \rightarrow\Exists{B}{\Figure}{{\mathcal Q}_B}}$ 
\end{tcolorbox}& 
\begin{tcolorbox}
[colback=MyGreen!5!white, frame hidden, enhanced, attach boxed title to top right={yshift=-5.2mm,xshift=-0.2mm}, 
boxed title style={size=small, frame hidden, colback=MyGreen},
title=\footnotesize\Postulate]
\footnotesize $\Forall{A}{\Figure}{{\mathcal P}_A \rightarrow\Build{B}{\Figure}{{\mathcal Q}_B}}$ 
\end{tcolorbox}
\\[-5mm]
\midrule
\text{\bf \Prop  Disjunctive proposition.} \em 
If figure $A$ has property ${\mathcal P}_A$, then it also has either property ${\mathcal L}_A$ or property ${\mathcal R}_A$.
&\text{\bf \Set Decision algorithm.} \em 
Given figure $A$ with property ${\mathcal P}_A$.
Decide which alternative is true  ${\mathcal L}_A$ or ${\mathcal R}_A$.
\\[2mm]
\begin{tcolorbox}
[colback=MyBlue!5!white, frame hidden, enhanced, attach boxed title to top right={yshift=-5.2mm,xshift=-0.2mm}, 
boxed title style={size=small, frame hidden, colback=MyBlue},
title=\footnotesize\Axiom]
\footnotesize $\Forall{A}{\Figure}{{\mathcal P}_A \rightarrow {{\mathcal L}_A}\lor{{\mathcal R}_A}}$ 
\end{tcolorbox}& 
\begin{tcolorbox}
[colback=MyGreen!5!white, frame hidden, enhanced, attach boxed title to top right={yshift=-5.2mm,xshift=-0.2mm}, 
boxed title style={size=small, frame hidden, colback=MyGreen},
title=\footnotesize\Postulate]
\footnotesize $\Forall{A}{\Figure}{{\mathcal P}_A \rightarrow \Or{{\mathcal L}_A}{{\mathcal R}_A}}$ 
\end{tcolorbox}
\\[-5mm]
\bottomrule
\end{tabular*}
  \caption{Existential and disjunctive propositions vs construction and decision algorithms.}
  \label{tab:DD}
\end{table*}

In the same way, we contrast weak logical disjunctive propositions with stronger decision algorithms, which in practice can distinguish one arrangement of figures on a plane from another.
The ability in practice to decide which of the two classes to assign a particular arrangement of figures on a plane entails the truth of corresponding logical disjunction, but the reverse is not true. The truth of disjunctive proposition does not give us a practical algorithm for the classification of figures on the plane.

\section{2. Lines}

\subsection{2.1. Conceptualization of line}

\mathnote{\Mrg
\begin{align}
&a\;b\;c\ldots : \Line\nonumber\\
&\quad\CC\text{  the symbols $a, b, c, \ldots$ denote straight lines}\nonumber\\[1.5\baselineskip]
&a = b \nonumber\\
&\quad \CC\text{  the lines $a$ and $b$ coincide with each other}\nonumber\\[0.5\baselineskip]
&a \neq b \nonumber\\
&\quad \CC\text{  the lines $a$ and $b$ are distinct from each other}\nonumber
\end{align}}
The ancient Greeks described a straight line as such a geometric figure, that has a length, but has no width.
Nowadays, we can no longer be satisfied with such a description, since the concepts of length and width used in it themselves need to be defined. 
Therefore, we consider the concept of a straight line, like that of a point, as undefined. 
The rules for constructing straight lines will be set by postulates, and their properties will be disclosed through axioms.

We shall denote straight lines by lowercase Latin letters: $a, b, c, \ldots$
As in the case of points, two different names can indicate both different straight lines and the same straight line.
In the latter case, it is said that the straight lines pointed to by these different names {\em coincide} with each other.
Propositions of coincidence or, conversely, distinction of two straight lines can be true or false, depending on their actual location on the plane.

\subsection{2.2. Incidence relation}

\begin{marginfigure}[0\baselineskip]\Mrg
\includegraphics[width=\linewidth]{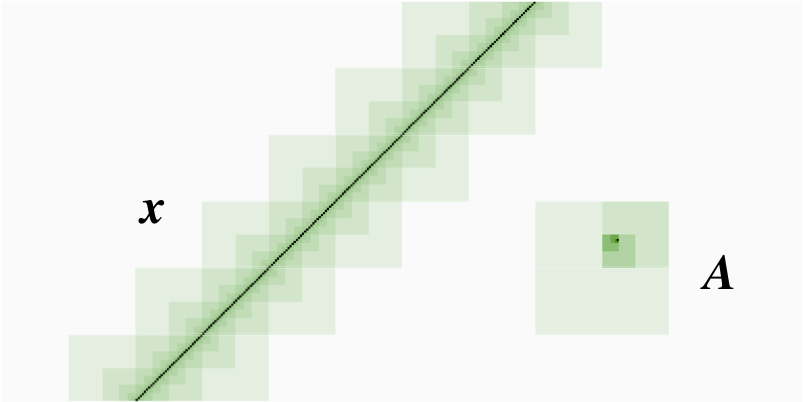}
\caption{The point $A$ does not belong to the line $x$.}
\label{fig:PointAndLine}
\end{marginfigure}

\mathnote{\Mrg
\begin{align}
&\ObjDouble{A\w B\ldots}{\Point}{x \w y\ldots}{\Line}{}\nonumber\\[0.5\baselineskip]
&\In{A\w B\ldots}{x}\equiv A\in x \land B\in x \land\ldots\nonumber\\
&\quad \CC\text{  points $A\cs B\cs \ldots$ lie on the line $x$}\nonumber\\[0.5\baselineskip]
&\In{A}{x\w y\ldots}\equiv A\in x \land A\in y \land\ldots\nonumber\\
&\quad \CC\text{  the lines $x\cs y\cs \ldots$ pass through point $A$}\nonumber
\end{align}}
The main undefined relation, that determines the relative position of points and straight lines on the plane, is the relation of {\em incidence}.
Recall that in real life, when plotting drawings on a sheet of paper, both points and lines are displayed only as graphite spots (Fig. \ref {fig:PointAndLine}). 
It is easy to see that if at some scale the spots that represent a point and a straight line do not intersect and do not touch each other, then they will be distinguishable on all larger scales. In this case, we say that the point is {\em apart} from the line. It is impossible to confirm {\em empirically} the fact that a point belongs to a straight line, because for this we would have to make sure that the graphite spots that represent this point and this straight line overlap or touch each other on all scales.

Note that expressions often used in geometry: "a point lies on a line", "a point belongs to a line", "a line passes through a point" --- all of them are equivalent to the proposition "a point is incident to a straight line".

\subsection{2.3. Drawing distinct points on a line and a point apart it}

We have already seen that a straight line is a geometric figure, which consists not of one, but of many points.
Later we shall prove that on a line one can draw an infinite number of points.
The basis of all these future constructions will be the following two postulates.
\vskip0.5\baselineskip
\begin{postulate}\label{pst:DrawPointOnLine}\Mrg
\mathnote{\Set
\begin{align}
&\Forall{x}{\Line}{}\nonumber\\
&\Build{A}{\Point}{A \in x}\nonumber
\end{align}}
Given an arbitrary straight line. Draw a point on the line.
\end{postulate}

\begin{postulate}\label{pst:DrawDistinctPointOnLine}\Mrg
\mathnote{\Set
\begin{align}
&\ForallDouble{A}{\Point}{x}{\Line}{}\nonumber\\
& A \in x \rightarrow \Build{B}{\Point}{A \neq B \land B \in x}\nonumber
\end{align}}
Given an arbitrary straight line and a point on it. Draw another point on this line distinct from the given one.
\end{postulate}
\begin{fullwidth}\vspace{-1\baselineskip}\end{fullwidth}

We also know that no straight line can coincide with the entire plane. 
In other words, whatever the line, there are points on the plane that do not belong to it. 
The following postulate allows us to draw one of these points.
\vskip0.5\baselineskip
\begin{postulate}\label{pst:DrawPointApartLine}\Mrg
\mathnote{\Set
\begin{align}
&\Forall{x}{\Line}{}\nonumber\\
&\Build{A}{\Point}{A \notin x}\nonumber
\end{align}}
Given an arbitrary line.
Draw a point that does not lie on this line.
\end{postulate}

So, no matter what the line is, you can draw at least two distinct points that are incident to this line, and one point apart from it.

\subsection{2.4. Drawing a straight line through two distinct points}

With the help of Postulates \ref{pst:DrawPoint} and \ref{pst:DrawDistinctPoint}, we can draw on the plane two distinct points. From our everyday experience, we know that through the points, using a pencil and a straightedge, we can draw a straight line. We formulate this our ability in the form of the next postulate.
\vskip0.5\baselineskip
\begin{postulate}[Point-Point extension]\label{pst:DrawExtensionLine}\Mrg
\mathnote{\Set
\begin{align}
&\Forall{A \w B}{\Point}{}\nonumber\\
& A \neq B\rightarrow \Build{x}{\Line}{A\w B \in x} \nonumber
\end{align}}
Given two distinct points on the plane.
Draw a straight line through both these points.
\end{postulate}
It is interesting to note that this postulate gives us {\em the only} way to draw a straight line on the plane --- to draw it through two distinct points. 
Therefore, we may assume that the Postulates \ref{pst:DrawPointOnLine} and \ref{pst:DrawDistinctPointOnLine} do not actually draw any new points on the line, but only restore those points through which the straight line was originally drawn.

As we have said, the postulates describe some elementary constructions that we can perform on the plane with the help of our ideal tools. 
Applying the postulates one after another one, we can draw on the plane new, more complex figures with the required properties.
Before applying each postulate, we must make sure that all the necessary prerequisites for its applicability are fulfilled.
With the help of the already introduced postulates, we can solve our first simple problems.
\vskip0.5\baselineskip
\begin{problem}\Mrg
\mathnote{\Set
\begin{align}
&\Forall{A}{\Point}{}\nonumber\\
&\Build{x}{\Line}{A\in x}\nonumber
\end{align}}Given an arbitrary point. Draw a line through this point.
\begin{proof}[\Set\DEQ]\Mrg
\proofnote{
\begin{align}
\Deq&\s[5]\Build{A}{\Point}{\top} \tag*{Given} \nonumber\\
\Ref{pbm:DrawLineThroughPoint_2}&\s[5]\Build{B}{\Point}{A \neq B} \tag*{\Postulate{\ref{pst:DrawDistinctPoint}}} \nonumber\\
\Ref{pbm:DrawLineThroughPoint_3}&\s[5]\Build{x}{\Line}{A\w B \in x} \tag*{\Postulate{\ref{pst:DrawExtensionLine}}}\nonumber\\
&\mquad \Build{x}{\Line}{A\in x}\tag*{\Qed} 
\end{align}}Given an arbitrary point $A$.
\begin{enumerate}
\item Draw a point $B$ distinct from the given point $A$ (Postulate \ref{pst:DrawDistinctPoint}). \label{pbm:DrawLineThroughPoint_2}
\item Draw a line $x$ through distinct points $A$ and $B$ (Postulate \ref{pst:DrawExtensionLine}). \label{pbm:DrawLineThroughPoint_3}
\end{enumerate}\noindent
Finally, we can erase the auxiliary point $B$ and leave on the plane only the required line $x$ that passes through the given point $A$.\Set\QED
\end{proof}
\label{pbm:DrawLineThroughPoint}
\end{problem}
\noindent On the margins we briefly outline the solution to the problem, i.e. our construction.
There, for brevity, we refer to Postulates, denoting them by the symbols \Postulate{n}, where $n$ is the number of the corresponding Postulate. 
The solutions of the previously solved problems are denoted by \Solution{m}, where $m$ is the number of the corresponding problem. 
So, if we ever need to draw a line passing through a given point, we may not repeat all the previous construction steps and arguments, and  immediately address the solution of the Problem \ref{pbm:DrawLineThroughPoint}.
\begin{problem}\Mrg\label{pbm:DrawLine}
\mathnote{\Set
\begin{align}
&\Build{x}{\Line}{\top}\nonumber
\end{align}} Draw an arbitrary straight line on the plane.
\begin{proof}[\Set\DEQ]\Mrg
\proofnote{
\begin{align}
\Deq&\s[5] \tag*{Given} \nonumber\\
\Ref{pbm:DrawLine_1}&\s[5]\Build{A}{\Point}{\top} \tag*{\Postulate{\ref{pst:DrawPoint}}} \nonumber\\
\Ref{pbm:DrawLine_2}&\s[5]\Build{x}{\Line}{A\in x} \tag*{\Solution{\ref{pbm:DrawLineThroughPoint}}}\nonumber\\
&\mquad\Build{x}{\Line}{\top}\tag*{\Qed} 
\end{align}}Given empty plane.
\begin{enumerate}
\item Draw an arbitrary point $A$ (Postulate \ref{pst:DrawPoint}). \label{pbm:DrawLine_1}
\item Draw a line $x$ passing through the point $A$ (Problem \ref{pbm:DrawLineThroughPoint}). \label{pbm:DrawLine_2}
\end{enumerate}\noindent
Finally, we can erase the auxiliary point $A$ and leave on the plane only the required line $x$.\Set\QED
\end{proof}
\end{problem}

\subsection{2.5. Uniqueness of the straight line through two distinct points}

\begin{marginfigure}[-5\baselineskip]\Mrg
\includegraphics[width=\linewidth]{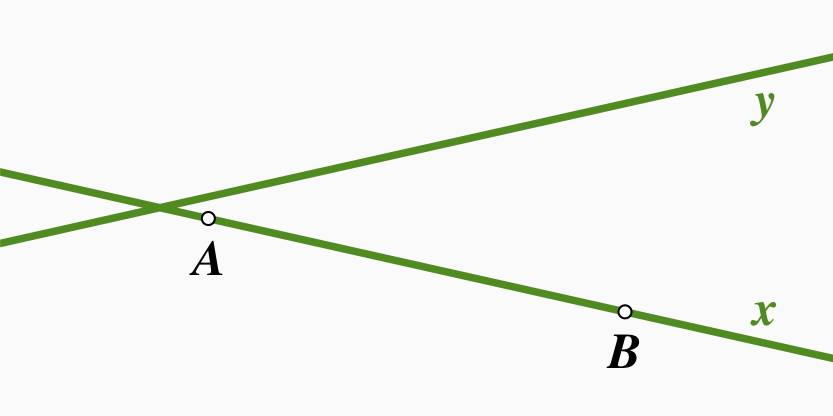}
\caption{Two distinct lines $x$ and $y$ intersect at two distinct points $A$ and $B$.} 
\label{fig:vonPlato}
\end{marginfigure}

Let us draw on the plane two distinct points $A$ and $B$ (Postulates \ref{pst:DrawPoint} and \ref{pst:DrawDistinctPoint}) and draw a line $x$ through these points (Postulates \ref{pst:DrawExtensionLine}). 
Consider on the plane some straight line $y$, distinct from the just constructed line $x$ (see Fig. \ref {fig:vonPlato}). 
We know from everyday experience that one can draw {\em only one} line through two distinct points.
This means that the points $A$ and $B$ cannot simultaneously belong to the line $y$. 
Indeed, in this case, two distinct lines $x$ and $y$ would pass through two distinct points $A$ and $B$. 
In other words, one of these points or both of them cannot not lie on the line $y$. 
The following postulate states that for any arrangement of points and lines on a plane we can always decide which of these two points, $A$ or $B$, {\em obviously} does not belong to $y$.

The following postulate was proposed by Jan von Plato and is similar to Heyting's postulate in the sense that it describes the observer's ability to distinguish between two possible positions of points and lines on the plane.
\vskip0.5\baselineskip
\begin{postulate}[Jan Von Plato]\label{pst:DecideUniqueExtensionLine}\Mrg
\mathnote{\Set
\begin{align}
&\ForallDouble{A \w B}{\Point}{x \w y}{\Line}{}\nonumber\\
& A \neq B \rightarrow \In{A\w B}{x} \rightarrow  x \neq y\rightarrow\Or{A \notin y}{B \notin y}\nonumber
\end{align}}
Given two distinct lines, on one of which two distinct points are marked. 
Decide which of these two points does not belong the other line.
\end{postulate}
\noindent This postulate directly implies that only one straight line can be drawn through two distinct points, and two distinct straight lines can have only one point in common.
\vskip0.5\baselineskip
\begin{theorem}\Mrg
a)\mathnote{\Prop
\begin{align}
&\ForallDouble{A \w B}{\Point}{x \w y}{\Line}{}\nonumber\\
& A\neq B\rightarrow \In{A \w B}{x \w y} \rightarrow x = y\nonumber
\end{align}} If two lines go through two distinct points, these lines coincide.\label{thm:DiPs_EqLs}\\
b)\mathnote{\Prop
\begin{align}
& x\neq y\rightarrow \In{A \w B}{x \w y} \rightarrow A = B\nonumber
\end{align}} If two distinct lines have two common points, these points coincide.
\label{thm:DeLs_EqPs}\\
 $~$\hfill{\Prop\rm \scriptsize / Hilbert, Chapter 1 : Theorem 1}
\end{theorem}

\vskip0.5\baselineskip
\begin{problem}\Mrg\label{pbm:DrawLineNotThroughPoint}
\mathnote{\Set
\begin{align}
&\Forall{A}{\Point}{}\nonumber\\
&\Build{x}{\Line}{A \notin x}\nonumber
\end{align}}
Given an arbitrary point. Draw a line that does not pass through this point.
\begin{proof}[\Set\DEQ]{\Mrg
\proofnote{\Mrg
\begin{align}
\Deq&\s[5]\Build{A}{\Point}{\top} \tag*{Given} \nonumber\\
\Ref{pbm:DrawLineNotThroughPoint_1}&\s[5]\Build{B}{\Point}{A \neq B} \tag*{\Postulate{\ref{pst:DrawDistinctPoint}}} \nonumber\\
\Ref{pbm:DrawLineNotThroughPoint_2}&\s[5]\Build{y}{\Line}{A\w B \in y} \tag*{\Postulate{\ref{pst:DrawExtensionLine}}} \nonumber\\
\Ref{pbm:DrawLineNotThroughPoint_3}&\s[5]\Build{C}{\Point}{C \notin y} \tag*{\Postulate{\ref{pst:DrawPointApartLine}}} \nonumber\\
\Ref{pbm:DrawLineNotThroughPoint_5}&\s[5]\Build{x}{\Line}{B \w C \in x} \tag*{\Postulate{\ref{pst:DrawExtensionLine}}}\nonumber\\
&\mquad\Build{x}{\Line}{A \notin x}\tag*{\Qed} 
\end{align}}Let us be given some point $ A $. Then:
\begin{enumerate}
\item Draw a point $B$ distinct from the point $A$  (Postulate \ref{pst:DrawDistinctPoint}). \label{pbm:DrawLineNotThroughPoint_1}
\item Draw a line $y$ through distinct points $A$ and $B$ (Postulate \ref{pst:DrawExtensionLine}). \label{pbm:DrawLineNotThroughPoint_2}
\item Draw a point $C$ not lying on the line $y$ (Postulate \ref{pst:DrawPointApartLine}). \label{pbm:DrawLineNotThroughPoint_3}
\item The points $B$ and $C$ are distinct, since line $y$ passes through $B$ and does not pass through $C$. \label{pbm:DrawLineNotThroughPoint_4} 
\item Draw a line $x$ through distinct points $B$ and $C$ (Postulate \ref{pst:DrawExtensionLine}). \label{pbm:DrawLineNotThroughPoint_5}
\item The lines $x$ and $y$ are distinct, since the the point $C$ lies on $x$ and does not lie on the $y$. \label{pbm:DrawLineNotThroughPoint_6} 
\item Point $A$ does not lie on line $x$. Otherwise, through two distinct points, $A$ and $B$, would pass two distinct lines, $x$ and $y$, which is impossible (Theorem \ref{thm:DiPs_EqLs}). \label{pbm:DrawLineNotThroughPoint_6} 
\end{enumerate}
\noindent So, we have built a line $x$, which does not pass through a given point $A$.}
\Set\QED
\end{proof}
\begin{marginfigure}[-3\baselineskip]\Mrg
\includegraphics[width=\linewidth]{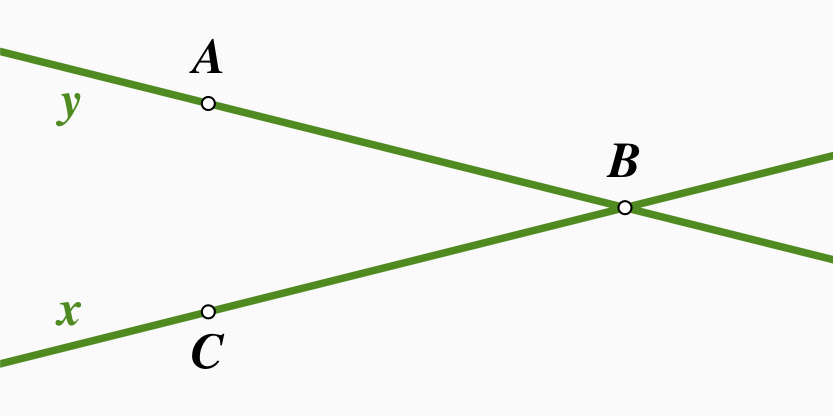}
\caption{Construction of a line $ x $ not passing through this point $A$.}
\label{fig:DrawLineNotThroughPoint}
\end{marginfigure}
\end{problem}

\vskip0.5\baselineskip
We see that the solution to this problem consists of successive steps of two types. 
On the steps of the first type, we, as before, are constructing new geometric figures on the plane with the help of postulates and problems that we have solved earlier. 
We summarize these construction steps on the margins. 
In fact, they define a program for constructing a figure with the required property that can be performed by any graphic robot. 
In this example, the construction steps are steps 1-3) and 5).

On the remaining steps, we prove prerequisite conditions that define the applicability of the postulates and previously solved problems. 
For example, to draw a line through two points on step 5), we must first ensure that these two points are distinct on step 4). 
In the last steps 6-7) we prove that the final figure does have the required property. 
Namely, we prove that the line that we have finally constructed does not pass through the initially given point.

\subsection{2.6. Drawing a points of intersection of two straight lines}

If two straight lines are distinct and somewhere on the plane {\em there exists} a point belonging to both of these lines, then we will say that these two lines {\em intersect}.
\vskip0.5\baselineskip
\begin{definition}
\label{def:IntersectingLines}\Mrg
\mathnote{\Def
\begin{align}
&\Obj{x\w y}{\Line},\nonumber\\
& x \nparallel y \equiv x \neq y \land \Exists{A}{\Point}{A \in x \w y}\nonumber\\
&\quad\CC\text{\rm  the lines $x$ and $y$ intersect each other}\nonumber
\end{align}}
It is said that two distinct lines intersect iff both of them pass through the same point of the plane.
\end{definition}
\noindent Theorem \ref{thm:DeLs_EqPs} ensures that the common point of the two intersecting lines is unique. 
The following postulate gives us the possibility to draw this point.
\vskip0.5\baselineskip
\begin{postulate}[Line-Line intersection]\label{pst:DrawIntersectionPoint}\Mrg
\mathnote{\Set
\begin{align}
&\Forall{x\w y}{\Line}{}\nonumber\\
&x\nparallel y\rightarrow\Build{A}{\Point}{A \in x \w y} \nonumber
\end{align}}
Given two intersecting lines. 
Draw the point of their intersection.
\end{postulate}
\noindent This postulate says that if the intersection point of two lines exists, then it can be constructed. 
This may cause some difficulties.
Probably the best way to understand this is to consider a model in which the points are represented by their coordinates, and the lines by the corresponding linear equations. In this model, the statement that the two given lines intersect is translated into the statement that the difference between the slope coefficients of these two lines is nonzero. If we take the decimal expansion of this difference, then we can perform a linear search that is guaranteed to find the nonzero digit in the expansion and terminate.
Having the leading digits of the expansion, we are able to calculate the coordinates of the intersection point.
\vskip0.5\baselineskip
\begin{problem}\Set
\mathnote{
\begin{align}
&\ForallDouble{A}{\Point}{x}{\Line}{}\nonumber\\
&A \in x \rightarrow \Build{y}{\Line}{x \nparallel y \land A \in y}\nonumber
\end{align}}Given an arbitrary point and straight line passing through it. 
Draw another straight line, distinct from the given one and passing through the same point.
\begin{proof}[\Set\DEQ]{\Mrg
\proofnote{
\begin{align}
\Deq&\s[5]\BuildDouble{A}{\Point}{x}{\Line}{A \in x} \tag*{Given}\nonumber\\
\Ref{pbm:DrawDistinctLineThroughPoint_1}&\s[5]\{B : \Point\vl B \notin x\} \tag*{\Postulate{\ref{pst:DrawPointApartLine}}}\nonumber\\
\Ref{pbm:DrawDistinctLineThroughPoint_3}&\s[5]\{y : \Line\vl A\w B \in y \} \tag*{\Postulate{\ref{pst:DrawExtensionLine}}} \nonumber\\
&\mquad\Build{y}{\Line}{x \neq y \land A \in x \w y}\tag*{\Qed} 
\end{align}}Let us be given a point $A$ and a line $x$ passing through it.
\begin{enumerate}
\item Draw a point $B$, that does not lie on $x$ (Postulate \ref{pst:DrawPointApartLine}). \label{pbm:DrawDistinctLineThroughPoint_1}
\item Points $A$ and $B$ are distinct, since $A$ lies on $x$, and $B$ does not. \label{pbm:DrawDistinctLineThroughPoint_2}
\item Draw a line $y$ through distinct points $A$ and $B$  (Postulate \ref{pst:DrawExtensionLine}). \label{pbm:DrawDistinctLineThroughPoint_3}
\item Lines $x$ and $y$ are distinct, since $B$ lies on $y$ and does not lie on $x$. \label{pbm:DrawDistinctLineThroughPoint_4}
\end{enumerate}\noindent
So, we have constructed a line $y$, which is distinct from the line $x$ and passes through a given point $A$ on this line.
The Definition \ref{def:IntersectingLines} also implies that lines $x$ and $y$ intersect with each other.}
\Set\QED
\end{proof}
\label{pbm:DrawDistinctLineThroughPoint}
\end{problem}
\vskip0.5\baselineskip
\begin{problem}\Set
\mathnote{
\begin{align}
&\Forall{x}{\Line}{}\nonumber\\
&\Build{y}{\Line}{x \nparallel y}\nonumber
\end{align}}Given an arbitrary straight line. Draw another line that intersect it.
\begin{proof}[\Set\DEQ]{\Mrg
\proofnote{
\begin{align}
\Deq&\s[5]\Build{x}{\Line}{\top} \tag*{Given}\nonumber\\
\Ref{pbm:DrawDistinctLine_1}&\s[5]\Build{A}{\Point}{A \in x} \tag*{\Postulate{\ref{pst:DrawPointOnLine}}}\nonumber\\
\Ref{pbm:DrawDistinctLine_2}&\s[5]\Build{y}{\Line}{x \nparallel y \land A \in y} \tag*{\Solution{\ref{pbm:DrawDistinctLineThroughPoint}}} \nonumber\\
&\mquad\Build{y}{\Line}{x \nparallel y}\tag*{\Qed} 
\end{align}}Let us be given some straight line $x$.
\begin{enumerate}
\item Draw a point $A$, lying on the line $x$ (Postulate \ref{pst:DrawPointOnLine}). \label{pbm:DrawDistinctLine_1}
\item Draw a line $y$, distinct from line $x$ and passing through the point $A$  (Problem \ref{pbm:DrawDistinctLineThroughPoint}). \label{pbm:DrawDistinctLine_2}
\end{enumerate}\noindent
Finally, we can erase the auxiliary point $A$ and leave on the plane only the required line $y$, which intersects with the given line $x$.}\Set\QED
\end{proof}
\label{pbm:DrawDistinctLine}
\end{problem}

\vskip0.5\baselineskip
\begin{problem}\Mrg
\mathnote{\Set
\begin{align}
&\ForallDouble{A \w B}{\Point}{x \w y}{\Line}{}\nonumber\\
&A \neq B \rightarrow \In{A}{x \w y} \rightarrow x \nparallel y\rightarrow \Or{B \notin x}{B \notin y}\nonumber
\end{align}}
Given two distinct points, through one of which pass two distinct lines.
 Decide which of these two lines does not pass through another point.
\begin{proof}[\Set\DEQ]{\Mrg
\proofnote{\Mrg
\begin{align}
\Deq&\s[5]\{\s A\w B \w x\w y \vl A \neq B\land A \in x\w y\land x \neq y \s\} \tag*{Given}\nonumber\\
\Ref{pbm:DecidePointApartTwoIntersectingLines_1}&\s[5]\Build{z}{\Line}{A \w B \in z} \tag*{\Postulate{\ref{pst:DrawExtensionLine}}}  \nonumber\\
\Ref{pbm:DecidePointApartTwoIntersectingLines_2}&\s[5]\Build{C}{\Point}{C \in x \land A \neq C} \tag*{\Postulate{\ref{pst:DrawDistinctPointOnLine}}} \nonumber\\
\Ref{pbm:DecidePointApartTwoIntersectingLines_3}&\s[5]\Build{D}{\Point}{D \in y \land A \neq D} \tag*{\Postulate{\ref{pst:DrawDistinctPointOnLine}}} \nonumber\\
\Ref{pbm:DecidePointApartTwoIntersectingLines_5}&\s[5]\Build{t}{\Line}{C \w D \in t} \tag*{\Postulate{\ref{pst:DrawExtensionLine}}}  \nonumber\\
\Ref{pbm:DecidePointApartTwoIntersectingLines_8}&\s[5]\Or{C \notin z}{D \notin z}\tag*{\Postulate{\ref{pst:DecideUniqueExtensionLine}}}\nonumber\\
&\Item{a} C \notin z \rightarrow B \notin x\nonumber\\
&\Item{b} D \notin z \rightarrow B \notin y\nonumber\\
&\mquad \Or{B \notin x}{B \notin y}\tag*{\Qed} 
\end{align}
}Let us be given two distinct points $A$ and $B$ and two distinct lines $x$ and $y$ passing through $A$.
\begin{enumerate}
\item Draw a line $z$ through distinct points $A$ and $B$ (Postulate \ref{pst:DrawExtensionLine}). 
\label{pbm:DecidePointApartTwoIntersectingLines_1}
\item Draw a point $C$ on the line $x$, distinct from $A$ (Postulate \ref{pst:DrawDistinctPointOnLine}). \label{pbm:DecidePointApartTwoIntersectingLines_2}
\item Draw a point $D$ on the line $y$, distinct from $A$ (Postulate \ref{pst:DrawDistinctPointOnLine}). \label{pbm:DecidePointApartTwoIntersectingLines_3}
\item The points $C$ and $D$ are distinct, since otherwise the lines $x$ and $y$ would coincide. (Theorem \ref{thm:DiPs_EqLs}). \label{pbm:DecidePointApartTwoIntersectingLines_4}
\item Draw a line $t$ through distinct points $C$ and $D$ (Postulate \ref{pst:DrawExtensionLine}). 
\label{pbm:DecidePointApartTwoIntersectingLines_5}
\item The point $A$ does not lie on the line $t$, because otherwise the lines $x$ and $y$ would coincide (Theorem \ref{thm:DiPs_EqLs}). \label{pbm:DecidePointApartTwoIntersectingLines_6}
\item The lines $z$ and $t$ are distinct, since the line $z$ passes through the point $A$, and the line $y$ does not (Theorem \ref{thm:DiPs_EqLs}). \label{pbm:DecidePointApartTwoIntersectingLines_7}
\item Decide which of these two points, $C$ or $D$, does not belong to the line $z$ (Postulate \ref{pst:DecideUniqueExtensionLine}). \label{pbm:DecidePointApartTwoIntersectingLines_8}
\begin{enumerate}
\item [a)] If the point $C$ does not lie on the line $z$, then the point $B$ cannot lie on the line $x$.
\item [b)] If the point $D$ does not lie on the line $z$, then the point $B$ cannot lie on the line $y$.
\end{enumerate}
\end{enumerate}\noindent
So, we have decided which of the two lines, $x$ or $y$, does not pass through the second point $B$.}\Set\qedhere
\end{proof}
\label{pol:DecidePointApartTwoIntersectingLines}
\end{problem}


\subsection{2.7. Collinearity of points}

Based on the notion of incidence of points and lines, we can introduce a new notion of {\it collinearity} of points.
\vskip0.5\baselineskip
\begin{definition}\label{def:Collinearity}\Mrg
\mathnote{\Def
\begin{align}
&\Obj{A\w B\w C \ldots}{\Point},\nonumber\\
&\ColThree{A}{B}{C \w \ldots} \equiv \Exists{x }{\Line}{A\w B\w C\w \ldots \in x} \nonumber\\
&\quad\CC\text{\rm  points $A\cs B\cs C\cs \ldots$ lie on the same line}\nonumber
\end{align}}
Points are called collinear if all of them lie on one straight line.
\end{definition}
\noindent Note that according to this definition, the collinearity property of points does not depend on the order in which the points are written.

This definition states that a straight line passing through collinear points {\em exists}. This, however, does not mean that such a straight line can be {\em constructed}. 
We shall illustrate this idea on the following example.
\vskip0.5\baselineskip
\begin{theorem}\label{cor:ExLineThroughTwoPoints}\Mrg
\mathnote{\Prop
\begin{align}
&\Forall{A \w B}{\Point}{}\nonumber\\
&\ColThree{A}{B}{}{\s[-1]}\nonumber
\end{align}}
Whatever the two points, there is a straight line passing through the points.
\begin{proof}[\Prop\DEQ]\Mrg
{According to the law of excluded middle, the two points $A$ and $B$ either coincide or are distinct from each other.
\begin{enumerate}
\item[a)]In the first case, when $A = B$, we can draw a straight line through this only point using construction from the Problem \ref{pbm:DrawLineThroughPoint}. 
\item[b)]In the second case, when $A \neq B$, we can draw a line through these two distinct points with the help of the Postulate \ref{pst:DrawExtensionLine}.
{\Prop\QED}
\end{enumerate}}
\end{proof}
\end{theorem}
\vskip0.5\baselineskip
So, we have proved that for any position of two points on the plane there exists some straight line passing through them.
However, we can not offer any algorithm to draw this line.
Indeed, as we discussed earlier, if we are given two arbitrary points on a plane, we cannot decide whether they are the same or are distinct from each other, and therefore a situation may arise when we do not know which tool to use: the solution of the Problem \ref{pbm:DrawLineThroughPoint} or the Postulate \ref{pst:DrawExtensionLine}.
Thus, even having proved the existence of a geometrical figure with the required properties, we, nevertheless, may not have an algorithm for its construction using our toolkit. 

At the same time, the mere fact of the existence of a geometric figure is often enough to analyze its properties. 
As an example, consider the following important property of non-collinear points which follows immediately from the above theorem.
\vskip0.5\baselineskip
\begin{theorem}\label{stm:nColPs_DiPs}\Mrg
\mathnote{\Prop
\begin{align}
&\Forall{A \w B \w C}{\Point}{}\nonumber\\
&\lnot \ColThree{A}{B}{C} \rightarrow A \neq B \land B \neq C \land A \neq C \nonumber
\end{align}} 
If three points do not lie on the same line, then they are all distinct.
\end{theorem}

This simple example illustrates why we should not easily give up the law of excluded middle, as intuitionists do. 
It can still be very useful for studying the properties of the constructed figures.

\subsection{2.8. Duality of points and lines in the incidence geometry}

Now we are going to show that in the incidence geometry, consisting only of nine postulates introduced so far, there is a subtheory with a deep internal symmetry.
To illustrate this duality, we shall write down the postulates together with some previously solved problems in the following Table. \ref{tab:duality}.
Here, each postulate in the left column of the table corresponds to exactly one postulate or problem in the right column, in which points and straight lines are swapped, and the distinction between two points is interchanged with the intersection relation between two straight lines.

We have placed in the table a slightly weaker version of the postulate of Jan von Plato.
Thus, the duality actually takes place not for the whole incidence geometry, but only for its certain subtheory.
However, it is well known that such duality plays a fundamental role in the projective geometry.

An attentive reader will notice that there is one problem not proved up to now. 
Namely, we have not yet proved Problem \ref{pbm:DecideIntersectingLines}, which is the dual analogue of Heyting's postulate. We will be able to prove it later, using also our other postulates and axioms.

\vskip0.5\baselineskip
\begin{table*}[p]
\begin{tabular*}{\textwidth}{p{0.475\textwidth}p{0.475\textwidth}}
\toprule\\[-4mm]
\text{\bf \Set Postulate \refcolor{MyGreen}\ref{pst:DrawPoint}.} \em Draw an arbitrary point on the plane.
& \text{\bf \Set Problem \refcolor{MyGreen}\ref{pbm:DrawLine}.} \em Draw an arbitrary line on the plane. \\[2mm]
\begin{tcolorbox}
[colback=MyGreen!5!white, frame hidden, enhanced, attach boxed title to top right={yshift=-5.2mm,xshift=-0.2mm}, 
boxed title style={size=small, frame hidden, colback=MyGreen},
title=\footnotesize\Postulate{\refcolor{white}\ref{pst:DrawPoint}}]
\footnotesize $\Build{A}{\Point}{\top}$ 
\end{tcolorbox}& 
\begin{tcolorbox}
[colback=MyGreen!5!white, frame hidden, enhanced, attach boxed title to top right={yshift=-5.2mm,xshift=-0.2mm}, 
boxed title style={size=small, frame hidden, colback=MyGreen},
title=\footnotesize\Solution{\refcolor{white}\ref{pbm:DrawLine}}]
\footnotesize $\Build{x}{\Line}{\top}$
\end{tcolorbox} 
\\[-5mm]
\midrule
\text{\bf \Set Postulate \refcolor{MyGreen}\ref{pst:DrawDistinctPoint}.} \em 
Given an arbitrary point. Draw another point distinct from the given one.
& \text{\bf \Set Problem \refcolor{MyGreen}\ref{pbm:DrawDistinctLine}.} \em 
Given an arbitrary line. Draw another line that intersect the given one.\\[2mm]
\begin{tcolorbox}
[colback=MyGreen!5!white, frame hidden, enhanced, attach boxed title to top right={yshift=-5.2mm,xshift=-0.2mm}, 
boxed title style={size=small, frame hidden, colback=MyGreen},
title=\footnotesize\Postulate{\refcolor{white}\ref{pst:DrawDistinctPoint}}]
\footnotesize $\Forall{A}{\Point}{}\\[1mm]
\Build{B}{\Point}{A \neq B}$ 
\end{tcolorbox}& 
\begin{tcolorbox}
[colback=MyGreen!5!white, frame hidden, enhanced, attach boxed title to top right={yshift=-5.2mm,xshift=-0.2mm}, 
boxed title style={size=small, frame hidden, colback=MyGreen},
title=\footnotesize\Solution{\refcolor{white}\ref{pbm:DrawDistinctLine}}]
\footnotesize $\Forall{x}{\Line}{}\\[1mm]
\Build{y}{\Line}{x \nparallel y}$
\end{tcolorbox} 
\\[-5mm]
\midrule
\text{\bf \Set Postulate \refcolor{MyGreen}\ref{pst:DecideDistinctPoints}.} \em 
Given two distinct points $A$ and $B$ and an arbitrary point  $C$. 
Decide whether the point $C$ is distinct from the point $A$ or from the point $B$.
&\text{\bf \Set Problem \refcolor{MyGreen}\ref{pbm:DecideIntersectingLines}.} \em 
Given two intersecting lines $x$ and $y$ and an arbitrary third line $z$. 
Decide whether the line $z$ intersects the line $x$ or the line $y$.\\[2mm]
\begin{tcolorbox}
[colback=MyGreen!5!white, frame hidden, enhanced, attach boxed title to top right={yshift=-5.2mm,xshift=-0.2mm}, 
boxed title style={size=small, frame hidden, colback=MyGreen},
title=\footnotesize\Postulate{\refcolor{white}\ref{pst:DecideDistinctPoints}}]
\footnotesize$\Forall{A \w B \w C}{\Point}{}\\[1mm]
A \neq B \rightarrow \Or{A \neq C}{C \neq B}$
\end{tcolorbox}& 
\begin{tcolorbox}
[colback=MyGreen!5!white, frame hidden, enhanced, attach boxed title to top right={yshift=-5.2mm,xshift=-0.2mm}, 
boxed title style={size=small, frame hidden, colback=MyGreen},
title=\footnotesize\Solution{\refcolor{white}\ref{pbm:DecideIntersectingLines}}]
\footnotesize $\Forall{x \w y \w z}{\Line}{}\\[1mm]
x\nparallel y \rightarrow \Or{x \nparallel z}{z \nparallel y}$
\end{tcolorbox}
\\[-5mm]
\midrule
\text{\bf \Set Postulate \refcolor{MyGreen}\ref{pst:DrawPointOnLine}.} \em Given an arbitrary straight line. Draw a point on the line.
& \text{\bf \Set Problem \refcolor{MyGreen} \ref{pbm:DrawLineThroughPoint}.} \em Given an arbitrary point. Draw a line through this point. \\[1mm]
\begin{tcolorbox}
[colback=MyGreen!5!white, frame hidden, enhanced, attach boxed title to top right={yshift=-5.2mm,xshift=-0.2mm}, 
boxed title style={size=small, frame hidden, colback=MyGreen},
title=\footnotesize\Postulate{\refcolor{white}\ref{pst:DrawPointOnLine}}]
\footnotesize $\Forall{x}{\Line}{}\\[1mm]
\Build{A}{\Point}{A \in x}$
\end{tcolorbox}& 
\begin{tcolorbox}
[colback=MyGreen!5!white, frame hidden, enhanced, attach boxed title to top right={yshift=-5.2mm,xshift=-0.2mm}, 
boxed title style={size=small, frame hidden, colback=MyGreen},
title=\footnotesize\Solution{\refcolor{white}\ref{pbm:DrawLineThroughPoint}}]
\footnotesize $\Forall{A}{\Point}{}\\[1mm]
\Build{x}{\Line}{A \in x}$
\end{tcolorbox}
\\[-5mm]
\midrule
\text{\bf \Set Postulate \refcolor{MyGreen}\ref{pst:DrawDistinctPointOnLine}.} \em Given an arbitrary straight line and a point on it. Draw another point on this line distinct from the given one.& \text{\bf \Set Problem \refcolor{MyGreen} \ref{pbm:DrawDistinctLineThroughPoint}.} \em Given an arbitrary point and straight line passing through it. Draw another straight line, distinct from the given one and passing through the same point. \\[1mm]
\begin{tcolorbox}
[colback=MyGreen!5!white, frame hidden, enhanced, attach boxed title to top right={yshift=-5.2mm,xshift=-0.2mm}, 
boxed title style={size=small, frame hidden, colback=MyGreen},
title=\footnotesize\Postulate{\refcolor{white}\ref{pst:DrawDistinctPointOnLine}}]
\footnotesize $\ForallDouble{A}{\Point}{x}{\Line}{}\\[1mm]
A \in x \rightarrow \Build{B}{\Point}{A \neq B \land B \in x}$
\end{tcolorbox}& 
\begin{tcolorbox}
[colback=MyGreen!5!white, frame hidden, enhanced, attach boxed title to top right={yshift=-5.2mm,xshift=-0.2mm}, 
boxed title style={size=small, frame hidden, colback=MyGreen},
title=\footnotesize\Solution{\refcolor{white}\ref{pbm:DrawDistinctLineThroughPoint}}]
\footnotesize $\ForallDouble{A}{\Point}{x}{\Line}{}\\[1mm]
A \in x \rightarrow \Build{y}{\Line}{x \nparallel y \land A \in y}$
\end{tcolorbox}
\\[-5mm]
\midrule
\text{\bf \Set Postulate \refcolor{MyGreen}\ref{pst:DrawPointApartLine}.} \em Given an arbitrary line.
Draw a point that does not lie on this line.
& \text{\bf \Set Problem \refcolor{MyGreen}\ref{pbm:DrawLineNotThroughPoint}.} \em 
Given an arbitrary point. Draw a line that does not pass through this point.\\[1mm]
\begin{tcolorbox}
[colback=MyGreen!5!white, frame hidden, enhanced, attach boxed title to top right={yshift=-5.2mm,xshift=-0.2mm}, 
boxed title style={size=small, frame hidden, colback=MyGreen},
title=\footnotesize\Postulate{\refcolor{white}\ref{pst:DrawPointApartLine}}]
\footnotesize $\Forall{x}{\Line}{}\\[1mm]
\Build{A}{\Point}{A \notin x}$
\end{tcolorbox}& 
\begin{tcolorbox}
[colback=MyGreen!5!white, frame hidden, enhanced, attach boxed title to top right={yshift=-5.2mm,xshift=-0.2mm}, 
boxed title style={size=small, frame hidden, colback=MyGreen},
title=\footnotesize\Solution{\refcolor{white}\ref{pbm:DrawLineNotThroughPoint}}]
\footnotesize $\Forall{A}{\Point}{}\\[1mm]
\Build{x}{\Line}{A \notin x}$
\end{tcolorbox} \\[-5mm]
\midrule
\text{\bf \Set Postulate \refcolor{MyGreen}\ref{pst:DrawExtensionLine}.} \em 
Given two distinct points on the plane. Draw a straight line through both these points.
&\text{\bf \Set Postulate \refcolor{MyGreen}\ref{pst:DrawIntersectionPoint}.} \em Given two intersecting lines. 
Draw the point of their intersection.\\[2mm]
\begin{tcolorbox}
[colback=MyGreen!5!white, frame hidden, enhanced, attach boxed title to top right={yshift=-5.2mm,xshift=-0.2mm}, 
boxed title style={size=small, frame hidden, colback=MyGreen},
title=\footnotesize\Postulate{\refcolor{white}\ref{pst:DrawExtensionLine}}]
\footnotesize $\Forall{A \w B}{\Point}{}\\[1mm]
A\neq B\rightarrow\Build{x}{\Line}{A \w B \in x}$ 
\end{tcolorbox}& 
\begin{tcolorbox}
[colback=MyGreen!5!white, frame hidden, enhanced, attach boxed title to top right={yshift=-5.2mm,xshift=-0.2mm}, 
boxed title style={size=small, frame hidden, colback=MyGreen},
title=\footnotesize\Postulate{\refcolor{white}\ref{pst:DrawIntersectionPoint}}]
\footnotesize $\Forall{x \w y}{\Line}{}\\[1mm]
x \nparallel y \rightarrow \Build{A}{\Point}{A \in x \w y}$ 
\end{tcolorbox}
\\[-5mm]
\midrule
\text{\bf \Set Postulate \refcolor{MyGreen}\ref{pst:DecideUniqueExtensionLine}*\!.} \em 
Given two {\em intersecting} lines, on one of which two distinct points are marked. 
Decide which of these two points does not belong the other line.
& \text{\bf \Set Problem \refcolor{MyGreen}\ref{pol:DecidePointApartTwoIntersectingLines}.} \em 
Given two distinct points, through one of which pass two distinct lines.
 Decide which of these two lines does not pass through another point.\\[2mm]
\begin{tcolorbox}
[colback=MyGreen!5!white, frame hidden, enhanced, attach boxed title to top right={yshift=-5.2mm,xshift=-0.2mm}, 
boxed title style={size=small, frame hidden, colback=MyGreen},
title=\footnotesize\Postulate{\refcolor{white}\ref{pst:DecideUniqueExtensionLine}^{*}\s[-4]}]
\footnotesize $\ForallDouble{A \w B}{\Point}{x \w y}{\Line}{}\\[1mm]
A \neq B \rightarrow \In{A \w B}{x} \rightarrow x \nparallel y\rightarrow \Or{A \notin y}{B \notin y}$
\end{tcolorbox}& 
\begin{tcolorbox}
[colback=MyGreen!5!white, frame hidden, enhanced, attach boxed title to top right={yshift=-5.2mm,xshift=-0.2mm}, 
boxed title style={size=small, frame hidden, colback=MyGreen},
title=\footnotesize\Solution{\refcolor{white}\ref{pol:DecidePointApartTwoIntersectingLines}}]
\footnotesize $\ForallDouble{A \w B}{\Point}{x \w y}{\Line}{}\\[1mm]
A \neq B \rightarrow \In{A}{x \w y} \rightarrow x \nparallel y\rightarrow \Or{B \notin x}{B \notin y}$
\end{tcolorbox} 
\\[-5mm]
\bottomrule
\end{tabular*}
  \caption{Duality between distinction of points and intersection of lines.}
  \label{tab:duality}
\end{table*}

\section{3. Circles}

\subsection{3.1. Conceptualization of circle}

To draw a circle on the ground, the ancient geometers used a simple rope and two pointed pegs. 
They drove one peg into the ground where the center of the circle was supposed to be, and tied one end of the rope to it. 
After that, the rope was pulled to the point on the ground, through which the circle had to pass, and exactly at this point, another peg was attached to the second end of the rope. 
Then they walked in a circle, slightly pulling the rope and leaving a mark on the ground with the help of a second peg. 
Such a primitive tool for constructing circles was called {\em rope} compass or {\em collapsing} compass.

From this description of the process of constructing a circle, we see that two arbitrary points on the plane define a circle. 
The first point is called the {\em center} of the circle,  and about the second one it is said that it {\em lies on} the circle. 
Following the ancient Greeks, we believe that such a compass gives us the {\em only} way to draw a circle on the plane. 
This, in particular, means that if a circle was drawn on a plane, then it was drawn by rope compass via two previously constructed points.

In order to describe the basic properties of the circle, we need to introduce two new concepts.
First of all, we note that any circle divides all points of the plane onto two classes: those that "lie inside" this circle, and those that "lie outside" from it. 
In the same way, any pair of points on the line allows us to select those points of this line that "lie between" two given points or those that "lie outside this segment". 
In the next section, we formulate the necessary axioms, that will determine the basic property of the "betweenness" relation.

Another important property of a circle is the "equidistance" of all its points from the center of the circle. 
We shall formalize this property after we first introduce the concept of a segment formed by two arbitrary points on the plane and then define its length.

\subsection{3.2. Betweenness relation}

Ancient geometers in their reasoning widely used figurative human language.
Relying only on their drawings, they could come to the conclusion that "these two points always lie on opposite sides of a given straight line" or that "the constructed point is inside this triangle", or that "one segment is shorter than the others", etc. 
In other words, they took some of the "obvious" truths directly from the drawings. 
Similarly, on the basis of diagrams and visual representations only, all the main geometric figures, such as the "segment", "half-plane" or "circle", were defined, and their simplest properties were formulated.

\begin{marginfigure}[-15.5\baselineskip]\Mrg
\includegraphics[width=\linewidth]{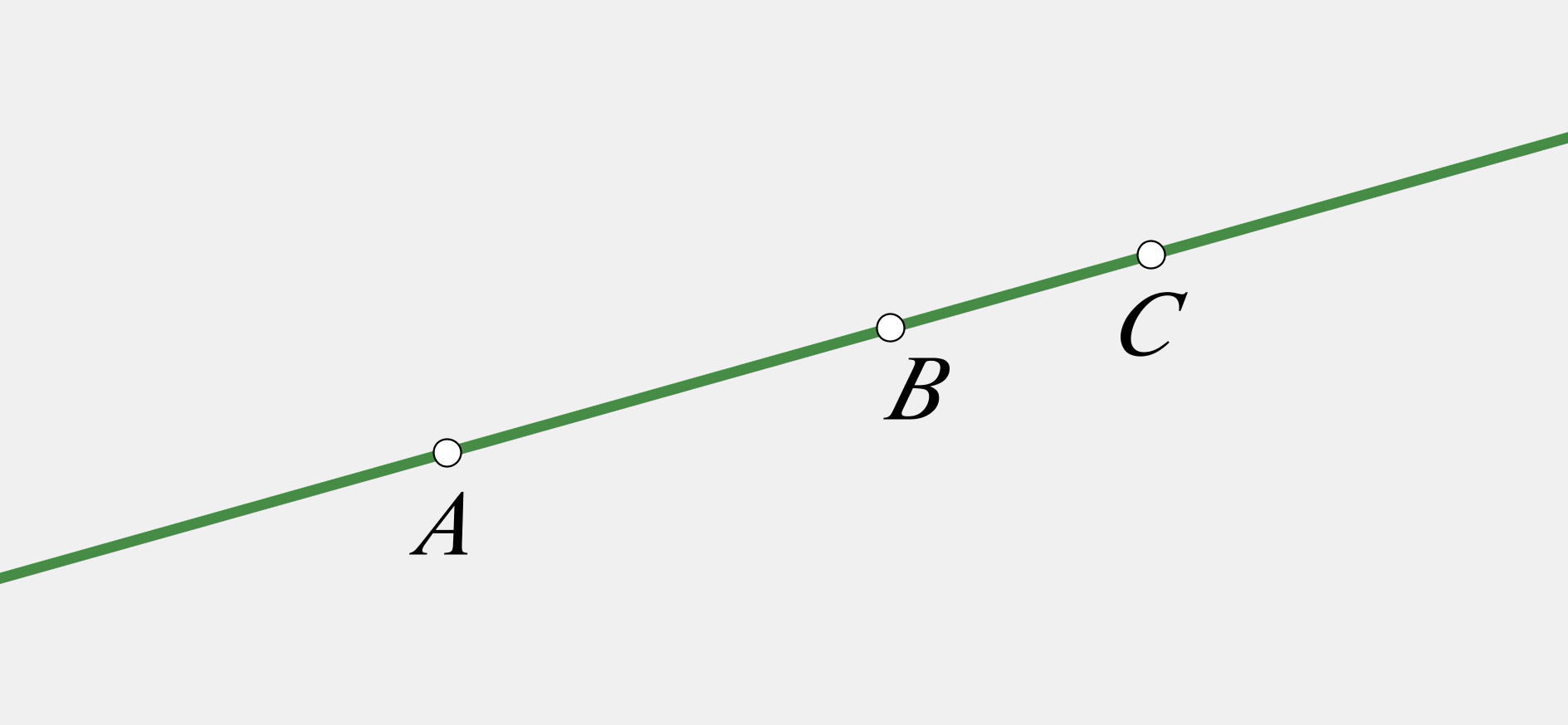}
\caption{Point $B$ lies between points $A$ and $C$.}
\label{fig:marginfig8}
\end{marginfigure}
Moritz Pasch have showed that as the main concept, that describes the mutual arrangement of points and straight lines on a plane, we can choose only one relationship "to lie between" and use it to replace all intuitive arguments of ancient geometers with strict definitions, axioms and proofs.

\mathnote{\Mrg
\begin{align}
&\Obj{A \w B \w C}{\Point},\nonumber\\
&\Bet{A}{B}{C}\nonumber\\
&\text{ ~~$\CC$  the point $B$ lies between the points $A$ and $C$}\nonumber\\
&\text{ ~~$\CC$  points $A$ and $C$ lie on the opposite half-lines}\nonumber\\
&\text{ ~~~~~from the point $B$}\nonumber
\end{align}}
Thus, the "betweenness" relation is undefined. 
This means that we must reveal the meaning of this relation through postulates and axioms. 
First of all, we note that whenever we say that one point lies between the other two, we mean that all three points in question are distinct and lie on one straight line. 
We also mean that the concept of "lying between" is symmetrical. Thus, we come to our first Axiom.
\vskip0.5\baselineskip
\begin{axiom}\Mrg
Given three points on the plane $A$, $B$ and $C$.
If the point $B$ lies between the points $A$ and $C$, then:
\mathnote{\Prop
\begin{align}
&\Forall{A \w B \w C}{\Point}{}\nonumber
\end{align}}
\begin{enumerate}
\item[a)] 
\mathnote{\Prop
\begin{align}
&\Bet{A}{B}{C}\rightarrow  A \neq C \nonumber
\end{align}}the points $A$ and $C$ are distinct;\label{axm:BetPs_DiPs}
\item[b)] 
\mathnote{\Prop
\begin{align}
&\Bet{A}{B}{C}\rightarrow \ColThree{A}{B}{C} \nonumber
\end{align}}all three points, $A$, $B$ and $C$ are collinear;\label{axm:BetPs_ColPs}
\item[c)] 
\mathnote{\Prop
\begin{align}
&\Bet{A}{B}{C}\rightarrow\Bet{C}{B}{A} \nonumber
\end{align}}the point $B$ also lies between the points $C$ and $A$;\label{axm:BetPs_sym}
\item[d)] 
\mathnote{\Prop
\begin{align}
&\Bet{A}{B}{C} \rightarrow \lnot \Bet{B}{A}{C}\nonumber
\end{align}}the point $A$ cannot lie between the points $B$ and $C$.\label{axm:BetPs_unique}
\end{enumerate}
\end{axiom}

\mathnote{\Mrg
\begin{align}
&\Obj{A\w B\w C\w D}{\Point}{} \nonumber\\
& \BetOX{A}{B}{C} \equiv \Bet{A}{B}{C} \lor (A = B \land B\neq C) \nonumber\\
& \BetXO{A}{B}{C} \equiv \Bet{A}{B}{C} \lor (A\neq B \land B = C) \nonumber\\
& \BetOO{A}{B}{C}\equiv \Bet{A}{B}{C} \lor A = B \lor B = C \nonumber
\end{align}}Sometimes it will be more convenient for us to use other betweenness symbols, that we define on margins. In these symbols we assume that some or all of the three points of the betweenness relation may coincide with each other. 
To indicate possible coincidence of points, we shall use white stars, while black stars will always indicate point's distinction. 

Two arbitrary points $A$ and $B$ define {\em segment} $\Seg{A}{B}$. 
\mathnote{\Mrg
\begin{align}
&\Obj{P \w A \w B}{\Point},\nonumber\\
&\BetOO{A}{P}{B} \nonumber\\
&\quad\CC\text{\rm  the point $P$ belongs to the segment $\Seg{A}{B}$}\nonumber
\end{align}}
The points $A$ and $B$ are called {\em endpoints} of the segment. 
Points lying between the endpoints of the segment, are called its {\em internal} points.
A point is said to belong to a segment iff it either lies between the endpoints of this segment, or coincides with one of its endpoints.

\subsection{3.3. Definition of segments and their lengths}

\mathnote{\Mrg
\begin{align}
&\Obj{s}{\text{Segment}}{}\nonumber\\
&s = \{A \w B\} \nonumber\\
&\quad A~\CC\text{ first point}\nonumber\\
&\quad B~\s[2]\CC\text{ second point}\nonumber
\end{align}}
\noindent 
Another important property of a circle is the "equidistance" of all its points from the center of the circle. 
To formalize this property, we must first introduce the concept of the length of a segment. 
This can be approached from two sides.

First, we can consider segments as pairs of points on a plane and assume, as Hilbert does, 
\mathnote{\Mrg
\begin{align}
&\Obj{s\w t}{\text{Segment}}{}\nonumber\\
&s\approx t\nonumber\\
&\quad\CC\text{ segments $s$ and $t$ are congruent}\nonumber\\
&{\mathbb L^{\rm +}}(s) \nonumber\\
&\quad\CC\text{ equivalence class for segment $s$}\nonumber\\[0.5\baselineskip]
&s\approx t\Iff {\mathbb L^{\rm +}}(s) = {\mathbb L^{\rm +}}(t)\nonumber
\end{align}}
that there is some equivalence (congruence) relation between different segments. This equivalence relation must be reflexive, symmetric, and transitive.
Based on this equivalence relation, we can define segment {\em lengths} as the corresponding classes of equivalent segments.

Alternativly, we can immediately enter the {\em length} of a segment as the basic undefined notion and assume that there is a mapping of segments to lengths (${\mathbb L^{\rm +}}$) and back (${\mathbb L^{\rm -}}$). In other words, each segment can be assigned its length, and for each length a corresponding segment on the plane can be found. 
\mathnote{\Mrg
\begin{align}
&{\mathbb L^{\rm +}} (s : \text{Segment}) : \text{Length}\nonumber\\
&{\mathbb L^{\rm -}} (d : \text{Length}) : \Build{s}{\text{Segment}}{{\mathbb L^{\rm +}}(s) = d}\nonumber
\end{align}}
Such mappings also define an equivalence relation for segments.

These two approaches are equivalent and we shall use both of them, however, for different purposes. 
In our opinion, the first is better suited for the formalization of new geometric figures on the plane with the required properties - circles, rays, flags, etc.
While the second approach is closer in spirit to the original Euclid's understanding of quantities --- lengths,  angle measures and areas.
 
\mathnote{
\begin{align}
&\Obj{A \w B \w C \w D}{\Point},\nonumber\\[0.5\baselineskip]
&\Len{A}{B} \equiv {\mathbb L^{\rm +}}({\{A \w B\}})\nonumber\\
&\quad \CC \text{  length of the segment $\Seg{A}{B}$}\nonumber\\[0.5\baselineskip]
&\Lcong{A}{B}{C}{D}\nonumber\\
&\quad \CC \text{  lengths of segments $\Seg{A}{B}$ and $\Seg{C}{D}$ are equal}  \nonumber
\end{align}}
In addition to the historical, there is yet another reason why for quantities we should adhere to the second approach and treat them as new undefined concepts. 
The thing is that one of the most important goals in the study of synthetic geometry is its arithmetization, as a result of which we will have to show that the quantities are similar to numbers and form some constructive field. In other words, we will have to show that quantities can be compared with each other, added to each other and even multiplied with each other (not in this paper) satisfying all necessary properties of a constructive filed.

We shall denote the length of the segment with endpoints $A$ and $B$ by double squared brackets $\Len{A}{B}$.

\mathnote{\Mrg
\begin{align}
&\Forall{A}{\Point}{}\nonumber\\
&\Len{A}{A} = 0 \nonumber
\end{align}}
A segment is called {\em null} if its endpoints coincide with each other.
The length of any such segment will be considered equal to {\em null}, and will be denoted by the symbol "$0$".

\vskip0.5\baselineskip
\begin{axiom}[Null segments]\Mrg
\mathnote{\Prop
\begin{align}
&\Forall{A \w B \w C}{\Point}{}\nonumber\\
&\Lcong{A}{B}{C}{C} \Iff A = B \nonumber
\end{align}}
The lengths of all null segments are equal to each other.
Any segment whose length is equal to the length of some null segment is itself a null segment. 
\label{axm:EqSs_EqPs}
\end{axiom}

\begin{marginfigure}[14\baselineskip]\Mrg
\includegraphics[width=\linewidth]{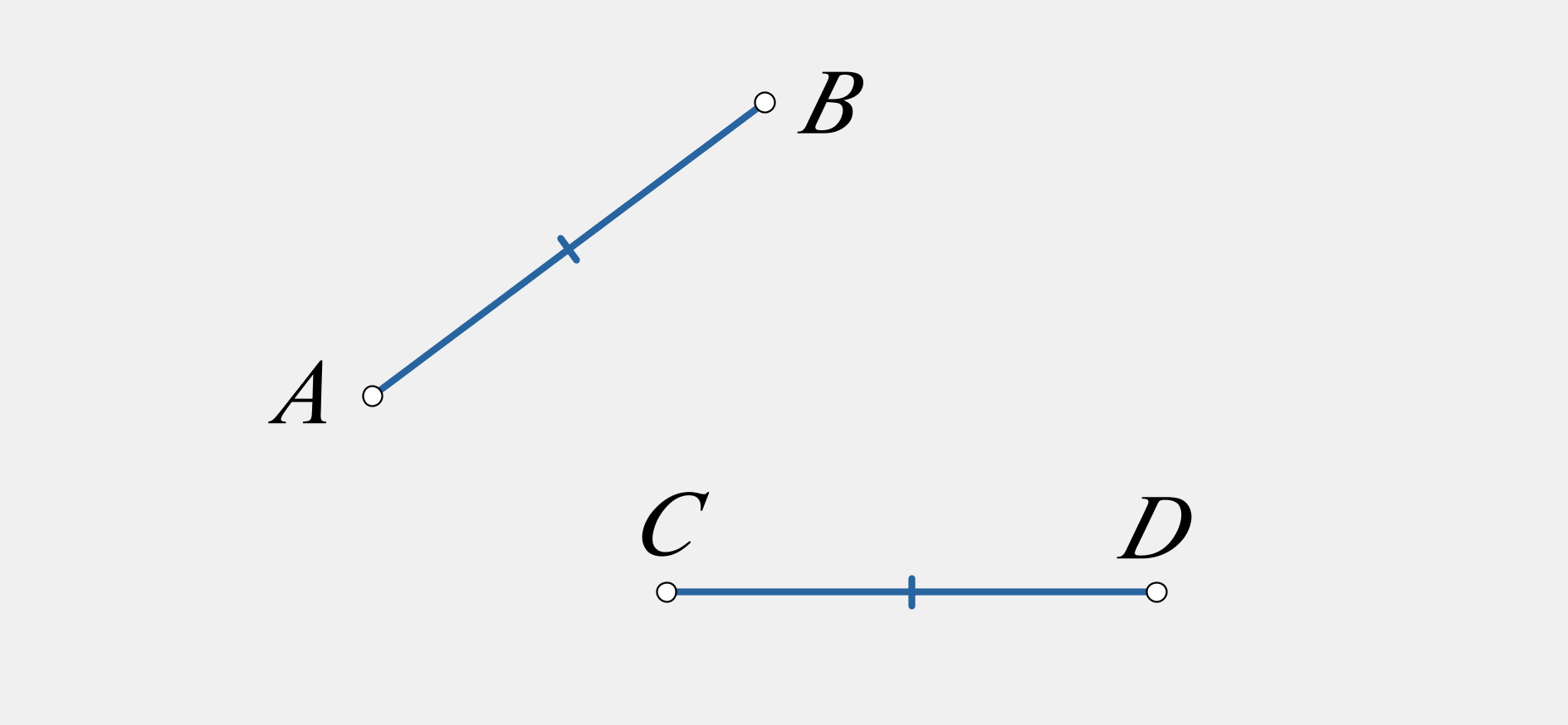}
\caption{An example of two segments of the same length ${\Lcong{A}{B}{C}{D}}$.}
\label{fig:EqSegments}
\end{marginfigure}
The endpoints of a {\em non-null} segment are by definition distinct from each other. 
A unique straight line can be drawn through these distinct points by Postulate \ref{pst:DrawExtensionLine}.
All internal points of the non-null segment will lie on this line by virtue of the Axiom \ref{axm:BetPs_ColPs}.
Thus, we can say that a non-null segment is a part of a line bounded by two points. 
All other points of this line are called {\em external} points of the segment.

Segments $\Seg{A}{B}$ and $\Seg{B}{A}$ both consist of the same points. 
This means that as geometric figures these two segments coincide with each other. 
This justifies the introduction of the following axiom:
\vskip0.5\baselineskip
\begin{axiom}\Mrg
\mathnote{\Prop
\begin{align}
&\Forall{A \w B}{\Point}{}\nonumber\\
&\Lcong{A}{B}{B}{A} \nonumber
\end{align}} 
The length of the segment does not depend on the order of its endpoints.
\label{axm:SegPs_sym}
\end{axiom}

\subsection{3.4. Drawing the point of intersection of a line with a circle}

\begin{marginfigure}[11\baselineskip]\Mrg
\includegraphics[width=\linewidth]{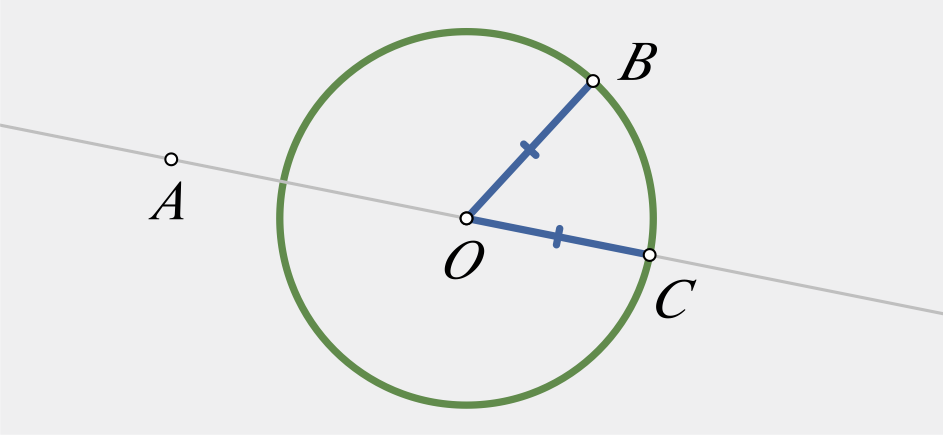}
\caption{Construction of a point $C$ of intersection of a straight line $\Lin{A}{O}$ and a circle $\Circ{O}{B}$.}
\label{fig:LineCircleIntersection}
\end{marginfigure}
So far, we have postulated the presence of only two distinct points on a given line.
From our everyday experience, we know that there should be many more distinct points on a line, perhaps even an infinite number of them. 
But how to draw them? 

In Euclidean geometry, there is one general "method" for constructing new points, which is based on the intersection of basic figures with each other. 
If the arrangement of two figures on a plane is such that their intersection point satisfying the necessary conditions is unique, then we can draw it.
We formalize this "method" in a series of postulates. 

We have already seen one such postulate --- a postulate about the intersection of two straight lines (Postulate \ref{pst:DrawIntersectionPoint}). 
Now we introduce another postulate, which gives us the opportunity to draw the point of intersection of the line and the circle.

\vskip0.5\baselineskip
\begin{postulate}[Line-Circle intersection]\Mrg
\mathnote{\Set
\begin{align}
&\Forall{O \w A \w B}{\Point}{}\nonumber\\
& A\neq O\rightarrow \Build{C}{\Point}{\BetOO{A}{O}{C} \land \Lcong{O}{B}{O}{C}}\nonumber
\end{align}}
Given two distinct points $A$ and $O$ and an arbitrary point $B$ on the plane. 
Draw the point $C$ of the intersection of the line $\Lin{A}{O}$ and the circle $\Circ{O}{B}$ such that the point $O$ belongs to the segment $\Seg{A}{C}$.
\label{pst:DrawIntersectionOfLineAndCircle}
\end{postulate}
Now, starting from two distinct points $A$ and $B$ on a straight line, with the help of our new tool --- compass --- we can build an unlimited number of points on this line.

We shall later formulate one more postulate stemming from the same general "method". 
Namely, the postulate of the intersection of two circles.

\subsection{3.5. Line divides plane into two half-planes}

\begin{marginfigure}[13.5\baselineskip]\Mrg
\includegraphics[width=\linewidth]{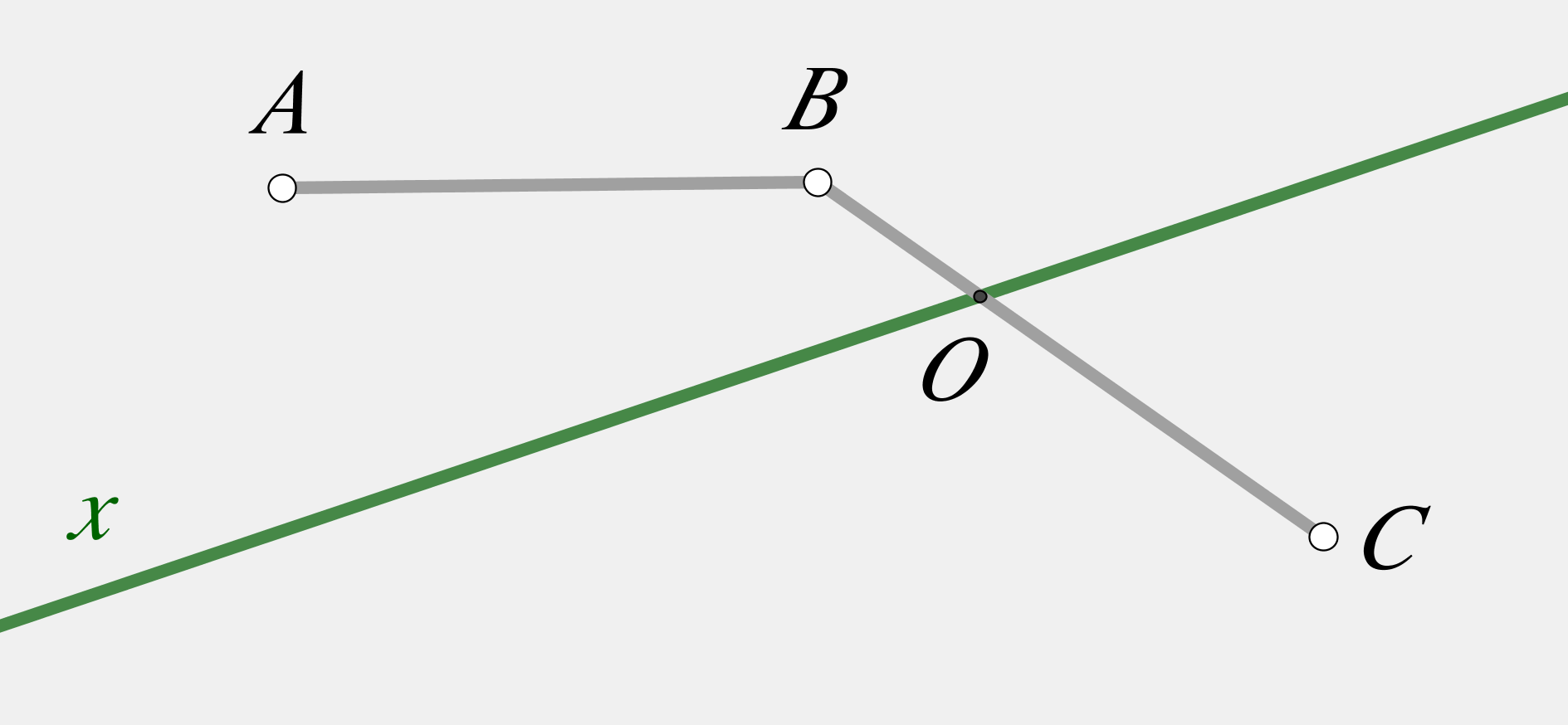}
\caption{The points $A$ and $B$ lie on same half-plane of the line $x$, and the points $B$ and $C$ lie on different half-planes of the line.}
\label{fig:SameAndOppositeSides_1}
\end{marginfigure}
Let us look at Fig. \ref{fig:SameAndOppositeSides_1}. 
In this figure, points $A$ and $B$ {\em lie on the same half-plane} of line $x$, and distinct points $B$ and $C$ {\em lie on different half-planes} of the line. 
This, however, is rather confusing, since we have not yet defined what a half-plane is and how many different half-planes are associated with a line on a plane.

To give these terms a precise meaning, we can also note that the line $x$ does not intersect the segment $\Seg{A}{B}$, but it does intersects the segment $\Seg{B}{C}$. 
Here, the notion "line intersects a segment" can be given a precise meaning. Namely, that there exists a point $O$ on the line $x$, which is also an internal point of the segment $\Seg{B}{C}$. Therefore, the point $O$ lies between the endpoints $B$ and $C$ of this segment.

Now, using the undefined "betweenness" relation, we can give precise definitions of the relations "lie on same half-plane" and "lie on different half-planes" from a given straight line.
\vskip0.5\baselineskip
\begin{definition}\Mrg
\mathnote{\Def
\begin{align}
&\ObjDouble{A \w B}{\Point}{x}{\Line} \nonumber
\end{align}}
Given a line $x$ and two points $A$ and $B$ not lying on this line.
The points $A$ and $B$ are said to:
\begin{enumerate}
\item[a)] 
\mathnote{\Def
\begin{align}
&\OppositeSide{A}{x}{B}\nonumber\\
&\quad\equiv \NotIn{A \w B}{x} \land \Exists{O}{\Point}{O\in x \land \Bet{A}{O}{B}} \nonumber
\end{align}}lie on different half-planes of the line $x$, iff there exists a point on this line between the points $A$ and $B$;
\item[b)] 
\mathnote{\Def
\begin{align}
&\SameSide{x}{A}{B} \nonumber\\
&\quad\equiv  \NotIn{A \w B}{x}\land \NotExists{O}{\Point}{O \in x \land \Bet{A}{O}{B}} \nonumber
\end{align}}lie on one same half-plane of the line $x$ iff no point of this line lies between the points $A$ and $B$.
\end{enumerate}
\label{def:SameSide}
\label{def:SameOppSide}
\end{definition}

\begin{marginfigure}[12.5\baselineskip]\Mrg
\includegraphics[width=\linewidth]{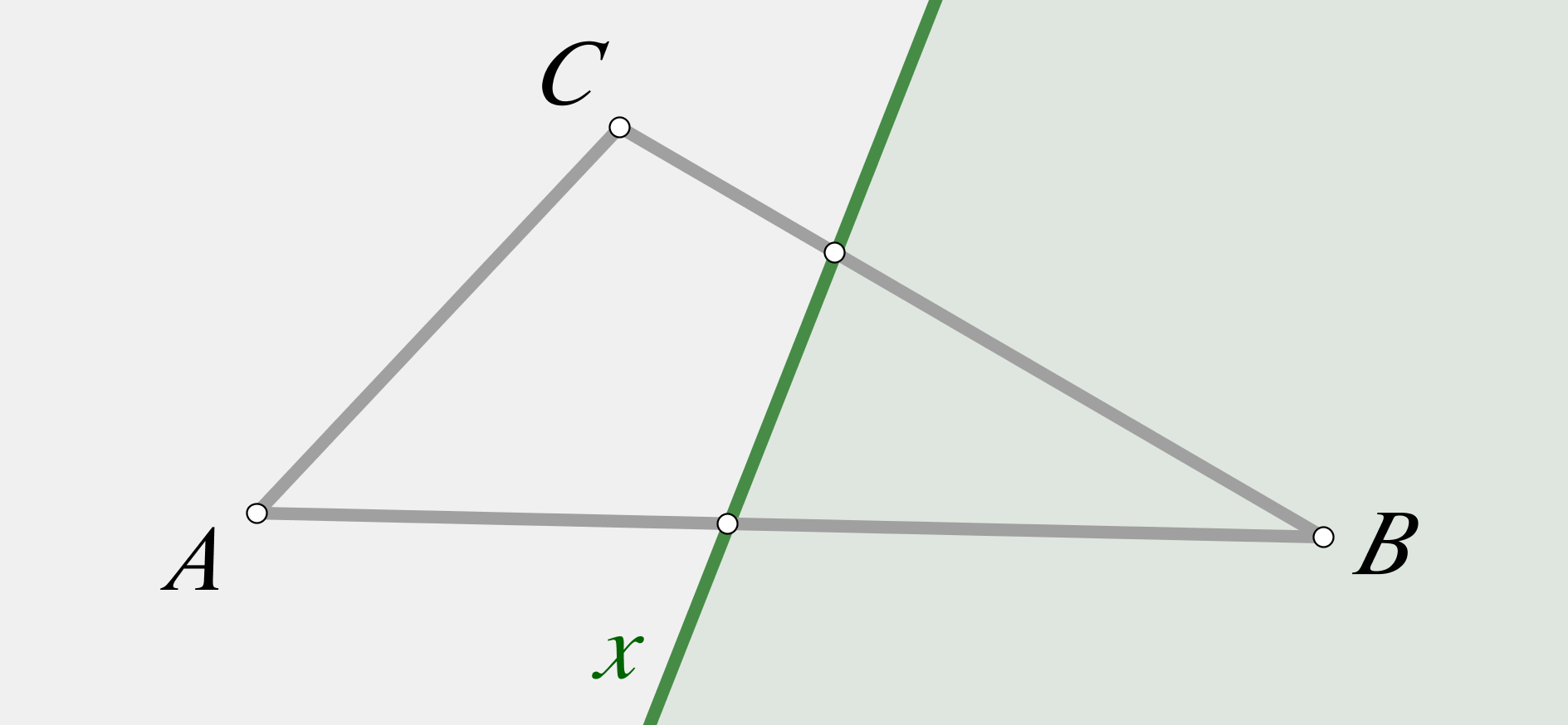}
\caption{Pasch Postulate.}
\label{fig:PaschPostulate}
\end{marginfigure}
Let us look at Fig. \ref {fig:PaschPostulate}. 
In this figure, the points $A$ and $B$ lie on different half-planes of the line $x$. 
According to the Definition \ref {def:SameOppSide}, in this case the points can not lie on the straight line itself. 
Let us now consider a point $C$, which also does not lie on the line $x$. 
Studying the configuration of points in our figure, we see that the points $B$ and $C$ lies on the different half-planes with respect to the straight line $x$. 
From our experience, we know that other options are also possible for the location of the point $C$ relative to the straight line $x$ on the plane.
However, in all cases it will necessarily lie with respect to a straight line $x$ on different half-planes either with $A$ or with $B$.
Moreover, we can always {\em decide} which of these two alternatives really takes place. 
Thus, we come to the following postulate:
\vskip0.5\baselineskip
\begin{postulate}[Morits Pasch]\Mrg
\mathnote{\Set
\begin{align}
&\ForallDouble{A \w B \w C}{\Point}{x}{\Line}{}\nonumber\\
&\OppositeSide{A}{x}{B} \rightarrow C \notin x \rightarrow \Or{\OppositeSide{A}{x}{C} }{\OppositeSide{C}{x}{B}}\nonumber
\end{align}}
Given three points $A, B, C$ and a line $x$ that does not pass through the points and intersects the segment $\Seg{A}{B}$.
Decide whether the line $x$ intersects the segment $\Seg{A}{C}$ or the segment $\Seg{C}{B}$.
\label{pst:DecideOppositeSide}
\end{postulate}
As in the case of Heyting and von Plato's postulates, the Pash postulate does not create any new objects on the plane. 
It describes the ability of an observer to decide which of the two alternatives really takes place for any valid configuration of points and lines on the plane.

From this postulate it follows that if a straight line intersects one side of the triangle and does not pass through its vertices, then it will necessarily intersect some other side of this triangle. Moreover, we can always decide which of the two other sides of the triangle was actually intersected.

\vskip0.5\baselineskip
\begin{theorem}\Mrg
\mathnote{\Prop
\begin{align}
&\ForallDouble{A\w B \w C}{\Point}{x}{\Line}{} \nonumber
\end{align}}
Given a line on the plane, the relation "to lie on the same half-plane" of this line is the equivalence relation.
\begin{description}
\item[\Prop\textsc{reflexivity.}]
\mathnote{\Prop
\begin{align}
&A \notin x \rightarrow \SameSide{x}{A}{A} \nonumber
\end{align}}If the point $A$ does not lie on line $x$, then the point $A$ lies on the half-plane of line $x$ with itself.
\item[\Prop\textsc{symmetry.}] 
\mathnote{\Prop
\begin{align}
&\SameSide{x}{A}{B}\rightarrow \SameSide{x}{B}{A}\nonumber
\end{align}}If the points $A$ and $B$ lie on the same half-plane of line $x$, then the points $B$ and $A$ also lie on the same half-plane.
\item[\Prop\textsc{transitivity.}]
\mathnote{\Prop
\begin{align}
&\SameSide{x}{A}{B}\rightarrow \SameSide{x}{B}{C}\rightarrow \SameSide{x}{A}{C}\nonumber
\end{align}}If the points $A$ and $B$ lie on the same half-plane of line $x$, and the points $B$ and $C$ lie on the same half-plane of line $x$, then the points $A$ and $C$ also lie on the same half-plane of line $x$.
\end{description}
\label{lem:ParallelEquivalence}
\end{theorem}
Similar theorem can also be proved for the relation "to lie on different half-planes" from the given line.
\begin{theorem}\Mrg
\mathnote{\Prop
\begin{align}
&\ForallDouble{A\w B\w C}{\Point}{x}{\Line}{} \nonumber
\end{align}}
Given a line on the plane, the relation "to lie on different half-planes" from this line has the following properties.
\begin{description}
\item[\Prop\textsc{irreflexivity.}]
\mathnote{\Prop
\begin{align}
&\OppositeSide{A}{x}{B} \rightarrow A \neq B\nonumber
\end{align}}If the points $A$ and $B$ lie on different half-planes of line $x$, then these points are distinct.
\item[\Prop\textsc{symmetry.}] 
\mathnote{\Prop
\begin{align}
&\OppositeSide{A}{x}{B}\rightarrow \OppositeSide{B}{x}{A}\nonumber
\end{align}}If the points $A$ and $B$ lie on different half-planes of line $x$, then the points $B$ and $A$ also lie on different half-planes of line $x$.
\item[\Prop\textsc{cotransitivity OOS.}] 
\mathnote{\Prop
\begin{align}
&\OppositeSide{A}{x}{B}\rightarrow \OppositeSide{B}{x}{C}\rightarrow \SameSide{x}{A}{C}\nonumber
\end{align}}If the points $A$ and $B$ lie on different half-planes of line $x$, and the points $B$ and $C$ lie on different half-planes of line $x$, then the points $A$ and $C$ lie on same half-plane of line $x$.
\item[\Prop\textsc{cotransitivity OSO.}] 
\mathnote{\Prop
\begin{align}
&\OppositeSide{A}{x}{B}\rightarrow \SameSide{x}{B}{C}\rightarrow \OppositeSide{A}{x}{C}\nonumber
\end{align}}If the points $A$ and $B$ lie on different half-planes of line $x$, and the points $B$ and $C$ lie on the same half-plane of line $x$, then the points $A$ and $C$ lie on different half-planes of line $x$.\\
 $~$\hfill{\Prop\rm \scriptsize / Hilbert, Chapter 1 : Theorem 8}
\end{description}
\label{thm:OppSideEq}
\end{theorem}

Finally, we can prove that any straight line divides the plane into exactly two half-planes.
This explains why these two different half-planes are often called {\em opposite sides} from a given line on a plane.
Thus, if a point does not lie on a straight line, then the Pasch's postulate gives us a tool to decide which of the two opposite sides of the line this point is actually on.
In the next section, we shall also prove that the half-plane is a convex figure.

\subsection{3.6. Point divides line into two half-lines}

Consider now points on a line.
The indefinite relationship "between" distinguishes such an arrangement of points on a straight line, in which the two extreme points "lie on different half-lines" from the midpoint. We can now define in a sense an opposite relation, which distinguishes  such an arrangement of points on a straight line, in which two points lie "on the same half-line" from a given point.\\ [-0.5\baselineskip]
\begin{definition}\Mrg
\mathnote{\Def
\begin{align}
&\Obj{O\w A \w B}{\Point},\nonumber\\
&\SameRay{O}{A}{B}\equiv O \neq A \land O\neq B \land \ColThree{O}{A}{B} \land \lnot \Bet{A}{O}{B} \nonumber\\
&\quad\CC\text{\rm  the points $A$ and $B$ lie on the same half-line }\nonumber\\
&\quad\quad\text{\rm ~from the point $O$}\nonumber
\end{align}} 
It is said that the points $A$ and $B$ lie on the same half-line from the point $O$, iff both these points are distinct from the point $O$, all three points are collinear, and the point $O$ does not lie between the points $A$ and $B$.
\label{def:SameRay}
\end{definition}

First we prove that this relation is an equivalence relation, i.e. it is reflexive, symmetric and transitive.
\vskip0.5\baselineskip
\begin{theorem}\Mrg
\mathnote{\Prop
\begin{align}
&\Forall{O\w A\w B\w C}{\Point}{}\nonumber
\end{align}} 
The relation "to lie on same half-line" is an equivalence relation.
\begin{description}[leftmargin=0pt]
\item[\Prop\textsc{reflexivity.}]
\mathnote{\Prop
\begin{align}
&O \neq A\rightarrow\SameRay{O}{A}{A}\nonumber
\end{align}}If the points $O$ and $A$ are distinct, then the point $A$ lies with itself on the same half-line from the point $O$.
\item[\Prop\textsc{symmetry.}] 
\mathnote{\Prop
\begin{align}
&\SameRay{O}{A}{B}\rightarrow\SameRay{O}{B}{A}\nonumber
\end{align}}If the points $A$ and $B$ lie on the same half-line from the point $O$, then the points $B$ and $A$ also lie on the same side of the point $O$.
\item[\Prop\textsc{transitivity.}]
\mathnote{\Prop
\begin{align}
&\SameRay{O}{A}{B}\rightarrow\SameRay{O}{B}{C}\rightarrow\SameRay{O}{A}{C}\nonumber
\end{align}}If the points $A$ and $B$ lie on the same half-line from the point $O$, and the points $B$ and $C$ lie on the same half-line from the point $O$, then the points $A$ and $C$ also lie on the same half-line from the point $O$.
\end{description}
\end{theorem}

Then, we prove that "betweenness" relation can equivalently be expressed via the relation "to lie on the same half-line".
\begin{theorem}\label{lem:BetPs_SR}\Mrg
\mathnote{\Prop
\begin{align}
&\Forall{A\w B\w C}{\Point}{}\nonumber\\
&\Bet{A}{B}{C} \Iff \SameRay{A}{B}{C}\land \SameRay{C}{B}{A} \nonumber
\end{align}} 
The point $B$ lies between the points $A$ and $C$ iff the points $B$ and $C$ lie on the same half-line from the point $A$, while the points $B$ and $A$ lie on same half-line from point $C$.
\end{theorem}

Finally, starting from two distinct points on a given straight line, (Postulates \ref{pst:DrawPointOnLine} and  \ref{pst:DrawDistinctPointOnLine}), we can construct infinitely many distinct points incident to this line.

\vskip0.5\baselineskip
\begin{problem}\Mrg
\mathnote{\Set
\begin{align}
&\Forall{A \w B}{\Point}{}\nonumber
\end{align}}Given two distinct points $A$ and $B$. 
\begin{enumerate}
\item[a)]
\mathnote{\Set
\begin{align}
& A \neq B\rightarrow  \Build{O}{\Point}{\Bet{A}{O}{B}}\nonumber
\end{align}}Draw a point $O$ between the distinct points $A$ and $B$.
\label{prb:DrawPointBetween}
\item[b)]
\mathnote{\Set
\begin{align}
& A\neq B\rightarrow  \Build{C}{\Ray}{\Bet{A}{B}{C}}\nonumber
\end{align}}Draw a point $C$ such that point $B$ is between $A$ and $C$.\\
 $~$\hfill{\Set\rm \scriptsize / Hilbert, Chapter 1 : Theorem 3}
\end{enumerate}
\end{problem}

It will also be convenient for us to use the four-point order relation.
\mathnote{\Mrg
\begin{align}
&\BetFour{A}{B}{C}{D} \nonumber\\
&\quad\equiv\Bet{A}{B}{C} \s[-2]\land\s[-2] \Bet{A}{B}{D} \s[-2]\land\s[-2] \Bet{A}{C}{D} \s[-2]\land\s[-2] \Bet{B}{C}{D}\nonumber
\end{align}}
We will say that the four given points are located on a straight line in the established order, if any three of them, taken in the same order, satisfy the betweenness relation.
It can be proved that "betweenness" relation satisfy some sort of transitivity properties.
\begin{theorem}\Mrg
\mathnote{\Prop
\begin{align}
&\Forall{A\w B\w C\w D}{\Point}{}\nonumber
\end{align}} 
The "betweenness" relation has the following properties:
\begin{enumerate}[leftmargin=*]
\item[a)]
\mathnote{\Prop
\begin{align}
& \Bet{A}{B}{C} \rightarrow  \Bet{A}{C}{D} \rightarrow \BetFour{A}{B}{C}{D}\nonumber
\end{align}}If the point $B$ lies between the points $A$ and $C$, and the point $C$ lies between the points $B$ and $D$, then all for points $A$, $B$, $C$ and $D$ are ordered.
\item[b)]
\mathnote{\Prop
\begin{align}
&\Bet{A}{B}{C} \rightarrow \Bet{B}{C}{D}\rightarrow \BetFour{A}{B}{C}{D}\nonumber
\end{align}}If the point $B$ lies between the points $A$ and $C$, and the point $C$ lies between the points $A$ and $D$, then all for points $A$, $B$, $C$ and $D$ are ordered.
\end{enumerate}
\end{theorem}

Finally, we can prove that whatever is a line, any point on the line divide it into exactly two half-lines.
This explains why these two different half-lines are often called {\em opposite sides} from a given point on a line.

\subsection{3.7. Definition of circles as equivalence classes}

\mathnote{\Mrg
\begin{align}
&\Obj{r}{\text{Circle}}{}\nonumber\\
&r = \{r_0\s[4] r_1\} \nonumber\\
&\quad r_0\CC\text{center point}\nonumber\\
&\quad r_1\CC\text{circumference point}\nonumber
\end{align}}
Another general "method" of Euclidean geometry is to define new geometric figures using some equivalence relations.
In fact, even the definition of a circle (as points equidistant from the center) requires an equivalence relation to formalize it.

To this end let us consider pairs of points which are not necessarily distinct, since we are also considering circles with null radius.
The first point of the pair we shall call its {\em center}, while the second -- its {\em circumference} point.
\mathnote{\Mrg
\begin{align}
&\Obj{r \w s}{\text{Circle}}{}\nonumber\\
&r \approx s\equiv (r_0 = s_0)\land \Lcong{r_0}{r_1}{s_0}{s_1}\nonumber
\end{align}}
The pairs are called equivalent iff their centers coincide and their circumference points are equidistant from the center.

Now we may consider {\em equivalence classes} with respect to this equivalence relation. 
\mathnote{\Mrg
\begin{align}
&\Obj{r\w s}{\text{Circle}}{}\nonumber\\
& r\approx s\Iff{\mathbb C^{\rm +}}(r) = {\mathbb C^{\rm +}}(s)\nonumber
\end{align}}
These equivalence classes ${\mathbb C^{\rm +}}(r)$ and ${\mathbb C^{\rm +}}(s)$ are nothing but {\em circles} as geometric figures, while $r$ and $s$ are pairs of points {\em representing} these circles. 
Unfortunately, Coq does not have instruments to produce quotient-types and to handle equivalence classes directly.
When programming in Coq, we always have to deal not with the circles themselves, but with pairs of points representing them, and with an equivalence relation on the pairs.

This also explains why we decided not to introduce circles into our formalization as basic objects on the same ground as points and lines.
For all practical purposes, it will be sufficient for us not to use the circles themselves, as basic geometric figures, but only the basic notions associated with them --- betweenness relation and the length of the segment. This way, we avoid introducing an additional undefined incidence relation (between points and circles) and the corresponding set of axioms that would ensure that the new concept of incidence is consistent with the concept of equidistance of points from the center of the circle.

Now we can easily define the incidence relation between a point and a circle, assuming that the latter is represented by two points on a plane, as discussed above.
\mathnote{\Mrg
\begin{align}
&O\w A\w P : \Point,\nonumber\\
&\Lcong{O}{P}{O}{A} \nonumber\\
&\quad\CC\text{\rm  the point $P$ lies on the circle $\Circ{O}{A}$} \nonumber
\end{align}}
It is said that the point $P$ lies on the circle $\Circ{O}{A}$ if it is located at the same distance from the center of the circle $O$ as the circumference point $A$ of the circle.

Let us now define some other notions and terms related to the {\em circle} as a geometric figure. 
The length of the segment $\Seg{O}{A}$ is called the {\em radius} of the circle.
A circle is called {\em degenerate} if its radius is zero. 
The line segment joining two distinct points on some non-degenerate circle is called a {\em chord} of the circle. 
\mathnote{\Mrg
\begin{align}
&O\w A\w B : \Point,\nonumber\\
&\Middle{A}{O}{B} \equiv \Bet{A}{O}{B} \land \Lcong{O}{A}{O}{B}\nonumber\\
&\quad \CC \text{  the point $O$ lies midway between}\nonumber\\
&\quad\quad\text{~~points $A$ and $B$} \nonumber
\end{align}}
If the center of the circle lies between the endpoints of the chord, then such chord is called {\em diameter} of the circle. 

We say that a point lies {\em outside} of a circle iff there is some point on circumference between the center of the circle and that point. 
We also say that a point lies {\em inside} a circle iff this point does not lie outside the circles and does not lie on the circle.

\section{4. Rays}

\subsection{4.1. Definition of rays as equivalence classes}

\mathnote{\Mrg
\begin{align}
&\Obj{a}{\Ray}{}\nonumber\\
&a = \{a_0\bs a_1\} \nonumber\\
&\quad a_0\CC\text{originating point}\nonumber\\
&\quad a_1\CC\text{directing point}\nonumber
\end{align}}
To define the rays on the plane, let us consider pairs of distinct points.
This time the first point of the pair we shall call the {\em origin} of the ray, while the second point -- its {\em director}.

On these pairs of distinct points, we can define equivalence relation, which we call the ray coincidence relation.
It is said that rays coincide with each other iff
\mathnote{\Mrg
\begin{align}
&\Obj{a \w b}{\Ray}{}\nonumber\\
&a \approx b\equiv (a_0 = b_0)\land\SameRay{a_0}{a_1}{b_1}\nonumber
\end{align}}
their originating points coincide, and their directing points lie on the same half-line from their common originating point.

\mathnote{\Mrg
\begin{align}
&\Obj{a\w b}{\Ray}{}\nonumber\\
&a\approx b\Iff {\mathbb R^{\rm +}}(a) = {\mathbb R^{\rm +}}(b)\nonumber
\end{align}}
Equivalence classes defined by this equivalence relation will be called {\em rays}. 
Intuitively, a ray is a figure consisting of an originating point and all points of a line lying on the same half-line from its origin.

\mathnote{\Mrg
\begin{align}
&\Obj{a \w b}{\Ray},\nonumber\\
&a\overset{*}{\approx} b \equiv (a_0 = b_0)\land\Bet{a_1}{a_0}{b_1}\nonumber\\
&\quad\CC\text{\rm  rays $a$ and $b$ are opposite to each other}\nonumber
\end{align}} 
We have already proved that any point on the line divides the line into two half-lines.
This suggest the following definition.
It is said that two rays are opposite to each other iff they both originate from the same point go in opposite directions.
In other words, their common point of origin is located between their directing points.
As a result, the opposing rays together form one straight line.

Finally, for every ray $a$ on the plane we can draw unique opposite ray which we denote $\Op{a}$.

\subsection{4.2. Direction on line}

Our discussion of opposite rays on a line also suggests that on each straight line we could consider two opposite {\em directions}.
To formalize the notion of direction on a line, let us consider another equivalence relation on pairs of distinct points.
\vskip0.5\baselineskip
\begin{definition}[Rays direction]\Mrg
\mathnote{\Def
\begin{align}
&\Obj{a\w b}{\Ray}{} \nonumber\\
&a \sim b\equiv (a_0 = b_0\land\SameRay{a_0}{a_1}{b_1})\lor\nonumber\\
&\s[45]\lor(\SameRay{a_0}{a_1}{b_0}\land\Bet{a_0}{b_0}{b_1})\lor\nonumber\\
&\s[45]\lor(\SameRay{b_0}{b_1}{a_0}\land\Bet{b_0}{a_0}{a_1})\nonumber
\end{align}}
It is said that that rays $a$ and $b$ have the same direction iff either of the three alternatives is valid:
\begin{enumerate}
\item[--] rays $a$ and $b$ coincide;
\item[--] the origin of the ray $b$ belongs to the ray $a$, while the origin of the ray $a$ belongs to the ray opposite to $b$;
\item[--] the origin of the ray $a$ belongs to the ray $b$, while the origin of the ray $b$ belongs to the ray opposite to $a$.
\end{enumerate}
\end{definition}

\mathnote{\Mrg
\begin{align}
&\Obj{a\w b}{\Ray}{}\nonumber\\
&a\sim b\Iff {\mathbb D^{\rm +}}(a) = {\mathbb D^{\rm +}}(b)\nonumber
\end{align}}
This equivalence holds only for those rays that lie on the same straight line. 
Therefore, equivalence classes for this equivalence relation define two opposite directions on each straight line.

Moreover, this equivalence relation is decidable, i.e. given any two rays on a line, we can decide if they have same or opposite directions.

\begin{problem}\Mrg
\mathnote{\Set
\begin{align}
&\Forall{a \w b}{\Ray}{}\nonumber\\
& \ColTwo{a}{b} \rightarrow\Or{a \sim b}{a\sim \Op{b}}\nonumber
\end{align}
}Given two collinear rays $a$ and $b$. Decide whether they have same or opposite direction.
\label{prb:DrawPointBetween}
\end{problem}

\begin{theorem}\Mrg
\mathnote{\Prop
\begin{align}
&\Forall{a \w b}{\Ray}{}\nonumber\\
&\lnot\s(\s a\sim b \land a\sim\Op{b}\s)\nonumber
\end{align}} 
Two rays cannot simultaneously have both same and opposite direction.
\end{theorem}

\subsection{4.3. Order of concurrent rays}

In analogy with the order of points on the straight line, we now would like to turn to the study of the order of concurrent rays, that is rays originating from one point.

The following notions are convenient to define order of rays on plane. This notions are somewhat similar to analogous relations for points.

\vskip0.5\baselineskip
\begin{definition}[Rays divergence]\Mrg
\mathnote{\Def
\begin{align}
&\Obj{a \w b}{\Ray},\nonumber\\
&\Apart{a}{b}\equiv (a_0 = b_0)\land\lnot\ColThree{a_0\w}{a_1\w}{b_1}\nonumber\\
&\quad\CC\text{\rm  rays $a$ and $b$ diverge from each other}\nonumber
\end{align}} 
It is said that two rays are divergent iff they are concurrent and non-collinear.
\end{definition}

\begin{definition}[Rays on same side]\label{def:RaysSameSide}\Mrg
\mathnote{\Def
\begin{align}
&\Obj{a \w b \w c}{\Ray}{}\nonumber\\
&\SameRay{a}{b}{c}\equiv (a_0 = b_0 = c_0)\land\SameHalfplane{a_0}{a_1}{b_1}{c_1}\nonumber
\end{align}
}
It is said that rays $b$ and $c$ are on the same side from ray $a$ iff all three rays are concurrent and direction points of rays $b$ and $c$ are on the same side of the line passing through the points of ray $a$. 
\end{definition}

First, we prove that this relation is an equivalence relation, i.e. it is reflexive, symmetric and transitive.

\vskip0.5\baselineskip
\begin{theorem}\Mrg
\mathnote{\Prop
\begin{align}
&\Forall{a\w b\w c\w d}{\Ray}{}\nonumber
\end{align}} 
To lay on the same side of ray is an equivalence relations.
\begin{description}[leftmargin=0pt]
\item[\Prop\textsc{reflexivity.}]
\mathnote{\Prop
\begin{align}
&\Apart{a}{b}\rightarrow\SameRay{a}{b}{b}\nonumber
\end{align}}Given two divergent rays $a$ and $b$, then the ray $b$ lies on the same side of the ray $a$ with itself .
\item[\Prop\textsc{symmetry.}] 
\mathnote{\Prop
\begin{align}
&\SameRay{a}{b}{c}\rightarrow\SameRay{a}{c}{b}\nonumber
\end{align}}If rays $b$ and $c$ lie on the same side of ray $a$, then rays $c$ and $b$ also lie on the same side of the ray $a$.
\item[\Prop\textsc{transitivity.}]
\mathnote{\Prop
\begin{align}
&\SameRay{a}{b}{c}\rightarrow\SameRay{a}{c}{d}\rightarrow\SameRay{a}{b}{d}\nonumber
\end{align}}If rays $b$ and $c$ lie on the same side of ray $a$, and rays $c$ and $d$ lie on the same side of the ray $a$, then  rays $b$ and $d$ lie on the same side of ray $a$.
\end{description}
\end{theorem}

Then, we introduce the "betweenness" relation of the rays and prove that this notion, although it does not literally coincide with "betweenness" relation for points, still has many similar properties.

\vskip0.5\baselineskip
\begin{definition}[Rays betweenness]\Mrg
\mathnote{\Def
\begin{align}
&\Obj{a \w b \w c}{\Ray}{}\nonumber\\
&\Bet{a}{b}{c}\equiv \SameRay{a}{b}{c}\land\SameRay{c}{b}{a}\nonumber
\end{align}}
It is said that that the ray $b$ is between rays $a$ and $c$ iff rays $b$ and $c$ are on the same side from the ray $a$ while rays $b$ and $a$ are on the same side from the ray $c$.
\end{definition}

\begin{theorem}\Mrg
If ray $b$ lies between rays $a$ and $c$, then:
\mathnote{\Prop
\begin{align}
&\Forall{a \w b \w c}{\Ray}{}\nonumber
\end{align}} 
\begin{enumerate}
\item[a)]
\mathnote{\Prop
\begin{align}
&\Bet{a}{b}{c} \rightarrow \Apart{a}{b} \land \Apart{b}{c} \land \Apart{a}{c} \nonumber
\end{align}}all three rays are divergent;
\item[b)]
\mathnote{\Prop
\begin{align}
&\Bet{a}{b}{c} \rightarrow\Bet{c}{b}{a} \nonumber
\end{align}}the ray $b$ also lies between rays $c$ and $a$;
\item[c)]
\mathnote{\Prop
\begin{align}
&\Bet{a}{b}{c} \rightarrow \lnot \Bet{b}{a}{c} \land \lnot \Bet{a}{c}{b} \nonumber
\end{align}}neither ray $a$, nor ray $c$ can lie between the other two rays.
\end{enumerate}
\end{theorem}

\begin{problem}\Mrg
\mathnote{\Set
\begin{align}
&\Forall{a \w b}{\Ray}{}\nonumber
\end{align}}Given two divergent rays $a$ and $b$. 
\begin{enumerate}
\item[a)]
\mathnote{\Set
\begin{align}
& \Apart{a}{b}\rightarrow  \Build{o}{\Ray}{\Bet{a}{o}{b}}\nonumber
\end{align}}Draw a ray $o$ between the rays $a$ and $b$.
\label{prb:DrawPointBetween}
\item[b)]
\mathnote{\Set
\begin{align}
& \Apart{a}{b}\rightarrow  \Build{c}{\Ray}{\Bet{a}{b}{c}}\nonumber
\end{align}}Draw a ray $c$ such that ray $b$ is between the rays $a$ and $c$.
\end{enumerate}
\end{problem}

It will also be convenient for us to use the four-point order relation for rays.
\mathnote{\Mrg
\begin{align}
&\Obj{a\w b\w c\w d}{\Ray}{} \nonumber\\
&\BetFour{a}{b}{c}{d}\nonumber\\
&\quad\equiv \Bet{a}{b}{c} \land \Bet{a}{b}{d} \land\Bet{a}{c}{d} \land \Bet{b}{c}{d}\nonumber
\end{align}}
We will say that the four given rays are located on a plane in the established order, if any three of them, taken in the same order, satisfy the ray betweenness relation.
It can be proved that ray "betweenness" relation satisfy some sort of transitivity property.

\begin{theorem}\Mrg
\mathnote{\Prop
\begin{align}
&\Forall{a\w b\w c\w d}{\Ray}{}\nonumber\\
& \Bet{a}{b}{c} \rightarrow  \Bet{a}{c}{d} \rightarrow \BetFour{a}{b}{c}{d}\nonumber
\end{align}}If the ray $b$ lies between the rays $a$ and $c$, and the ray $c$ lies between the points $b$ and $d$, then all four rays $a$, $b$, $c$ and $d$ are ordered.
\end{theorem}

\section{5. Flags}

\subsection{5.1. Definition of flags as equivalence classes}

In the previous section we defined rays with the help of the equivalence relation between ordered pairs of distinct points.
Now we can go further and consider similar equivalence relation between ordered pairs of divergent rays.
\mathnote{\Mrg
\begin{align}
&\Obj{X}{\Flag}{}\nonumber\\
&X = \{X_0\bs X_1\} \nonumber\\
&\quad X_0\CC\text{initial ray}\nonumber\\
&\quad X_1\CC\text{terminal ray}\nonumber
\end{align}}
We will say that these pairs of divergent rays represent new geometric figures called {\em flags}.
The point from which the rays originate is called the {\em vertex} of the flag.
The first ray of the pair is called the {\em initial} ray of the flag, while the second ray is called its {\em terminal} ray.

\mathnote{\Mrg
\begin{align}
&\Obj{X\w Y}{\Flag}{}\nonumber\\
&X \approx Y\equiv (X_0 \approx Y_0)\land\SameRay{X_0}{X_1}{Y_1}\nonumber
\end{align}}
Flags themselves as geometric figures are defined as equivalence classes of the following {\em coincidence} relation on the pairs of divergent rays.
Namely, it is said that flags $X$ and $Y$ coincide with each other iff their initial rays $X_0$ and $Y_0$ coincide, and their terminal rays $X_1$ and $Y_1$ both lie on the same side from their common initial ray.  

\mathnote{\Mrg
\begin{align}
&\Obj{X\w Y}{\Flag}{}\nonumber\\
&X\approx Y\Iff {\mathbb F^{\rm +}}(X) = {\mathbb F^{\rm +}}(Y)\nonumber
\end{align}}
Intuitively, a flag is a figure consisting of an initial ray and a half-plane attached to it.
The choice of one of two half-plane is determined by the location of the terminal ray on the plane.

Note that coincidence symbol $\approx$ is overloaded in our notations. 
Depending on its arguments it may refer either to coincidence of rays or to coincidence of flags.

\subsection{5.2. Orientation on the plane}

We already know that every line divides the plane onto two opposite half-planes.
This suggests that on the plane we could consider two opposite {\em orientations} of flags.
We are going to formalize this concept of orientation of flags on the plane in exactly the same way as we previously formalized the concept of the direction of rays on a line.
Namely, we are going to introduce yet another flag equivalence relation, which will be proved to be decidable only with two different equivalence classes.
These equivalence classes will correspond to the {\em left} and {\em right} orientations of flags on the plane.
However, the required equivalence relation for the orientation of the flags becomes somewhat more complicated than a similar relation for the direction of the rays.
We shall define it in three steps.

First, we define a rotation equivalence relation for those flags whose rays come from the same vertex point.
\vskip0.5\baselineskip
\begin{definition}[Rotational equivalence]\Mrg
\label{def:FlagRotationEquivalence}
\mathnote{\Def
\begin{align}
&\Obj{X \w Y}{\Flag}{}\nonumber\\
&X \overset{\text{\tiny\tt rot}}{\sim} Y\equiv (X_0 \approx Y_0\land\SameRay{X_0}{X_1}{Y_1})\lor\nonumber\\
&\s[55]\lor(X_0 \approx \Op{Y}_0\land\SameRay{X_0}{Y_1}{\Op{X}_1})\lor\nonumber\\
&\s[55]\lor \SameRay{X_0}{X_1}{Y_0}\land\SameRay{Y_0}{Y_1}{\Op{X}_0} \lor\nonumber\\
&\s[55]\lor \SameRay{Y_0}{Y_1}{X_0}\land\SameRay{X_0}{X_1}{\Op{Y}_0}\nonumber
\end{align}
}
It is said that that flags $X$ and $Y$ are rotationally equivalent iff either of the four alternatives is valid:
\begin{enumerate}
\item[--] flags $X$ and $Y$ coincide;
\item[--] the flag $X$ coincide with flag $\Op{Y}$ whose both rays are opposite to that of $Y$;
\item[--] the initial ray of flag $Y$ belongs to the flag $X$, while the terminal ray of flag $X$ does not belong to flag $Y$;
\item[--] the initial ray of flag $X$ belongs to the flag $Y$, while the terminal ray of flag $Y$ does not belong to flag $X$.
\end{enumerate}
\end{definition}

\noindent 
This relation of rotational equivalence implies that if we are given two equivalent flags with a common vertex, then we can rotate one of the flags around the vertex so that the flags will eventually coincide with each other. In other words, both the initial rays of the flags and the half-planes attached to them will coincide.

We can prove that this relation of rotational equivalence is decidable, and defines only two classes of equivalent flags.
Thus, for any vertex point on the plane we obtain a local definition of the orientation of the flags.

\begin{marginfigure}[-6\baselineskip]\Mrg
\includegraphics[width=\linewidth]{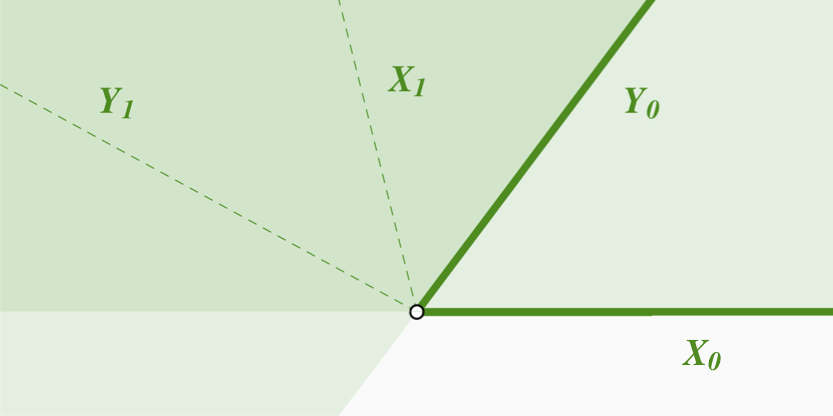}
\caption{Two rotationally equivalent flags $X \overset{\text{\tiny\tt rot}}{\sim} Y$ on the plane.}
\label{fig:DrawRotationalFlags}
\end{marginfigure}

Next, we need to synchronize two local definitions of rotationally equivalent flags at distinct points.
To this end, we define the translational equivalence relation for flags whose initial rays are collinear and have the same direction.
\vskip0.5\baselineskip
\begin{definition}[Translational equivalence]\Mrg
\label{def:FlagTranslationEquivalence}
\mathnote{\Def
\begin{align}
&\Obj{X \w Y}{\Flag}{}\nonumber\\
&X \overset{\text{\tiny \tt tr}}{\sim} Y\equiv (X_0 \sim Y_0)\land\SameHalfplane{X_{00}}{X_{01}}{X_{11}}{Y_{11}}\nonumber
\end{align}}
It is said that that flags $X$ and $Y$ are translationally equivalent iff their initial rays have the same direction on line, while their terminal rays are on the same side of this line.
\end{definition}
\begin{marginfigure}[9.8\baselineskip]\Mrg
\includegraphics[width=\linewidth]{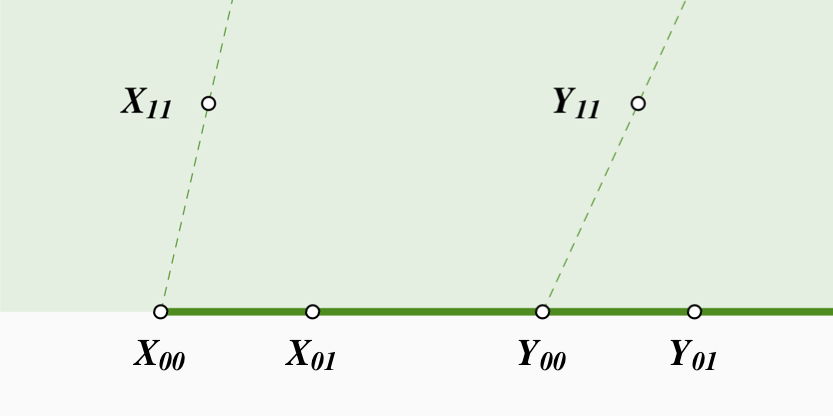}
\caption{Two translationally equivalent flags $X \overset{\text{\tiny \tt tr}}{\sim} Y$ on the plane.}
\label{fig:DrawTranslationalFlags}
\end{marginfigure}
The relation of translational equivalence implies that if we are given two equivalent flags, whose initial rays are collinear and have the same direction, then we can translate one of the flags along their common direction so that it will eventually coincide with the other flag. 
In other words, both the initial rays of the flags and the half-planes attached to them will coincide.

We can prove that this relation of translational equivalence is decidable, and defines only two classes of equivalent flags.
Thus, for any directed line on the plane we obtain yet another local definition of the orientation of the flags.

Finally, we define the universal relation of orientational equivalence of flags as a combination of the previous two equivalences.
\vskip0.5\baselineskip
\begin{definition}[Orientational equivalence]\Mrg
\label{def:FlagOrientationalEquivalence}
\mathnote{\Def
\begin{align}
&\Obj{X \w Y}{\Flag}{}\nonumber\\
&X \sim Y\nonumber\\
&\quad\equiv \Exists{\Prime{X}\w\Prime{Y}}{\Flag}{
\s[4]X\overset{\text{\tiny\tt rot}}{\sim}\Prime{X}
\s[4]\land\s[4]\Prime{X}\overset{\text{\tiny \tt tr}}{\sim}\Prime{Y}
\s[4]\land\s[4]\Prime{Y}\overset{\text{\tiny\tt rot}}{\sim} Y}\nonumber
\end{align}}
It is said that that flags $X$ and $Y$ are equivalent iff there exist two flags $\Prime{X}$ and  $\Prime{Y}$ such that pairs of flags $X$ and $\Prime{X}$ as well as $Y$ and $\Prime{Y}$ are rotationally equivalent while flags $\Prime{X}$ and $\Prime{Y}$ are translationally equivalent.
\end{definition}
The relation of orientational equivalence implies that if we are given two equivalent flags on the plane, then we can rotate them and translate so that they will eventually coincide with each other. 
\mathnote{\Mrg
\begin{align}
&\Obj{X\w Y}{\Flag}{}\nonumber\\
& X\sim Y\Iff{\mathbb O^{\rm +}}(X) = {\mathbb O^{\rm +}}(Y)\nonumber
\end{align}}
We can prove that the relation of orientational equivalence is decidable and there are only two classes of equivalent flags on the plane.
In other words, given any two flags on the plane, we can decide if they have the same or opposite orientation. 
We shall call the equivalence classes of this equivalence relation as {\em orientations} of flags on the plane. 

\mathnote{\Mrg
\begin{align}
&\Obj{X}{\Flag}{}\nonumber\\
&\Flip{X} = \{X_1\bs X_0\} \nonumber\\
&\quad \CC\Flip{X}\text{ is a flip of } X\nonumber
\end{align}}
If in a pair of divergent rays representing a flag, we swap the initial and terminal rays, then the resulting flag will be called a flip of the original one.
The orientation of the flipped flag is opposite to the orientation of the original flag.
The flip of a flag will be denoted by a hat symbol ($\w\Flip{~}\w$) in our formulas on the margins.

\begin{problem}\Mrg
\mathnote{\Set
\begin{align}
&\Forall{X \w Y}{\Flag}{}\nonumber\\
&\Or{X \approx Y}{X \approx \Flip{Y}}\nonumber
\end{align}}Given two flags $X$ and $Y$. Decide whether they have same or opposite orientation.
\label{prb:DrawPointBetween}
\end{problem}

\begin{theorem}\Mrg
\mathnote{\Prop
\begin{align}
&\Forall{X \w Y}{\Flag}{}\nonumber\\
&\lnot\s(\s X \approx Y \land X \approx \Flip{Y}\s)\nonumber
\end{align}} 
Two flags cannot simultaneously have both same and opposite direction.
\end{theorem}

Since we have decidable equivalence relation with just two equivalence classes, we can define a {\em function} that maps every flag to its orientation (boolean value).
To this end, we may first draw three non-collinear points on the plane with the help of Postulates \ref{pst:DrawPoint}, \ref{pst:DrawDistinctPoint}, \ref{pst:DrawExtensionLine} and \ref{pst:DrawPointApartLine}. Taking any of the three points at the vertex, and the other two for the directing points of diverging rays, we can build on them the reference flag, the orientation of which we, by definition, will consider {\tt left}. Opposite orientation we shall call {\tt right}. Since we have just two mutually exclusive orientations, we can treat them as boolean values.
\mathnote{\Mrg
\begin{align}
&{\mathbb O^{\rm +}} :  \Flag \rightarrow\text{Bool}\nonumber
\end{align}}
Finally, given an arbitrary flag, we can compare it with the reference flag and, thus, determine its orientation.

With the help of this function, we can also expand the definition of orientation from flags to triangles.

\section{6. Angles} %

\subsection{6.1. Definition of oriented angles and their measures}

\mathnote{\Mrg
\begin{align}
&\Obj{A}{\Angle}{}\nonumber\\
&A = \{a \cdot b\} \nonumber\\
&\quad a\CC\text{ initial ray}\nonumber\\
&\quad b\CC\text{ terminal ray}\nonumber
\end{align}}
The figure formed by two concurrent rays is called the {\em angle}.
It is important to note that the definition of the angle is wider than a similar definition of flags.
Indeed, in the latter definition, flags are represented by a pair of not just concurrent, but also non-collinear (divergent) rays.
Thus, pair of divergent rays representing a flag also makes some angle, but the reverse is not true.

Just like for flags, the point from which the rays originate is called the {\em vertex} of this angle, 
the first ray of the pair is called the {\em initial} ray, and the second ray is called the {\em terminal} ray of the angle.
\mathnote{\Mrg
\begin{align}
&\Obj{A\w B}{\Angle}{}\nonumber\\
&A\approx B\nonumber\\
&\quad\CC\text{ angles $A$ and $B$ are congruent}\nonumber\\
&{\mathbb M^{\rm +}}(A) \nonumber\\
&\quad\CC\text{ equivalence class for angle $A$}\nonumber\\[0.5\baselineskip]
&A\approx B\Iff {\mathbb M^{\rm +}}(A) = {\mathbb M^{\rm +}}(B)\nonumber
\end{align}}
The concurrency of rays is indicated by the dot symbol ($\cdot$) in our formulas on the margins.

As for the definition of length, we could assume, following Hilbert, that there is some equivalence (congruence) relation between angles. This equivalence relation must be reflexive, symmetric, and transitive. Then, based on this equivalence relation, we can define angle {\em measures} as the corresponding classes of congruent angles.

\mathnote{\Mrg
\begin{align}
&{\mathbb M^{\rm +}} (A : \Angle) : \Measure\nonumber\\
&{\mathbb M^{\rm -}} (\alpha : \Measure) : \Build{A}{\Angle}{{\mathbb M^{\rm +}}(A) = \alpha}\nonumber
\end{align}}
Alternatively, we can immediately introduce the {\em measure} of an angle as the basic undefined notion and assume that there is a mapping of angles (pairs of concurrent rays) to their measures and back. 
In other words, each angle can be assigned its measure, and for each measure a corresponding angle can be found. 
Such mappings also define an equivalence relation for angles.

\mathnote{
\begin{align}
&\Obj{a \w b \w c \w d}{\Ray},\nonumber\\[0.5\baselineskip]
&\AngRs{a}{b} = {\mathbb M^{\rm +}} (\{a \cdot b\})\nonumber\\
&\quad \CC \text{  measure of angle ${\{a \cdot b\}}$}\nonumber\\[0.5\baselineskip]
&\AngCongRs{a}{b}{c}{d}\nonumber\\
&\quad \CC \text{  measures of angles ${\{a \cdot b\}}$ and ${\{c \cdot d\}}$ are equal}\nonumber
\end{align}}
We shall denote the angle measure for the pair of concurrent rays by a double angular brackets. 

In order to complete the program of full arithmetization of synthetic geometry, we will need to prove later that the angular measures, like the lengths of segments, constitute a constructive field. In other words, the angular measures can be compared with each other, added to each other (modulo the full angle) and even multiplied with each other (not in this paper) satisfying all necessary properties of a constructive filed.

\subsection{6.2. Axioms of oriented angles}

Now we are going to formulate axioms that explore the meaning of the concept of angular measure, and also associate it with the concept of flag orientation.

\mathnote{\Mrg
\begin{align}
&\Forall{a}{\Ray}{}\nonumber\\
&\AngRs{a}{a} = 0^\circ\nonumber
\end{align}}An angle is called {\em null} if its initial and terminal rays coincide with each other. 
The measure of any null angle will be denoted by the symbol "$0^\circ$". 
\vskip0.5\baselineskip
\begin{axiom}[Null angles]\Mrg
\mathnote{\Prop
\begin{align}
&\Forall{a \w b}{\Ray}{}\nonumber\\
&\AngCongRs{a}{a}{b}{b} \nonumber
\end{align}}All null angles on the plane are equal to each other.
\label{axm:EqRs_EqRs_EqAs}
\end{axiom}

\mathnote{\Mrg
\begin{align}
&\Forall{a}{\Ray}{}\nonumber\\
&\AngRs{a}{\Op{a}} = 180^\circ\nonumber
\end{align}}An angle is called {\em straight} if its initial and terminal rays are opposite, i.e. they form complete line. 
The measure of any such angle will be denoted by the symbol "$180^\circ$". 
\vskip0.5\baselineskip
\begin{axiom}[Straight angles]\Mrg
\mathnote{\Prop
\begin{align}
&\Forall{a \w b}{\Ray}{}\nonumber\\
&\AngCongRs{a}{\Op{a}}{b}{\Op{b}} \nonumber
\end{align}}All straight angles on the plane are equal to each other. 
\label{axm:OpRs_OpRs_EqAs}
\end{axiom}
Let us commit the position of the initial ray. 
Intuitively, we assume that as the terminal ray rotates, the measure of the angle formed by these two rays will change  monotonously. 
In other words, there are no two distinct positions of the terminal ray for which the measures of the corresponding angles would be equal.
\vskip0.5\baselineskip
\begin{axiom}[Unique angles]\Mrg
\mathnote{\Prop
\begin{align}
&\Forall{a \w b \w c}{\Ray}{}\nonumber\\
&\AngCongRs{a}{b}{a}{c} \Iff b \approx c \nonumber
\end{align}} 
Given two angles with the same initial ray. 
Then the terminal rays of these angles coincide iff their angle measures are equal.
\label{axm:EqAs_EqRs}
\end{axiom}

Informally speaking, at each point of the plane, taken as a vertex, we can define its own system of angles.
Then we have to make sure that for each point of the plane and for each such system of angles, all the "full angles" are equal.

\vskip0.5\baselineskip
\begin{axiom}[Explementary angles]\Mrg
\mathnote{\Prop
\begin{align}
&\Forall{a \w b \w c \w d}{\Ray}{}\nonumber\\
&\AngCongRs{a}{b}{c}{d} \rightarrow \AngCongRs{b}{a}{d}{c} \nonumber
\end{align}} If two angles are equal, then their explementary angles are also equal. 
\label{axm:EqAs_EqExpAs}
\end{axiom}

We also have to make sure that for every point of the plane considered as the vertex of the local system of angles, all angles are measured in the same direction.
In other words angle measures and their orientations should be consistent with each other.

\vskip0.5\baselineskip
\begin{axiom}[Angle orientation]\Mrg
\mathnote{\Prop
\begin{align}
&\Forall{X \w Y}{\Flag}{}\nonumber\\
&\AngCongRs{X_0}{X_1}{Y_0}{Y_1} \rightarrow {\mathbb O^{\rm +}}(X) = {\mathbb O^{\rm +}}(Y)\nonumber
\end{align}}
Given two flags. If angular measures of corresponding angles are equal, then orientations of these flags should also be equal. 
\label{axm:EqAs_EqTs}
\end{axiom}
With this axiom we can lift the definition of orientation up from flags (pairs of rays) to angle measures ${\mathbb O^{\rm +}}(\alpha)$, provided that angle measure $\alpha$ is not equal to $0^\circ$ or $180^\circ$. 

This axiom completes our definition of {\em oriented} angles on the plane. 
In the next paragraph we define {\em non-oriented} or {\em convex} angles that will be necessary to formulate our other axioms of Euclidean geometry.

\subsection{6.3. Definition of non-oriented angles}

The sides of the angle divide the whole plane into two regions. 
One of them is called the inner region of the angle, and the other is its outer region.
Usually, the inner region is considered to be the one in which the segment with endpoints on the sides of the angle is placed entirely.
We will also call this angle {\em convex}. The other part of the plane is said to form {\em reflex} angle.

Since we defined orientation of all flags using the reference flag, we can assume, without loss of generality, that the oriented angle of the reference flag corresponds to the convex region on the plane. This means that all angles with {\tt left} orientation will be considered convex, while those with {\tt right} orientation will be considered reflex. Finally, with the help of Axiom \ref{axm:EqAs_EqTs} we can lift the definition of convex and reflex angles up to the angle measures, provided that the angle measures under consideration are not equal to $0^\circ$ or $180^\circ$.
\vskip0.5\baselineskip
\begin{definition}\Mrg
\label{def:FlagEquivalence}
\mathnote{\Def
\begin{align}
&\Obj{\alpha}{\Measure},\nonumber\\
&{\rm Convex~}\alpha \equiv \alpha \neq 0^\circ \land\alpha \neq 180^\circ\land  {\mathbb O^{\rm +}}(\alpha) = \text{\tt left}\nonumber\\
&{\rm Reflex~}\alpha \equiv \alpha \neq 0^\circ \land\alpha \neq 180^\circ\land  {\mathbb O^{\rm +}}(\alpha) = \text{\tt right}\nonumber
\end{align}} 
Consider non-null and and non-straight angles. 
They are called convex iff orientation of their rays is {\tt left}. 
They are called reflex iff orientation of their rays is {\tt right}.
\end{definition}

Whatever the flag, it has its own orientation and the corresponding angular measure.
Let us now look at the flipped flag. It has the opposite orientation and explementary angular measure.
If the measure of the original flag is reflex, then the measure of the flipped flag will be convex and vice versa.
\mathnote{\Mrg
\begin{align}
&\Obj{a \w b}{\Ray}{}\nonumber\\
& \IntAngRs{a}{b}\equiv{\tt if~ }  {\mathbb O^{\rm +}} (\{a \bs b\})\nonumber\\
&\s[75]{\tt ~then~ }\AngRs{a}{b} \nonumber\\
&\s[75]{\tt ~else~ }\AngRs{b}{a}\nonumber
\end{align}}
Thus, for any flag, we can define corresponding {\em non-oriented} angular measure as follows.
If orientation of the original flag is {\tt left}, then non-oriented angular measure, by definition, coincides with the oriented one.
Otherwise, when orientation of the flag is {\tt right}, non-oriented angular measure is defined to coincide with the oriented measure of the flipped flag.

\mathnote{\Mrg
\begin{align}
&\Obj{O \w A \w B}{\Point}{}\nonumber\\
& \IntAngPs{A}{O}{B}\equiv\IntAngRs{\{O\bs A\}}{\{O\bs B\}}\nonumber
\end{align}}
Given three non-collinear points $O$, $A$ and $B$ we can define a flag with vertex at the points $O$, initial ray $OA$ and terminal ray $OB$.
To this flag we can associate the corresponding non-oriented angular measure which we shall denote $\IntAngPs{A}{O}{B}$.
Thus, whenever we use such a non-oriented angular measure for three points, we assume that these three points are non-collinear.

\section{7. Segment arithmetic} %

In this section we will continue our discussion of the properties of circles, and the main question for us will be the following. Under what conditions do two circles intersect or, conversely, do not intersect?

\subsection{7.1. Concentric circles}

\begin{marginfigure}[0.5\baselineskip]\Mrg
\includegraphics[width=\linewidth]{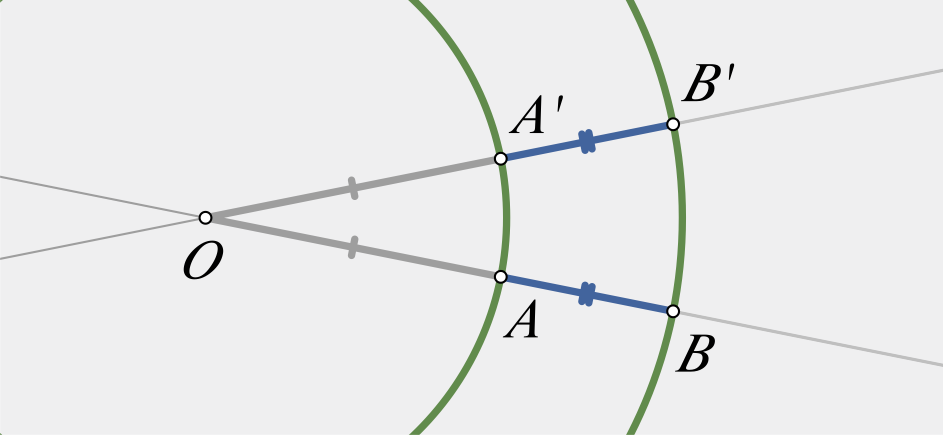}
\caption{The intersection of concurrent rays and  concentric circles.}
\label{fig:raysandcircles}
\end{marginfigure}

The simplest arrangement of circles on a plane in which they obviously do not intersect is when the circles have the same center and different radii. 
Such circles are called {\em concentric}. 
Speaking informally, such circles are akin to parallel lines in the sense that they do not intersect each other and are equally spaced from each other.
In other words, if some point of one of the concentric circles lies inside the other, then the entire first circle will lie inside the second (see Fig. \ref{fig:raysandcircles}). 

This property of concentric circles can be formulated a little differently. 
Take two concentric circles and consider a pair of intersection points of these circles with the rays originating from the center of these circles. 
Then the order in which each ray meets with concentric circles is the same. 
Moreover, the lengths of the segments formed by the intersection points of each ray with concentric circles are equal.  
This property of concentric circles can be formalized as the following axiom:

\vskip0.5\baselineskip
\begin{axiom}[Concentric circles and rays]\Mrg
\mathnote{\Prop
\begin{align}
&\Forall{O \w A \w \Prime{A} \w B \w \Prime{B}}{\Point}{}\nonumber\\
& \Bet{O}{A}{B} \rightarrow \SameRay{O}{\Prime{A}}{\Prime{B}} \rightarrow\nonumber\\
&\quad \rightarrow\Lcong{O}{A}{O}{\Prime{A}} \rightarrow \Lcong{O}{B}{O}{\Prime{B}}  \rightarrow \nonumber\\
&\quad \rightarrow \Bet{O}{\Prime{A}}{\Prime{B}}\land\Lcong{A}{B}{\Prime{A}}{\Prime{B}}
\nonumber
\end{align}}
Consider two concentric circles and two rays originating from the center of the circles $O$.
Let points $A$ and $\Prime{A}$ lie on one of the circles, and points $B$ and $\Prime {B}$ on another.
If in this case the point $A$ lies on one of the rays between the points $O$ and $B$, and points $\Prime{A}$ and $\Prime{B}$ on another ray, then the point $\Prime{A}$ lies between the points $O$ and $\Prime{B}$ and the lengths of the segments $\Seg{A}{B}$ and $\Seg{\Prime{A}}{\Prime{B}}$ are equal.
\label{axm:CircleSegmentSubtraction}
\end{axiom}

\subsection{7.2. Segment transfer}

The most general condition for the intersection of two circles can be formulated as follows: two circles intersect with each other iff one point of the first circle lies inside the second, and vice versa.

\begin{marginfigure}[4.5\baselineskip]\Mrg
\includegraphics[width=\linewidth]{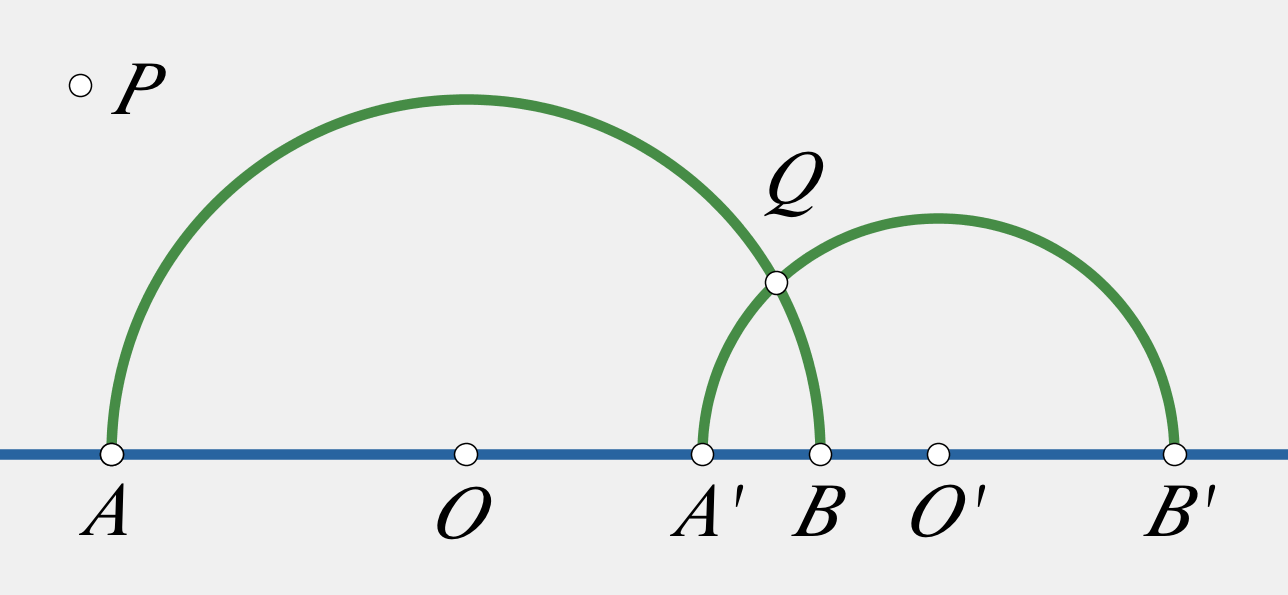}
\caption{The intersection of two circles.}
\label{fig:TwoCirclesIntersection}
\end{marginfigure}
In order to simplify and formalize this condition, let's look at Fig. \ref{fig:TwoCirclesIntersection}. We see two circles, with centers at distinct points $O$ and $\Prime{O}$. We can draw a unique line through the points. This line intersects every circle at two points. Let the intersection points of the line with the first circle be denoted $A$ and $B$, and with the second one --- $\Prime{A}$ and $\Prime{B}$. Two circles intersects if the order of these points on a straight line can always be chosen so that the point $\Prime{A}$ lies between the points $A$ and $B$, and the point $B$ --- between the points $\Prime{A}$ and $\Prime{B}$. We can formalize this as our final postulate.

\vskip0.5\baselineskip
\begin{postulate}[Circle-Circle intersection]\Mrg
\mathnote{\Set
\begin{align}
&\Forall{O \w A \w B \w \Prime{O} \w \Prime{A} \w \Prime{B} \w P}{\Point}{}\nonumber\\
&\lnot\ColThree{O}{\Prime{O}}{P} \rightarrow\BetFour{A}{\Prime{A}}{B}{\Prime{B}}\rightarrow\nonumber\\
&\quad \rightarrow \Middle{A}{O}{B} \rightarrow  \Middle{\Prime{A}\s[-2]}{\Prime{O}\s[-2]}{\Prime{B}}\rightarrow\nonumber\\
&\quad \rightarrow\{\s \Obj{Q}{\Point} \vl \SameHalfplane{O}{\Prime{O}}{P}{Q} \land 
\nonumber\\
&\quad\qquad\qquad \land \Lcong{O}{Q}{O}{A} \land \Lcong{\Prime{O}}{Q}{\Prime{O}}{\Prime{A}} \s\}\nonumber
\end{align}}
Given two circles and a line passing through the centers of these circles and intersecting each of them at two points.
Let one of the two points of intersection of the first circle lies between the two points of intersection of the second circle and vice versa.
Draw a point of intersection of these two circles with each other, which will lie on a given side of this line.
\label{pst:DrawIntersectionOfTwoCircles}
\end{postulate}

With  the help of this postulate we can solve the following problems from Euclid's {\em Elements}.

\begin{problem}\Mrg
\mathnote{\Set
\begin{align}
&\Forall{A\w B\w P}{\Point}{}\nonumber\\
&\lnot\ColThree{A}{B}{P}\rightarrow  \nonumber\\
&\quad\rightarrow\Build{C}{\Point}{\SameHalfplane{A}{B}{C}{P} \land \LcongThree{A}{B}{B}{C}{A}{C}}\nonumber
\end{align}}Given a segment on the border of a half-plane. Build an equilateral triangle on a given segment in the given half-plane.\\
 $~$\hfill{\Set\rm \scriptsize / Euclid, Book I : Proposition 1}
\label{prb:DrawEquilateralTriangle}
\end{problem}

Finally, we can solve the problem of transferring a segment from one place on the plane to another. 
In other words, regardless of the segment, we can draw an equal segment on any given ray.

\begin{problem}\Mrg
\mathnote{\Set
\begin{align}
&\Forall{O\w A\w B\w C}{\Point}{}\nonumber\\
&O\neq A\rightarrow\Build{D}{\Point}{\SameRayO{O}{A}{D}\land \Lcong{O}{D}{B}{C}}\nonumber
\end{align}}
Given the segment and the ray. Draw on the ray from its origin a segment equal to the given one.\\
 $~$\hfill{\Set\rm \scriptsize / Euclid, Book I : Proposition 2}
\label{prb:DrawSegmentOnRay}
\end{problem}

\subsection{7.3. Segment order}

Using the betweenness relation, we can define the strict order relation for the lengths of segments.
\vskip0.5\baselineskip
\begin{definition}\Mrg
\mathnote{\Def
\begin{align}
&\Obj{a \w b}{\Length},\nonumber\\
&a < b\nonumber\\
&\quad\equiv\Exists{A \w B\w C}{\Point}{\LC{A}{B}{a}\land \LC{A}{C}{b} \land\BetOX{A}{B}{C}}\nonumber\\
&\quad\CC\text{\rm  segment length $a$ is less then length $b$}\nonumber
\end{align}} 
It is said that the length of $a$ is less than the length of $b$, if there are three points $A$, $B$ and $C$ on the plane such that the length of the segment $AB$ is $a$, the length of the segment $AC$ is $b$, and, moreover, the point $B$ lies between the points $A$ and $C$ or coincides with the point $A$.
\end{definition}

Based on this definition, we can prove all the necessary basic properties of a strict order relation.

\vskip0.5\baselineskip
\begin{theorem}[Non-negative length]\Mrg
\mathnote{\Prop
\begin{align}
&\Forall{a}{\Length}{}\nonumber\\
& \lnot (a < 0)\nonumber
\end{align}
} 
The segment length can not be less than the length of the null segment.
\end{theorem}

\begin{theorem}[Segment trichotomy]\Mrg
\mathnote{\Prop
\begin{align}
&\Forall{a \w b}{\Length}{}\nonumber\\
& (a < b) \lor (a = b) \lor (b < a)\nonumber
\end{align}} 
Whatever the two lengths of the segments, they are either equal to each other, or one of them is less than the other.
\end{theorem}

\vskip0.5\baselineskip
\begin{theorem}[Properties of strict order]\Mrg
\mathnote{\Prop
\begin{align}
&\Forall{a\w b\w c}{\Length}{}\nonumber
\end{align}}
The relation "less then" defines a  strict order on the lengths of segment.
\begin{description}
\item[\Prop\textsc{irreflexivity.}]
\mathnote{\Prop
\begin{align}
& \lnot (a < a)\nonumber
\end{align}}The length of any segment is not less than itself.
\item[\Prop\textsc{transitivity.}]
\mathnote{\Prop
\begin{align}
& (a < b) \rightarrow (b < c) \rightarrow (a < c) \nonumber
\end{align}}Three segments of length $a$, $b$ and $c$ are given. 
If $a$ is less than $b$ and $b$ is less than $c$, then $a$ is less than $c$.
\end{description}
\label{thm:StrictOrder}
\end{theorem}

\begin{problem}\Mrg
\mathnote{\Set
\begin{align}
&\Forall{a \w b \w c}{\Length}{} \nonumber\\
&a < b\rightarrow\Or{a < c}{c < b}\nonumber
\end{align}}Let three lengths of segments $a$, $b$, and $c$ be given, so that $a$ is less than $b$. 
Decide whether the length of the segment $c$ is more than $a$, or less than $b$.
\end{problem}

\subsection{7.4. Segment addition}

With the help of our postulates, we can solve the following problem, which gives us the ability to add segments to each other.

\begin{problem}[Segment addition]\Mrg
\mathnote{\Set
\begin{align}
&\Forall{a \w b}{\Length}{} \nonumber\\
&\Build{c}{\Length}{\Forall{A \w B \w C}{\Point}{\BetOO{A}{B}{C} \rightarrow\nonumber\\
&\quad\quad\quad\quad\quad\rightarrow\LC{A}{B}{a}\rightarrow\LC{B}{C}{b}\rightarrow\LC{A}{C}{c}}}\nonumber\\[0.5\baselineskip]
&\quad\CC\text{\rm  the resulting length $c$ is called the sum of two}\nonumber\\
&\quad\text{\rm  ~~~~lengths $a$ and $b$, i.e. we denote $c = a + b$}\nonumber
\end{align}}
Given two segments of length $a$ and $b$.
Draw another segment of length $c$ with the following property. 
Whatever are three points $A$, $B$ and $C$ such that the point $B$ belongs to the segment $AC$, if the length of segment $AB$ is equal to $a$ and the length of segment $BC$ is equal to $b$, then the length of segment $AC$ is also equal to $c$.
\end{problem}

Now we can prove all the necessary basic properties of the segment addition operation.

\begin{theorem}[Properties of addition]
Addition of segments has the following properties:
\mathnote{\Prop
\begin{align}
&\Forall{a\w b\w c}{\Length}{}\nonumber
\end{align}
} 
\begin{description}
\item[\Prop\textsc{null length.}]
\mathnote{\Prop
\begin{align}
& a + 0 = a\nonumber
\end{align}}When adding null to any length, the quantity does not change.
\item[\Prop\textsc{commutativity.}]
\mathnote{\Prop
\begin{align}
& a + b = b + a \nonumber
\end{align}}Addition is commutative --- one can change the order of the terms in a sum, and the result is the same.
\item[\Prop\textsc{associativity.}]
\mathnote{\Prop
\begin{align}
& a + (b + c) = (a + b) + c \nonumber
\end{align}}Addition is associative --- when adding three or more lengths, the order of operations does not matter.
\item[\Prop\textsc{subtraction.}]
\mathnote{\Prop
\begin{align}
& a + c = b + c \s[4]\rightarrow\s[4] a = b \nonumber
\end{align}}One can subtract the same length from two equal sums, and the remaining lengths will also be equal.
\end{description}
\end{theorem}

\section{8. Congruence of triangles} %

Triangle is figure on the plane that is completely determined by three non-collinear points $A$, $B$ and $C$, called {\em vertices} of the triangle.
These vertices form segments $AB$, $BC$ and $CA$, which are called {\em sides} of the triangle. 
These segments, connected at the vertices, form the angles $ABC$, $BCA$ and $CAB$.
Vertexes, sides and angles together form the {\em triangle} which is denoted as $\Triangle{A}{B}{C}$.

We can consider the following triangles equivalence relation.
\vskip0.5\baselineskip
\begin{definition}\Mrg
\mathnote{\Def
\begin{align}
&A\w B\w C \w D\w E \w F : \Point,\nonumber\\
&\CongTriangles{A}{B}{C}{D}{E}{F} \nonumber\\
&\quad \equiv \Lcong{A}{B}{D}{E}\land \IntAngCongPs{A}{B}{C}{D}{E}{F} \land \nonumber\\
&\quad\quad\s[6] \Lcong{B}{C}{E}{F} \land \IntAngCongPs{B}{C}{A}{E}{F}{D}\land\nonumber\\
&\quad\quad\s[6] \Lcong{C}{A}{F}{D} \land \IntAngCongPs{C}{A}{B}{F}{D}{E} \nonumber\\
&\quad\CC\text{\rm the triangles $\Triangle{A}{B}{C}$ and $\Triangle{D}{E}{F}$ are equal} \nonumber
\end{align}}
Two triangles are called congruent iff lengths of their respective sides and (non-oriented) measures of their respective angles are equal.\label{def:SuperpositionPrinciple}
\end{definition}

It is easy to prove that this relation has all the necessary properties --- reflexivity, symmetry and transitivity. 
Thus, we can talk about classes of equivalent triangles that have the same side lengths and angle measures.
\subsection{8.1. Method of superposition}

To prove the theorems of congruence of triangles, Euclid used the so-called superposition method.
Namely, Euclid suggested that one triangle could be moved and superimposed on another, so that their equal sides and angles coincide.
Then he argued that all the other sides and angles of these triangles should also coincide.

Several attempts were made to formalize the Euclidean superposition method. 
Some authors considered superpositions as isometric mappings of one geometrical figure onto another, and then studied the properties of these mappings. 
Another approach, which we also adhere to, is used, for example, in the popular in Russia textbook on geometry by A.V. Pogorelov \cite{Pogorelov1974}.
In this case, the superposition method is considered as a weak statement only about the mere {\em existence} of an equal triangle in any part of the plane. 
Later we will be able to prove a constructive version of the axiom of superposition and build an equal triangle anywhere on the plane using our standard tools -- a pencil, a straightedge and a compass.

This latter approach is closest in spirit to the original understanding of the method of superposition by Euclid himself and can be formalized as the following axiom:

\vskip0.5\baselineskip
\begin{axiom}[Triangle superposition]\Mrg
\mathnote{\Prop
\begin{align}
&\Forall{A \w B \w C \w D \w E \w F}{\Point}{}\nonumber\\
&\lnot \ColThree{A}{B}{C} \rightarrow \lnot \ColThree{D}{E}{F} \rightarrow\Lcong{A}{B}{D}{E}\rightarrow\nonumber\\
&\quad\rightarrow \Exists{F'}{\Point}{\w \CongTriangles{A}{B}{C}{D}{E}{\Prime{F}} \land \SameHalfplane{D}{E}{F}{\Prime{F}}}\nonumber
\end{align}}
Given two triangles $\Triangle{A}{B}{C}$ and $\Triangle{D}{E}{F}$ whose bases  $\Len{A}{B}$ and $\Len{D}{E}$ are equal, then there exists a point $\Prime{F}$ located on the same side of the base as point $F$, so that triangles  $\Triangle{D}{E}{\Prime{F}}$ and $\Triangle{A}{B}{C}$ are equal.
\label{axm:SuperpositionPrinciple}
\end{axiom}

Once again, we note that this axiom does not provide a procedure for how to actually construct a congruent triangle in a given place.
The simple existence of a congruent triangle is sufficient to prove all congruence theorems (SAS, ASA and SSS) and to provide the prerequisites for the applicability of the necessary postulates for the further construction of an equal triangle.

\begin{theorem}[SAS]\Mrg
\mathnote{\Prop
\begin{align}
&\Forall{A \w B \w C \w D \w E \w F}{\Point}{}\nonumber\\
& \Lcong{A}{B}{D}{E}\rightarrow \Lcong{B}{C}{E}{F}\rightarrow\IntAngCongPs{A}{B}{C}{D}{E}{F} \rightarrow \nonumber\\
&\quad \rightarrow  \CongTriangles{A}{B}{C}{D}{E}{F}\nonumber
\end{align}} 
If the two sides and the angle between them of one triangle are equal respectively to the two sides and the angle between them of the other triangle, then such triangles are equal. \\
 $~$\hfill{\Prop\rm \scriptsize / Euclid, Book I : Proposition 4 / Hilbert, Chapter 1 : Theorem 12}
\end{theorem}

\begin{theorem}[ASA]\Mrg
\mathnote{\Prop
\begin{align}
&\Forall{A \w B \w C \w D \w E \w F}{\Point}{}\nonumber\\
& \lnot\ColThree{A}{B}{C}\rightarrow\lnot\ColThree{D}{E}{F}\rightarrow \Lcong{A}{C}{D}{F}\rightarrow\nonumber\\
&\quad \rightarrow \IntAngCongPs{B}{C}{A}{E}{F}{D}\rightarrow\IntAngCongPs{C}{A}{B}{F}{D}{E}\rightarrow\nonumber\\
&\quad \rightarrow  \CongTriangles{A}{B}{C}{D}{E}{F}\nonumber
\end{align}} 
If the side and two adjacent angles of one triangle are equal respectively to the side and two adjacent angles of the other triangle, then such triangles are equal. \\
 $~$\hfill{\Prop\rm \scriptsize / Euclid, Book I : Proposition 26\,a / Hilbert, Chapter 1 : Theorem 13}
\end{theorem}

\begin{theorem}[SSS]\Mrg
\mathnote{\Prop
\begin{align}
&\Forall{A \w B \w C \w D \w E \w F}{\Point}{}\nonumber\\
& \lnot\ColThree{A}{B}{C}\rightarrow\lnot\ColThree{D}{E}{F}\rightarrow\Lcong{A}{B}{D}{E}\rightarrow\nonumber\\
&\quad \rightarrow \Lcong{B}{C}{E}{F}\rightarrow\Lcong{C}{A}{F}{D}\rightarrow \nonumber\\
&\quad \rightarrow  \CongTriangles{A}{B}{C}{D}{E}{F}\nonumber
\end{align}} 
If three sides of one triangle are equal to three sides of another triangle, then such triangles are equal. \\
 $~$\hfill{\Prop\rm \scriptsize / Euclid, Book I : Proposition 8 / Hilbert, Chapter 1 : Theorem 18}
\end{theorem}

Thus, we have proved that the triangle is a rigid figure, and all its angles are completely determined by its sides.

\subsection{8.2. Triangle congruence}
In the previous section, we started with the existential axiom, which formalizes the Euclidean method of superposition, and on the basis of this axiom we proved all the main theorems on the congruence of triangles. However, we could also choose an alternative approach to the presentation of the material, which is in some sense opposite to the previous one. Namely, we can first accept the following axiom, which is essentially equivalent to SAS and SSS congruence theorems.
\vskip0.5\baselineskip
\begin{axiombis}[Triangle congruence]\Mrg
\mathnote{\Prop
\begin{align}
&\Forall{A \w B \w C \w D \w E \w F}{\Point}{}\nonumber\\
&\lnot \ColThree{A}{B}{C} \rightarrow \lnot \ColThree{D}{E}{F} \rightarrow \nonumber\\
&\quad\rightarrow\Lcong{A}{B}{D}{E}\rightarrow\Lcong{B}{C}{E}{F}\rightarrow\nonumber\\
&\quad\rightarrow\Lcong{C}{A}{F}{D}\Iff \IntAngCongPs{A}{B}{C}{D}{E}{F}\nonumber
\end{align}}
Let two sides of one triangle be equal respectively to two sides of another triangle.
Then the third sides of these triangles will be equal, iff the angles opposite to them are equal.
\label{axm:CongruencePrinciple}
\end{axiombis}
Proceeding from this axiom, we prove that for any triangle we can draw an equal triangle on a given side of a given semi-axis.
In other words, from the axiom of the congruence of triangles, we can obtain a constructive analogue of the Euclidean superposition method.

It is interesting to note that Hilbert rejected the "superposition principle" of Euclid, but left himself the opportunity to transfer angles and impose them onto each other.
To do this, Hilbert had to accept the first theorem on the congruence of triangles (SAS) as a new axiom.
From our formalization of Euclidean. geometry, we see that if we go further and also give up the possibility of transferring and overlapping angles, then we also have to accept the third triangle congruence theorem (SSS) as a new axiom.
In other words, the Euclidean method of superposition is exactly equivalent to two (SAS and SSS) theorems on the congruence of triangles.

First, we prove so-called {\em triangle inequality} theorem.
\begin{theorem}[Triangle inequality]\Mrg
\mathnote{\Prop
\begin{align}
&\Forall{A \w B \w C}{\Point}{}\nonumber\\
&\lnot \ColThree{A}{B}{C} \rightarrow \LessPlus{A}{C}{A}{B}{B}{C}\nonumber
\end{align}
}Any side of the triangle is less than the sum of its two other sides.\\
$~$\hfill{\Prop\rm \scriptsize / Euclid, Book I : Proposition 20}
\end{theorem}

Nex, we prove the statement, in a sense, the opposite of our Postulate \ref{pst:DrawIntersectionOfTwoCircles} about the intersection of two circles. This theorem states that if we are given two intersecting circles, then a straight line passing through the centers of these circles will intersect each of the circles at exactly two points. These two pairs of points will be located on the same line, so that one of the two points of the first circle will be located between two points of the second circle and vice versa.
\begin{theorem}\Mrg
\mathnote{\Prop
\begin{align}
&\Forall{O \w A \w B \w P \w \Prime{O} \s \Prime{A} \s \Prime{B}}{\Point}{}\nonumber\\
&\lnot \ColThree{O}{\Prime{O}}{P}\rightarrow \Middle{A}{O}{B} \rightarrow \Middle{\Prime{A}}{\Prime{O}}{\Prime{B}} \rightarrow\nonumber\\
&\quad \rightarrow \Lcong{O}{A}{O}{P}\rightarrow \Lcong{\Prime{O}}{\Prime{A}}{\Prime{O}}{P}\rightarrow\nonumber\\
&\quad \rightarrow \BetFour{A}{O}{\Prime{O}}{\Prime{B}}\rightarrow\BetFour{A}{\Prime{A}}{B}{\Prime{B}}\nonumber
\end{align}} 
Given two circles and a straight line passing through the centers of these circles. 
Then, if the circles intersect each other, then two pair of points of intersection of the circles with the line are ordered in such a way that one of the two points of intersection of the first circle with the line is between the two points of intersection of the second circle with this line and vice versa.
\end{theorem}

Finally, we prove the constructive superposition theorem. The general idea of the proof is as follows. Consider an arbitrary triangle. We build two circles with centers at two vertices that form the base of the triangle and pass through the third vertex of this triangle. These two circles intersect with each other at the third vertex. They also cross the straight line passing through the centers of the two circles. Consequently, by the theorem above, the intersection points of these circles with the straight line will be ordered.

Now, using the theorem on the transfer of segments, we can transfer all points from this line to any other. 
In this case, all points on this other line will be arranged in exactly the same order as on the initial line. 
Finally, using the postulate about the intersection of two circles, we can build the intersection point in the given half-plane. 
This intersection point is nothing but the third vertex of a required triangle which is equal to the given one.

\begin{problem}[Constructive superposition]\Mrg
\mathnote{\Set
\begin{align}
&\Forall{A \w B \w C\w D\w E\w F}{\Point}{}\nonumber\\
&\lnot \ColThree{A}{B}{C} \rightarrow \lnot \ColThree{D}{E}{F} \rightarrow\Lcong{A}{B}{D}{E}\rightarrow \nonumber\\
&\quad\rightarrow  \Build{\Prime{F}}{\Point}{\w \CongTriangles{A}{B}{C}{D}{E}{\Prime{F}} \land \SameHalfplane{D}{E}{F}{\Prime{F}}}\nonumber
\end{align}}Given two triangles $\Triangle{A}{B}{C}$ and $\Triangle{D}{E}{F}$ whose bases  $\Len{A}{B}$ and $\Len{D}{E}$ are equal.
Draw a point $\Prime{F}$ located on the same side of the base as point $F$, so that triangles  $\Triangle{D}{E}{\Prime{F}}$ and $\Triangle{A}{B}{C}$ are equal.
\end{problem}

Thus we have proved equivalence of the superposition axiom and the axiom of triangle congruence.

\section{9. Angle arithmetic} %

\subsection{9.1. Angle transfer}

Using a constructive variant of the superposition method, we can transfer triangles from one location on a plane to another.
With the help of triangles, we can also transfer convex angles.

\begin{problem}\Mrg
\mathnote{\Set
\begin{align}
&\Forall{A \w O \w B\w \Prime{A}\w \Prime{O}\w \Prime{B}}{\Point}{}\nonumber\\
&\lnot \ColThree{A}{O}{B} \rightarrow\lnot \ColThree{\Prime{A}}{\Prime{O}}{\Prime{B}} \rightarrow\nonumber\\
&\quad\rightarrow  \Build{P}{\Point}{\SameHalfplane{O}{A}{B}{P} \land \IntAngCongPs{A}{O}{P}{\Prime{A}}{\Prime{O}}{\Prime{B}}}\nonumber
\end{align}}
 At a given point on a given straight line and to a given half-plane, to construct a convex angle equal to a given angle.\\
$~$\hfill{\Set\rm \scriptsize / Euclid, Book I : Proposition 23 / Hilbert, Chapter 1 : Axiom $\rm III_4$}
\end{problem}

Then we can extend this construction to draw an equal oriented angle from a given ray.

\begin{problem}[Draw angle clockwise]\Mrg
\mathnote{\Set
\begin{align}
&\ForallDouble{a}{\Ray}{\alpha}{\Measure}\nonumber\\
&\Build{b}{\Ray}{ \AC{a}{b}{\alpha}}\nonumber
\end{align}}
Given an arbitrary angle and an arbitrary ray, draw an equal angle from the ray with given orientation.
\end{problem}

\subsection{9.2. Angle order}

Using the rays betweenness relation, we can define the strict order relation for the angular measures.

\vskip0.5\baselineskip
\begin{definition}\Mrg
\mathnote{\Def
\begin{align}
&\Obj{\alpha \w \beta}{\Measure},\nonumber\\
&\alpha < \beta\equiv(\alpha = 0^\circ \land\beta\neq 0^\circ)\lor\nonumber\\
&\quad\lor({\rm Convex~}\alpha\land{\rm Convex~}\beta\land\Exists{a \w b\w c}{\Ray}{\nonumber\\
&\quad\quad\s[6]\AC{a}{b}{\alpha}\land \AC{a}{c}{\beta} \land\Bet{a}{b}{c}})\lor\nonumber\\
&\quad\lor({\rm Convex~}\alpha\land\beta=180^\circ)\lor\nonumber\\
&\quad\lor({\rm Convex~}\alpha\land {\rm Reflex~}\beta)\lor\nonumber\\
&\quad\lor(\alpha=180^\circ\land {\rm Reflex~}\beta)\lor\nonumber\\
&\quad\lor({\rm Reflex~}\alpha\land{\rm Reflex~}\beta\land\Exists{a \w b\w c}{\Ray}{\nonumber\\
&\quad\quad\s[6]\AC{a}{b}{\alpha}\land \AC{a}{c}{\beta} \land\Bet{a}{c}{b}})\nonumber\\[0.5\baselineskip]
&\quad\CC\text{\rm  angle $\alpha$ is less then angle $\beta$}\nonumber
\end{align}} 
It is said that the angle $\alpha$ is less the the angle $\beta$ iff either of the six alternatives is valid:
\begin{enumerate}
\item[--] angle $\alpha$ is null while angle $\beta$ is non-null;
\item[--] both angles are convex and there exists system of three rays such that $\alpha$ is placed inside $\beta$;
\item[--] $\alpha$ is convex angle while $\beta$ is straight angle;
\item[--] $\alpha$ is convex angle while $\beta$ is reflex angle;
\item[--] $\alpha$ is straight angle while $\beta$ is reflex angle;
\item[--] both angles are reflex and there exists system of three rays such that $\beta$ is placed inside $\alpha$.
\end{enumerate}
\end{definition}

Based on this definition, we can prove all the necessary basic properties of a strict order relation.

\vskip0.5\baselineskip
\begin{theorem}[Non-negative angles]\Mrg
\mathnote{\Prop
\begin{align}
&\Forall{\alpha}{\Measure}{}\nonumber\\
& \lnot (\alpha < 0^\circ)\nonumber\end{align}
} 
The angle measure can not be less than the measure of the null angle.
\end{theorem}

\begin{theorem}[Angle trichotomy]\Mrg
\mathnote{\Prop
\begin{align}
&\Forall{\alpha \w \beta}{\Measure}{}\nonumber\\
& (\alpha < \beta) \lor (\alpha = \beta) \lor (\beta < \alpha)\nonumber
\end{align}} 
Whatever the two angle measures, they are either equal to each other, or one of them is less than the other.
\end{theorem}

\vskip0.5\baselineskip
\begin{theorem}[Properties of strict order]\Mrg
\mathnote{\Prop
\begin{align}
&\Forall{\alpha\w \beta\w \gamma}{\Measure}{}\nonumber
\end{align}}
The relation "less then" defines a  strict order on the angle measures.
\begin{description}
\item[\Prop\textsc{irreflexivity.}]
\mathnote{\Prop
\begin{align}
& \lnot (\alpha < \alpha)\nonumber
\end{align}}Measure of any angle is not less than itself.
\item[\Prop\textsc{transitivity.}]
\mathnote{\Prop
\begin{align}
& (\alpha < \beta) \rightarrow (\beta < \gamma) \rightarrow (\alpha < \gamma) \nonumber
\end{align}}Three angle measures $\alpha$, $\beta$ and $\gamma$ are given. 
If $\alpha$ is less than $\beta$ and $\beta$ is less than $\gamma$, then $\alpha$ is less than $\gamma$.
\end{description}
\label{thm:StrictOrder}
\end{theorem}

\begin{theorem}\Mrg
\mathnote{\Prop
\begin{align}
&\Forall{\alpha}{\Measure}{}\nonumber\\
& {\rm Convex~}\alpha \Iff 0^\circ < \alpha < 180^\circ\nonumber\\
& {\rm Reflex~}\alpha \Iff 180^\circ < \alpha\nonumber
\end{align}} 
Convex angles are those whose measure is greater than null, but less than that of the straight angle.
Reflex angles are those whose measure is greater than that of the straight angle.
\end{theorem}

\subsection{9.3. Angle addition}

With the help of our postulates, we can solve the following problem, which gives us the ability to add angles to each other.

\begin{problem}[Angle addition]\Mrg
\mathnote{\Set
\begin{align}
&\Forall{\alpha \w \beta}{\Measure}{} \nonumber\\
&\Build{\gamma}{\Angle}{\Forall{a \w b \w c}{\Ray}{\nonumber\\
&\quad\quad\quad\quad\quad\AC{a}{b}{\alpha}\rightarrow\AC{b}{c}{\beta}\rightarrow\AC{a}{c}{\gamma}}}\nonumber\\[0.5\baselineskip]
&\quad\CC\text{\rm  the resulting angle $\gamma$ is called the sum of two}\nonumber\\
&\quad\text{\rm  ~~~~angles $\alpha$ and $\beta$, i.e. we denote $\gamma = \alpha + \beta$}\nonumber
\end{align}}
Given two angle measure $\alpha$ and $\beta$.
Construct another angle measure $\gamma$ such that for any three rays $a$, $b$ and $c$, if the measure of angle $\angle (ab)$ is equal to $\alpha$ and the measure of angle $\angle (bc)$ is equal to $\beta$, then the measure of angle $\angle (ac)$ is equal to $\gamma$.
\end{problem}

Now we can prove all the necessary basic properties of the angle addition operation.

\begin{theorem}[Properties of addition]
Addition of segments has the following properties:
\mathnote{\Prop
\begin{align}
&\Forall{\alpha\w \beta\w \gamma}{\Measure}{}\nonumber
\end{align}
} 
\begin{description}
\item[\Prop\textsc{null length.}]
\mathnote{\Prop
\begin{align}
& \alpha + 0^\circ = \alpha\nonumber
\end{align}}When adding null to any angle measure, the quantity does not change.
\item[\Prop\textsc{commutativity.}]
\mathnote{\Prop
\begin{align}
& \alpha + \beta = \beta + \alpha \nonumber
\end{align}}Addition is commutative --- one can change the order of the terms in a sum, and the result is the same.
\item[\Prop\textsc{associativity.}]
\mathnote{\Prop
\begin{align}
& \alpha + (\beta + \gamma) = (\alpha + \beta) + \gamma \nonumber
\end{align}}Addition is associative --- when adding three or more angle measure, the order of operations does not matter.
\item[\Prop\textsc{subtraction.}]
\mathnote{\Prop
\begin{align}
& \alpha + \gamma = \beta + \gamma \rightarrow \alpha = \beta \nonumber
\end{align}}One can subtract the same angle measure from two equal sums, and the remaining angle measures will also be equal.
\end{description}
\end{theorem}

\section{10. Parallel lines}

Up to this point, we have studied so-called neutral geometry, i.e. planar geometry without axiom of parallel lines.
This means that all our earlier results are valid not only for Euclidean, but also for Lobachevskian geometry.

To proceed further we need to define the relation of parallel lines.
\vskip0.5\baselineskip
\begin{definition}\label{def:ParallelLines}\Mrg
\mathnote{\Def
\begin{align}
& \Obj{x\w y}{\Line}\nonumber\\
& x \parallel y \equiv \lnot x \nparallel y \nonumber\\
&\quad\CC\text{\rm  the lines $x$ and $y$ are parallel}\nonumber
\end{align}}
Two lines are called parallel iff they do not intersect.
\end{definition}

Now we can introduce our last axiom. 
\vskip0.5\baselineskip
\begin{axiom}[John Playfair]\label{axm:ParallelAxiom}\Mrg
\mathnote{\Prop
\begin{align}
&\ForallDouble{A}{\Point}{x}{\Line}{}\nonumber\\
&x \parallel z \rightarrow y \parallel z \rightarrow \In{A}{x \w y} \rightarrow x = y \nonumber
\end{align}}For every line $x$ and for every point $A$ on the plane there exists no more than one line through $A$ that is parallel to $x$.
\end{axiom}

\noindent With this axiom the relation of parallel lines becomes an {\em equivalence} relation.

\vskip0.5\baselineskip
\begin{theorem}\Mrg
The parallelism of lines is an equivalence relation and has the following properties:
\mathnote{\Prop
\begin{align}
&\Forall{x\w y\w z}{\Line}{} \nonumber
\end{align}}
\begin{description}
\item[\Prop\textsc{reflexivity.}]
\mathnote{\Prop
\begin{align}
&x \parallel x \nonumber
\end{align}}Every line is parallel to itself.
\item[\Prop\textsc{symmetry.}] 
\mathnote{\Prop
\begin{align}
& x \parallel y\rightarrow y \parallel x\nonumber
\end{align}}If line $x$ is parallel to line $y$, then line $y$ is parallel to line $x$.
\item[\Prop\textsc{transitivity.}]
\mathnote{\Prop
\begin{align}
& x \parallel y \rightarrow y \parallel z \rightarrow x \parallel z\nonumber
\end{align}}If line $x$ is parallel to line $y$ and line $y$ is parallel to line $z$, then line $x$ is parallel to line $z$.\\
 $~$\hfill{\Prop\rm \scriptsize / Euclid, Book I : Proposition 30}
\end{description}
\label{lem:ParallelEquivalence}
\end{theorem}

The following results are simple consequences of the axiom of parallel lines.

\begin{theorem}\Mrg
\mathnote{\Prop
\begin{align}
&\Forall{A \w B \w C}{\Point}\nonumber\\
& \lnot \ColThree{A}{B}{C} \rightarrow\IntAngPs{A}{B}{C} + \IntAngPs{B}{C}{A}+ \IntAngPs{C}{A}{B} = 180^\circ\nonumber
\end{align}
}In any triangle the sum of the three interior angles
 of the triangle equals two right angles.\\
$~$\hfill{\Prop\rm \scriptsize / Euclid, Book I : Proposition 32 / Hilbert, Chapter 1 : Theorem 31}
\end{theorem}

\begin{problem}\Mrg
\mathnote{\Set
\begin{align}
&\ForallDouble{A}{\Point}{x}{\Line}\nonumber\\
&\Build{y}{\Line}{A\in y\land {x}\parallel{y}}\nonumber
\end{align}
}Given line $x$ and an arbitrary point $A$, draw another line $y$ through the point $A$ that is parallel to the line $x$.
\end{problem}

\begin{problem}\Mrg
\mathnote{\Set
\begin{align}\
&\Forall{x \w y \w z}{\Line}{}\nonumber\\
& x\nparallel z \rightarrow \Or{x\nparallel y}{y\nparallel z}\nonumber
\end{align}
}Given two intersecting lines $x$ and $y$ and an arbitrary third line $z$.
Decide which of the two lines $x$ or $y$ is crossed by the line $z$.
\label{pbm:DecideIntersectingLines}
\end{problem}

\begin{problem}\Mrg
\mathnote{\Set
\begin{align}\
&\Forall{x \w y \w z}{\Line}{}\nonumber\\
& x\neq z \rightarrow \Or{x\neq y}{y\neq z}\nonumber
\end{align}
}Given two distinct lines $x$ and $y$ and an arbitrary third line $z$.
Decide which of the two lines $x$ or $y$ is distinct from the line $z$.
\end{problem}

\noindent The last two problems imply important properties of the line distinction and intersection relations, analogous to the Hayting postulate for the distinction of points (Postulate \ref{pst:DecideDistinctPoints}). 

\section{11. Comparison with Hilbert's system of axioms}

In conclusion, let us compare Hilbert's axioms line by line with the postulates and axioms of the proposed constructive-deductive method for planar Euclidean geometry.

\begin{table*}[p]
\begin{tabular*}{\textwidth}{p{0.475\textwidth}p{0.475\textwidth}}
\toprule\\[-4mm]
\text{\bf \Prop Axiom {$\text{\rm I}_1$}.} \em For every two distinct points $A, B$ there exists a line $x$ that contains each of the points  $A, B$.
&\text{\bf \Set Postulate \refcolor{MyGreen}\ref{pst:DrawExtensionLine}.} \em 
Given two distinct points $A$ and $B$. Draw a straight line $x$ through both these points.\\[2mm]
\begin{tcolorbox}
[colback=MyBlue!5!white, frame hidden, enhanced, attach boxed title to top right={yshift=-5.2mm,xshift=-0.2mm}, 
boxed title style={size=small, frame hidden, colback=MyBlue},
title=\footnotesize {$\text{\rm I}_1$}]
\footnotesize $\Forall{A \w B}{\Point}{A\neq B\rightarrow\Exists{x}{\Line}{A \w B \in x}}$ 
\end{tcolorbox}& 
\begin{tcolorbox}
[colback=MyGreen!5!white, frame hidden, enhanced, attach boxed title to top right={yshift=-5.2mm,xshift=-0.2mm}, 
boxed title style={size=small, frame hidden, colback=MyGreen},
title=\footnotesize\Postulate{\refcolor{white}\ref{pst:DrawExtensionLine}}]
\footnotesize $\Forall{A \w B}{\Point}{A\neq B\rightarrow\Build{x}{\Line}{A \w B \in x}}
$ 
\end{tcolorbox}
\\[-5mm]
\midrule
\text{\bf \Prop Axiom {$\text{\rm I}_2$}.} \em For every two distinct points $A, B$ there exists not more than one line that contains each of the points $A, B$.
&\text{\bf \Set Postulate \refcolor{MyGreen}\ref{pst:DecideUniqueExtensionLine}} (Jan von Plato){\bf \Set .} \em 
Given two distinct lines, on one of which two distinct points are marked. 
Decide which of these two points does not belong the other line.
\\[2mm]
\begin{tcolorbox}
[colback=MyBlue!5!white, frame hidden, enhanced, attach boxed title to top right={yshift=-5.2mm,xshift=-0.2mm}, 
boxed title style={size=small, frame hidden, colback=MyBlue},
title=\footnotesize {$\text{\rm I}_2$}]
\footnotesize $\ForallDouble{A \w B}{\Point}{x \w y}{\Line}{}\\[1mm]
A\neq B\rightarrow \In{A \w B}{x \w y} \rightarrow x = y$ 
\end{tcolorbox}& 
\begin{tcolorbox}
[colback=MyGreen!5!white, frame hidden, enhanced, attach boxed title to top right={yshift=-5.2mm,xshift=-0.2mm}, 
boxed title style={size=small, frame hidden, colback=MyGreen},
title=\footnotesize\Postulate{\refcolor{white}\ref{pst:DecideUniqueExtensionLine}}]
\footnotesize $\ForallDouble{A \w B}{\Point}{x \w y}{\Line}{}\\[1mm]
A \neq B \rightarrow \In{A\w B}{x} \rightarrow  x \neq y\rightarrow\Or{A \notin y}{B \notin y}$ 
\end{tcolorbox}
\\[-5mm]
\midrule
\text{\bf \Prop Axiom {$\text{\rm I}_{3}~a)$}.} 
\em There exists at least two distinct points on a line. 
&
\text{\bf \Set Postulate \refcolor{MyGreen}\ref{pst:DrawPointOnLine}.} \em 
Given an arbitrary straight line. Draw a point on the line.
\\[2mm]
\begin{tcolorbox}
[colback=MyBlue!5!white, frame hidden, enhanced, attach boxed title to top right={yshift=-5.2mm,xshift=-0.2mm}, 
boxed title style={size=small, frame hidden, colback=MyBlue},
title=\footnotesize {$\text{\rm I}_{3a}$}]
\footnotesize $\Forall{x}{\Line}{\Exists{A\w B}{\Point}{A\neq B \land (A \w B \in x)}}$ 
\end{tcolorbox}& 
\begin{tcolorbox}
[colback=MyGreen!5!white, frame hidden, enhanced, attach boxed title to top right={yshift=-5.2mm,xshift=-0.2mm}, 
boxed title style={size=small, frame hidden, colback=MyGreen},
title=\footnotesize\Postulate{\refcolor{white}\ref{pst:DrawPointOnLine}}]
\footnotesize $\Forall{x}{\Line}{\Build{A}{\Point}{A \in x}}$ 
\end{tcolorbox}
\\[-5mm]
&\text{\bf \Set Postulate \refcolor{MyGreen}\ref{pst:DrawDistinctPointOnLine}.} \em 
Given an arbitrary straight line and a point on it. 
Draw another point on this line distinct from the given one.
\\[2mm]
& 
\begin{tcolorbox}
[colback=MyGreen!5!white, frame hidden, enhanced, attach boxed title to top right={yshift=-5.2mm,xshift=-0.2mm}, 
boxed title style={size=small, frame hidden, colback=MyGreen},
title=\footnotesize\Postulate{\refcolor{white}\ref{pst:DrawDistinctPointOnLine}}]
\footnotesize $\ForallDouble{A}{\Point}{x}{\Line}{}\\[1mm]
A \in x \rightarrow \Build{B}{\Point}{A \neq B \land B \in x}$ 
\end{tcolorbox}
\\[-5mm]
\midrule
\text{\bf \Prop Axiom {$\text{\rm I}_{3} ~b)$}.} 
\em There exists at least three points that do not lie on a line.
&\text{\bf \Set Postulate \refcolor{MyGreen}\ref{pst:DrawPoint}.} \em 
Draw a point on the plane.\\[1.5\baselineskip]
\begin{tcolorbox}
[colback=MyBlue!5!white, frame hidden, enhanced, attach boxed title to top right={yshift=-5.2mm,xshift=-0.2mm}, 
boxed title style={size=small, frame hidden, colback=MyBlue},
title=\footnotesize {$\text{\rm I}_{3b}$}]
\footnotesize $\Exists{A\w B \w C}{\Point}{\lnot\ColThree{A}{B}{C}}$ 
\end{tcolorbox}& 
\begin{tcolorbox}
[colback=MyGreen!5!white, frame hidden, enhanced, attach boxed title to top right={yshift=-5.2mm,xshift=-0.2mm}, 
boxed title style={size=small, frame hidden, colback=MyGreen},
title=\footnotesize\Postulate{\refcolor{white}\ref{pst:DrawPoint}}]
\footnotesize $\{\w A : \Point \w\}$ 
\end{tcolorbox}
\\[-5mm]
&\text{\bf \Set Postulate \refcolor{MyGreen}\ref{pst:DrawDistinctPoint}.} \em 
Given an arbitrary point on the plane.
Draw another point distinct from the given one.
\\[2mm]
& 
\begin{tcolorbox}
[colback=MyGreen!5!white, frame hidden, enhanced, attach boxed title to top right={yshift=-5.2mm,xshift=-0.2mm}, 
boxed title style={size=small, frame hidden, colback=MyGreen},
title=\footnotesize\Postulate{\refcolor{white}\ref{pst:DrawDistinctPoint}}]
\footnotesize $\Forall{A}{\Point}{\Build{B}{\Point}{A \neq B}}$ 
\end{tcolorbox}
\\[-5mm]
&\text{\bf \Set Postulate \refcolor{MyGreen}\ref{pst:DrawPointApartLine}.} \em 
Given an arbitrary line.
Draw a point that does not lie on this line.
\\[2mm]
& 
\begin{tcolorbox}
[colback=MyGreen!5!white, frame hidden, enhanced, attach boxed title to top right={yshift=-5.2mm,xshift=-0.2mm}, 
boxed title style={size=small, frame hidden, colback=MyGreen},
title=\footnotesize\Postulate{\refcolor{white}\ref{pst:DrawPointApartLine}}]
\footnotesize $\Forall{x}{\Line}{\Build{A}{\Point}{A \notin x}}$ 
\end{tcolorbox}
\\[-5mm]
\midrule
There is no analog of Heyting's postulate in Hilbert's system of axioms.
&\text{\bf \Set Postulate \refcolor{MyGreen}\ref{pst:DecideDistinctPoints}} (Arend Heyting){\bf \Set .} \em 
Given two distinct points $A$ and $B$ and an arbitrary point  $C$. 
Decide the distinction of the point $C$ either from the point $A$ or from the point $B$.
\\[2mm]
& 
\begin{tcolorbox}
[colback=MyGreen!5!white, frame hidden, enhanced, attach boxed title to top right={yshift=-5.2mm,xshift=-0.2mm}, 
boxed title style={size=small, frame hidden, colback=MyGreen},
title=\footnotesize\Postulate{\refcolor{white}\ref{pst:DecideDistinctPoints}}]
\footnotesize $\Forall{A \w B \w C}{\Point}{A\neq B \rightarrow \Or{A \neq C}{C \neq B}}$ 
\end{tcolorbox}
\\[-5mm]
\midrule
There is no analog of line-line intersection postulate in Hilbert's system of axioms.
&\text{\bf \Set Postulate \refcolor{MyGreen}\ref{pst:DrawIntersectionPoint}} (Line-line intersection){\bf \Set .} \em 
Given two intersecting lines. 
Draw the point of their intersection.
\\[2mm]
& 
\begin{tcolorbox}
[colback=MyGreen!5!white, frame hidden, enhanced, attach boxed title to top right={yshift=-5.2mm,xshift=-0.2mm}, 
boxed title style={size=small, frame hidden, colback=MyGreen},
title=\footnotesize\Postulate{\refcolor{white}\ref{pst:DrawIntersectionPoint}}]
\footnotesize $\Forall{x\w y}{\Line}{x\nparallel y\rightarrow\Build{A}{\Point}{A \in x \w y}}$ 
\end{tcolorbox}
\\[-5mm]
\bottomrule
\end{tabular*}
  \caption{Group I: Axioms of incidence.}
  \label{tab:GI}
\end{table*}

\begin{table*}[p]
\begin{tabular*}{\textwidth}{p{0.475\textwidth}p{0.475\textwidth}}
\toprule\\[-4mm]
\text{\bf \Prop Axiom {$\text{\rm II}_{1,2}$}.} \em 
If a point $B$ lies between a point $A$ and a point $C$ then:
\begin{enumerate}
\item[a)] the points $A, B, C$ are three distinct points;
\item[b)] the points $A, B, C$ are collinear;
\item[c)] the point $B$ also lies between the points $C$ and $A$;
\item[d)] the point $A$ cannot lie between the points $B$ and $C$.
\end{enumerate}
&\text{\bf \Prop Axiom \refcolor{MyGreen}\ref{axm:BetPs_unique}.} \em 
Given three points on the plane $A$, $B$ and $C$.
If the point $B$ lies between the points $A$ and $C$, then:
\begin{enumerate}
\item[a)] the points $A$ and $C$ are distinct;
\item[b)] all three points, $A$, $B$ and $C$ are collinear;
\item[c)] the point $B$ also lies between the points $C$ and $A$;
\item[d)] the point $A$ cannot lie between the points $B$ and $C$.
\end{enumerate}\\[-3mm]
\begin{tcolorbox}
[colback=MyBlue!5!white, frame hidden, enhanced, attach boxed title to top right={yshift=-5.2mm,xshift=-0.2mm}, 
boxed title style={size=small, frame hidden, colback=MyBlue},
title=\footnotesize {$\text{\rm II}_{1,2}$}]
\footnotesize $\Forall{A \w B \w C}{\Point}{}\\[1mm]
\Bet{A}{B}{C}\rightarrow  A \neq C \land\ColThree{A}{B}{C}\land\Bet{C}{B}{A}\land \lnot \Bet{B}{A}{C}$
\end{tcolorbox}& 
\begin{tcolorbox}
[colback=MyBlue!5!white, frame hidden, enhanced, attach boxed title to top right={yshift=-5.2mm,xshift=-0.2mm}, 
boxed title style={size=small, frame hidden, colback=MyBlue},
title=\footnotesize\Axiom{\refcolor{white}\ref{axm:BetPs_unique}}]
\footnotesize $\Forall{A \w B \w C}{\Point}{}\\[1mm]
\Bet{A}{B}{C}\rightarrow  A \neq C \land\ColThree{A}{B}{C}\land\Bet{C}{B}{A}\land \lnot \Bet{B}{A}{C}$
\end{tcolorbox}
\\[-5mm]
\midrule
\text{\bf \Prop Axiom {$\text{\rm II}_3$}.} \em 
For two distinct points $A$ and $B$, there always exists at least one point $C$ such that point $B$ lies between $A$ and $C$.
&  This Hilbert's axiom follows immediately from the line-circle intersection postulate (Postulate \ref{pst:DrawIntersectionOfLineAndCircle}).
\\[5mm]
\begin{tcolorbox}
[colback=MyBlue!5!white, frame hidden, enhanced, attach boxed title to top right={yshift=-5.2mm,xshift=-0.2mm}, 
boxed title style={size=small, frame hidden, colback=MyBlue},
title=\footnotesize {$\text{\rm II}_3$}]
\footnotesize $\Forall{A \w B}{\Point}{(A\neq B\rightarrow \Exists{C}{\Point}{\Bet{A}{B}{C}})}$ 
\end{tcolorbox}& 
\\[-5mm]
\midrule
\text{\bf \Prop Axiom {$\text{\rm II}_4$}}(Morits Pasch){\bf\bf \Prop .} \em 
Let $A, B, C$ be three points that do not lie on a line and let $x$ be a line which does not meet any of the points $A, B, C$. 
If the line $x$ passes through a point of the segment $\Seg{A}{B}$, it also passes through a point of the segment $\Seg{A}{C}$, or through a point of the segment $\Seg{C}{B}$.
&\text{\bf \Set Postulate \refcolor{MyGreen}\ref{pst:DecideOppositeSide}} (Morits Pasch){\bf \Set .} \em 
Given three points $A, B, C$ and a line $x$ that does not pass through any of these points and intersects the segment $\Seg{A}{B}$.
Decide the intersection of the line $x$ either with the segment $\Seg{A}{C}$ or with the segment $\Seg{C}{B}$.
\\[2mm]
\begin{tcolorbox}
[colback=MyBlue!5!white, frame hidden, enhanced, attach boxed title to top right={yshift=-5.2mm,xshift=-0.2mm}, 
boxed title style={size=small, frame hidden, colback=MyBlue},
title=\footnotesize {$\text{\rm II}_4$}]
\footnotesize $\ForallDouble{A \w B \w C}{\Point}{x}{\Line}{}\\[1mm]
\lnot\ColThree{A}{B}{C}\rightarrow \OppositeSide{A}{x}{B} \rightarrow  C \notin x \rightarrow\OppositeSide{A}{x}{C} \lor \OppositeSide{C}{x}{B}$ 
\end{tcolorbox}& 
\begin{tcolorbox}
[colback=MyGreen!5!white, frame hidden, enhanced, attach boxed title to top right={yshift=-5.2mm,xshift=-0.2mm}, 
boxed title style={size=small, frame hidden, colback=MyGreen},
title=\footnotesize\Postulate{\refcolor{white}\ref{pst:DecideOppositeSide}}]
\footnotesize $\ForallDouble{A \w B \w C}{\Point}{x}{\Line}{}\\[1mm]
\OppositeSide{A}{x}{B} \rightarrow C \notin x \rightarrow \Or{\OppositeSide{A}{x}{C} }{\OppositeSide{C}{x}{B}}$ 
\end{tcolorbox}
\\[-5mm]
\bottomrule
\end{tabular*}
  \caption{Group II: Axioms of order.}
  \label{tab:GII}
\end{table*}

\begin{table*}[p]
\begin{tabular*}{\textwidth}{p{0.475\textwidth}p{0.475\textwidth}}
\toprule\\[-4mm]
\text{\bf \Prop Axiom {$\text{\rm III}_1$}.} \em 
If $A, B$ are two distinct points, and if $C$ is a point on some line $x$ then it is always possible to find a point $D$ on a given side of the line line $x$ through $C$ such that the segment $\Seg{A}{B}$ is congruent or equal to the segment $\Seg{C}{D}$.
&\text{\bf \Set Postulate \refcolor{MyGreen}\ref{pst:DrawIntersectionOfLineAndCircle}} (Line-Circle intersection){\bf \Set .} \em 
Given two distinct points $A$ and $O$ and an arbitrary point $B$ on the plane. 
Draw the point $C$ of the intersection of the line $\Lin{A}{O}$ and the circle $\Circ{O}{B}$ such that the point $O$ belongs to the segment $\Seg{A}{C}$.
\\[2mm]
\begin{tcolorbox}
[colback=MyBlue!5!white, frame hidden, enhanced, attach boxed title to top right={yshift=-5.2mm,xshift=-0.2mm}, 
boxed title style={size=small, frame hidden, colback=MyBlue},
title=\footnotesize {$\text{\rm III}_1$}]
\footnotesize $\Forall{A \w B \w C\w P}{\Point}{}\\[1mm]
A\neq B\rightarrow C\neq P\rightarrow\Exists{D}{\Point}{\SameRay{C}{P}{D}\land \Seg{A}{B} \cong \Seg{C}{D}}$ 
\end{tcolorbox}& 
\begin{tcolorbox}
[colback=MyGreen!5!white, frame hidden, enhanced, attach boxed title to top right={yshift=-5.2mm,xshift=-0.2mm}, 
boxed title style={size=small, frame hidden, colback=MyGreen},
title=\footnotesize\Postulate{\refcolor{white}\ref{pst:DrawIntersectionOfLineAndCircle}}]
\footnotesize $\Forall{O \w A \w B}{\Point}{}\\[1mm]
A\neq O\rightarrow \Build{C}{\Point}{\BetXO{A}{O}{C} \land \Lcong{O}{B}{O}{C}}$ 
\end{tcolorbox}
\\[-5mm]
\midrule
\text{\bf \Prop Axiom {$\text{\rm III}_2$}.} \em 
If a segment $\Seg{C}{D}$ and a segment $\Seg{E}{F}$, are congruent to the same segment $\Seg{A}{B}$, then the segment $\Seg{C}{D}$ is congruent to the segment $\Seg{E}{F}$.
&Instead of congruence relation we introduce congruence classes as lengths of segments (pairs of points). Thus we introduce two functions:
\\[2mm]
\begin{tcolorbox}
[colback=MyBlue!5!white, frame hidden, enhanced, attach boxed title to top right={yshift=-5.2mm,xshift=-0.2mm}, 
boxed title style={size=small, frame hidden, colback=MyBlue},
title=\footnotesize {$\text{\rm III}_2$}]
\footnotesize $\Forall{A \w B \w C\w D\w E\w F}{\Point}{}\\[1mm]
\Seg{A}{B}\cong \Seg{C}{D} \rightarrow \Seg{A}{B}\cong \Seg{E}{F}\rightarrow \Seg{C}{D}\cong \Seg{E}{F}$ 
\end{tcolorbox}& 
\begin{tcolorbox}
[colback=MyGreen!5!white, frame hidden, enhanced, attach boxed title to top right={yshift=-5.2mm,xshift=-0.2mm}, 
boxed title style={size=small, frame hidden, colback=MyGreen}]
\footnotesize ${\mathbb L^{\rm +}} (s : \text{Segment}) : \text{Length}\\[1mm]
{\mathbb L^{\rm -}} (d : \text{Length}) : \Build{s}{\text{Segment}}{{\mathbb L^{\rm +}}(s) = d}$ 
\end{tcolorbox}
\\[-5mm]
\midrule
\text{\bf \Prop Axiom {$\text{\rm III}_3$}.} \em 
On the line $x$ let $\Seg{A}{B}$ and $\Seg{B}{C}$ be two segments which except for $B$ have no point in common. Furthermore, on the same or of another line $y$
let $\Seg{\Prime{A}}{\Prime{B}}$ and $\Seg{\Prime{B}}{\Prime{C}}$ be two segments which except for $\Prime{B}$ also have no point in common. 
In that case, if $\Seg{A}{B} \cong \Seg{\Prime{A}}{\Prime{B}}$ and $\Seg{B}{C} \cong \Seg{\Prime{B}}{\Prime{C}}$, then $\Seg{A}{C} \cong \Seg{\Prime{A}}{\Prime{C}}$.
&\text{\bf \Prop Axiom \refcolor{MyGreen}\ref{axm:CircleSegmentSubtraction}} (Concentric circles and rays){\bf\bf \Prop .} \em 
Consider two concentric circles and two rays originating from the center of the circles $O$.
Let points $A$ and $\Prime{A}$ lie on one of the circles, and points $B$ and $\Prime {B}$ on another.
If in this case the point $A$ lies on one of the rays between the points $O$ and $B$, and points $\Prime{A}$ and $\Prime{B}$ on another ray, then the point $\Prime{A}$ lies between the points $O$ and $\Prime{B}$ and the lengths of the segments $\Seg{A}{B}$ and $\Seg{\Prime{A}}{\Prime{B}}$ are equal.
\\[2mm]
\begin{tcolorbox}
[colback=MyBlue!5!white, frame hidden, enhanced, attach boxed title to top right={yshift=-5.2mm,xshift=-0.2mm}, 
boxed title style={size=small, frame hidden, colback=MyBlue},
title=\footnotesize {$\text{\rm III}_3$}]
\footnotesize $\Forall{A \w B \w C\w \Prime{A} \w \Prime{B}\w \Prime{C}}{\Point}{}\\[1mm]
\Bet{A}{B}{C}\rightarrow \Bet{\Prime{A}}{\Prime{B}}{ \Prime{C}}\rightarrow\\[1mm]
~~~ \rightarrow\Seg{A}{B} \cong \Seg{\Prime{A}}{\Prime{B}} \rightarrow\Seg{B}{C} \cong \Seg{\Prime{B}}{ \Prime{C}}\rightarrow\Seg{A}{C} \cong \Seg{\Prime{A}}{ \Prime{C}}$ 
\end{tcolorbox}& 
\begin{tcolorbox}
[colback=MyBlue!5!white, frame hidden, enhanced, attach boxed title to top right={yshift=-5.2mm,xshift=-0.2mm}, 
boxed title style={size=small, frame hidden, colback=MyBlue},
title=\footnotesize\Axiom{\refcolor{white}\ref{axm:CircleSegmentSubtraction}}]
\footnotesize $\Forall{O \w A \w \Prime{A} \w B \w \Prime{B}}{\Point}{}\\[1mm]
\Bet{O}{A}{B} \rightarrow \SameRay{O}{\Prime{A}}{\Prime{B}} \rightarrow\Lcong{O}{A}{O}{\Prime{A}} \rightarrow\\[1mm]
~~~ \rightarrow \Lcong{O}{B}{O}{\Prime{B}}  \rightarrow \Bet{O}{\Prime{A}}{\Prime{B}}\land\Lcong{A}{B}{\Prime{A}}{\Prime{B}}$ 
\end{tcolorbox}
\\[-5mm]
\midrule
This axiom is implicit in Hilbert's system of axioms since he defines segment as a set of two distinct points.
&\text{\bf \Prop Axiom \refcolor{MyBlue}\ref{axm:SegPs_sym}}{\bf \Prop .} \em 
The length of the segment does not depend on the order of its endpoints.
\\[2mm]
\begin{tcolorbox}
[colback=MyBlue!5!white, frame hidden, enhanced, attach boxed title to top right={yshift=-5.2mm,xshift=-0.2mm}, 
boxed title style={size=small, frame hidden, colback=MyBlue},
title=\footnotesize {$\text{\rm III}_0$}]
\footnotesize $\Forall{A \w B}{\Point}{\Seg{A}{B} \cong \Seg{B}{A}}$ 
\end{tcolorbox}& 
\begin{tcolorbox}
[colback=MyBlue!5!white, frame hidden, enhanced, attach boxed title to top right={yshift=-5.2mm,xshift=-0.2mm}, 
boxed title style={size=small, frame hidden, colback=MyBlue},
title=\footnotesize\Axiom{\refcolor{white}\ref{axm:SegPs_sym}}]
\footnotesize $\Forall{A \w B}{\Point}{\Lcong{A}{B}{B}{A}}$ 
\end{tcolorbox}
\\[-5mm]
\midrule
Null segments are not defined in Hilbert's system of axioms.
&\text{\bf \Prop Axiom \refcolor{MyBlue}\ref{axm:EqSs_EqPs}} (Null segments){\bf \Prop .} \em 
The lengths of all null segments are equal to each other.
Any segment whose length is equal to the length of some null segment is itself a null segment. 
\\[2mm]
& \begin{tcolorbox}
[colback=MyBlue!5!white, frame hidden, enhanced, attach boxed title to top right={yshift=-5.2mm,xshift=-0.2mm}, 
boxed title style={size=small, frame hidden, colback=MyBlue},
title=\footnotesize\Axiom{\refcolor{white}\ref{axm:EqSs_EqPs}}]
\footnotesize $\Forall{A \w B \w C}{\Point}{\Lcong{A}{B}{C}{C} \Iff A = B}$ 
\end{tcolorbox}
\\[-5mm]
\bottomrule
\end{tabular*}
  \caption{Group IIIa: Axioms of segment congruence.}
  \label{tab:GIIIa}
\end{table*}

\begin{table*}[p]
\begin{tabular*}{\textwidth}{p{0.475\textwidth}p{0.475\textwidth}}
\toprule\\[-4mm]
\text{\bf \Prop Axiom {$\text{\rm III}_4 ~a)$}.} \em 
If $\angle ABC$ is an angle and if $\RayObj{\Prime{B}}{\Prime{A}}$ is a ray, then there exists a ray $\RayObj{\Prime{B}}{\Prime{C}}$ on given side of line $\Lin{\Prime{B}}{\Prime{A}}$ such that $\angle ABC\cong\angle \Prime{A}\Prime{B}\Prime{C}$.
&\text{\bf \Set Postulate \refcolor{MyGreen}\ref{pst:DrawIntersectionOfTwoCircles}} (Circle-Circle intersection){\bf \Set .} \em 
Draw a point of intersection of two circles with each other.
\\[-3mm]
\begin{tcolorbox}
[colback=MyBlue!5!white, frame hidden, enhanced, attach boxed title to top right={yshift=-5.2mm,xshift=-0.2mm}, 
boxed title style={size=small, frame hidden, colback=MyBlue},
title=\footnotesize {$\text{\rm III}_{4a}$}]
\footnotesize $\Forall{A \w B \w C\w \Prime{A}\w \Prime{B}\w P}{\Point}{\lnot\ColThree{A}{B}{C}\rightarrow \lnot\ColThree{\Prime{A}}{\Prime{B}}{P}\rightarrow}\\[1mm]
~~~\rightarrow\Exists{\Prime{C}}{\Point}{\SameHalfplane{\Prime{A}}{\Prime{B}}{P}{\Prime{C}}\land \angle ABC \cong \angle \Prime{A} \Prime{B} \Prime{C}}$ 
\end{tcolorbox}& 
\begin{tcolorbox}
[colback=MyGreen!5!white, frame hidden, enhanced, attach boxed title to top right={yshift=-5.2mm,xshift=-0.2mm}, 
boxed title style={size=small, frame hidden, colback=MyGreen},
title=\footnotesize\Postulate{\refcolor{white}\ref{pst:DrawIntersectionOfTwoCircles}}]
\footnotesize $\Forall{O \w A \w B \w \Prime{O} \w \Prime{A} \w \Prime{B} \w P}{\Point}{\lnot\ColThree{O}{\Prime{O}}{P} \rightarrow}\\[1mm]
~~~\rightarrow\BetFour{A}{\Prime{A}}{B}{\Prime{B}}\rightarrow\Middle{A}{O}{B}\rightarrow\Middle{\Prime{A}\s[-2]}{\Prime{O}\s[-2]}{\Prime{B}}\rightarrow\\[1mm]
~~~\rightarrow\{\s \Obj{Q}{\Point} \vl \SameHalfplane{O}{\Prime{O}}{P}{Q} \land\\[1mm]
~~~~~~~~~~~~~~~~~~~~\land \Lcong{O}{Q}{O}{A} \land \Lcong{\Prime{O}}{Q}{\Prime{O}}{\Prime{A}} \s\}$ 
\end{tcolorbox}
\\[-5mm]
\midrule
\text{\bf \Prop Axiom {$\text{\rm III}_4 ~b)$}.} \em 
There exists not more then one angle from the same ray in the same half-plane.
&\text{\bf \Prop Axiom \refcolor{MyBlue}\ref{axm:EqAs_EqRs}} (Unique angles){\bf \Prop .} \em 
Given two angles initial sides of which coincide. Then the terminal sides of this angles coincide iff the angles are equal.
\\[-4mm]
\begin{tcolorbox}
[colback=MyBlue!5!white, frame hidden, enhanced, attach boxed title to top right={yshift=-5.2mm,xshift=-0.2mm}, 
boxed title style={size=small, frame hidden, colback=MyBlue},
title=\footnotesize {$\text{\rm III}_{4b}$}]
\footnotesize $\Forall{A \w B \w C\w \Prime{A}\w \Prime{B} \w P \w \Prime{C} \w \Prime{\Prime{C}}}{\Point}{\lnot\ColThree{\Prime{A}}{\Prime{B}}{P}\rightarrow}\\[1mm]
\rightarrow \angle ABC \cong \angle \Prime{A} \Prime{B} \Prime{C}\rightarrow \angle ABC \cong \angle \Prime{A} \Prime{B} \Prime{\Prime{C}}\rightarrow\\[1mm]
\rightarrow\SameHalfplane{\Prime{A}}{\Prime{B}}{P}{\Prime{C}} \rightarrow\SameHalfplane{\Prime{A}}{\Prime{B}}{P}{\Prime{\Prime{C}}} \rightarrow \SameRay{\Prime{B}}{\Prime{C}}{\Prime{\Prime{C}}}$ 
\end{tcolorbox}& 
\begin{tcolorbox}
[colback=MyBlue!5!white, frame hidden, enhanced, attach boxed title to top right={yshift=-5.2mm,xshift=-0.2mm}, 
boxed title style={size=small, frame hidden, colback=MyBlue},
title=\footnotesize\Postulate{\refcolor{white}\ref{axm:EqAs_EqRs}}]
\footnotesize $\Forall{a \w b \w c}{\Ray}{}\\[1mm]
\AngCongRs{a}{b}{a}{c} \Iff b \approx c$ 
\end{tcolorbox}
\\[-5mm]
\midrule
\text{\bf \Prop Axiom {$\text{\rm III}_{4}~c)$}.} \em 
Every angle is congruent to itself, i.e. $\angle(a, b) \cong \angle (b, a)$ is always true.
&Instead of congruence relation we introduce congruence classes as measures of angles (ordered pairs of rays). 
Thus we introduce two functions:
\\[2mm]
\begin{tcolorbox}
[colback=MyBlue!5!white, frame hidden, enhanced, attach boxed title to top right={yshift=-5.2mm,xshift=-0.2mm}, 
boxed title style={size=small, frame hidden, colback=MyBlue},
title=\footnotesize {$\text{\rm III}_{4c}$}]
\footnotesize $\Forall{a \w b}{\Ray}{}\\[1mm]
\angle(\Seg{a}{b}) \cong \angle(\Seg{a}{b})$ 
\end{tcolorbox}& 
\begin{tcolorbox}
[colback=MyGreen!5!white, frame hidden, enhanced, attach boxed title to top right={yshift=-5.2mm,xshift=-0.2mm}, 
boxed title style={size=small, frame hidden, colback=MyGreen}]
\footnotesize ${\mathbb M^{\rm +}} (A : \Angle) : \Measure\\[1mm]
{\mathbb M^{\rm -}} (\alpha : \Measure) : \Build{A}{\Angle}{{\mathbb M^{\rm +}}(A) = \alpha}$ 
\end{tcolorbox}
\\[-5mm]
\midrule
\text{\bf \Prop Axiom {$\text{\rm III}_5$}.} \em 
If for two triangles $\Triangle{A}{B}{C}$ and $\Triangle{\Prime{A}}{\Prime{B}}{\Prime{C}}$ the congruences 
$\Seg{A}{B} \cong \Seg{\Prime{A}}{\Prime{B}}$ , $\Seg{A}{C} \cong \Seg{\Prime{A}}{\Prime{C}}$ and $\angle BAC \cong \angle\Prime{B}\Prime{A}\Prime{C}$
hold, then the congruence $\angle ABC \cong \angle\Prime{A}\Prime{B}\Prime{C}$ is also satisfied.
&\text{\bf \Prop Axiom \refcolor{MyGreen}\ref{axm:CongruencePrinciple}} (Triangle congruence){\bf\bf \Prop .} \em 
Let two sides of one triangle be equal respectively to two sides of another triangle.
Then the third sides of these triangles will be equal, iff the angles opposite to them are equal.
\\[2mm]
\begin{tcolorbox}
[colback=MyBlue!5!white, frame hidden, enhanced, attach boxed title to top right={yshift=-5.2mm,xshift=-0.2mm}, 
boxed title style={size=small, frame hidden, colback=MyBlue},
title=\footnotesize {$\text{\rm III}_5$}]
\footnotesize $\Forall{A \w B \w C \w \Prime{A} \w \Prime{B} \w \Prime{C}}{\Point}{\lnot \ColThree{A}{B}{C} \rightarrow \lnot \ColThree{\Prime{A}}{\Prime{B}}{\Prime{C}} \rightarrow}\\[1mm]
~~~\rightarrow\Seg{A}{B} \cong \Seg{\Prime{A}}{\Prime{B}}\rightarrow\Seg{A}{C} \cong \Seg{\Prime{A}}{\Prime{C}}\rightarrow\angle BAC \cong \angle\Prime{B}\Prime{A}\Prime{C}\rightarrow\\[1mm]
~~~\rightarrow\angle ABC \cong \angle\Prime{A}\Prime{B}\Prime{C}$ 
\end{tcolorbox}& 
\begin{tcolorbox}
[colback=MyBlue!5!white, frame hidden, enhanced, attach boxed title to top right={yshift=-5.2mm,xshift=-0.2mm}, 
boxed title style={size=small, frame hidden, colback=MyBlue},
title=\footnotesize\Axiom{\refcolor{white}\ref{axm:CongruencePrinciple}}]
\footnotesize $\Forall{A \w B \w C \w D \w E \w F}{\Point}{\lnot \ColThree{A}{B}{C} \rightarrow \lnot \ColThree{D}{E}{F} \rightarrow}\\[1mm]
~~~\rightarrow\Lcong{A}{B}{D}{E}\rightarrow\Lcong{B}{C}{E}{F}\rightarrow\\[1mm]
~~~\rightarrow\IntAngCongPs{A}{B}{C}{D}{E}{F}\Iff \Lcong{C}{A}{F}{D}$ 
\end{tcolorbox}
\\[-5mm]
\midrule
This axiom is implicit in Hilbert's system of axioms since he defines angles as a set of two distinct rays.
&\text{\bf \Prop Axiom \refcolor{MyBlue}\ref{axm:EqAs_EqExpAs}} (Explementary angles){\bf \Prop .} \em 
If two angles are equal, then their explementary angles are also equal. 
\\[2mm]
\begin{tcolorbox}
[colback=MyBlue!5!white, frame hidden, enhanced, attach boxed title to top right={yshift=-5.2mm,xshift=-0.2mm}, 
boxed title style={size=small, frame hidden, colback=MyBlue},
title=\footnotesize {$\text{\rm III}_0$}]
\footnotesize $\Forall{a \w b}{\Ray}{\angle(\Seg{a}{b}) \cong \angle(\Seg{b}{a})}$ 
\end{tcolorbox}& 
\begin{tcolorbox}
[colback=MyBlue!5!white, frame hidden, enhanced, attach boxed title to top right={yshift=-5.2mm,xshift=-0.2mm}, 
boxed title style={size=small, frame hidden, colback=MyBlue},
title=\footnotesize\Axiom{\refcolor{white}\ref{axm:EqAs_EqExpAs}}]
\footnotesize $\Forall{a \w b \w c \w d}{\Ray}{\AngCongRs{a}{b}{c}{d} \rightarrow \AngCongRs{b}{a}{d}{c}}$ 
\end{tcolorbox}
\\[-5mm]
\midrule
Null angles are not defined in Hilbert's system of axioms.
&\text{\bf \Prop Axiom \refcolor{MyBlue}\ref{axm:EqRs_EqRs_EqAs}} (Null angles){\bf \Prop .} \em 
All null angles on the plane are equal to each other.
\\[2mm]
& 
\begin{tcolorbox}
[colback=MyBlue!5!white, frame hidden, enhanced, attach boxed title to top right={yshift=-5.2mm,xshift=-0.2mm}, 
boxed title style={size=small, frame hidden, colback=MyBlue},
title=\footnotesize\Axiom{\refcolor{white}\ref{axm:EqRs_EqRs_EqAs}}]
\footnotesize $\Forall{a \w b}{\Ray}{\AngCongRs{a}{a}{b}{b}}$ 
\end{tcolorbox}
\\[-5mm]
\midrule
Straight angles are not defined in Hilbert's system of axioms.
&\text{\bf \Prop Axiom \refcolor{MyBlue}\ref{axm:OpRs_OpRs_EqAs}} (Straight angles){\bf \Prop .} \em 
All straight angles on the plane are equal to each other. 
\\[2mm]
& 
\begin{tcolorbox}
[colback=MyBlue!5!white, frame hidden, enhanced, attach boxed title to top right={yshift=-5.2mm,xshift=-0.2mm}, 
boxed title style={size=small, frame hidden, colback=MyBlue},
title=\footnotesize\Axiom{\refcolor{white}\ref{axm:OpRs_OpRs_EqAs}}]
\footnotesize $\Forall{a \w b}{\Ray}{\AngCongRs{a}{\Op{a}}{b}{\Op{b}}}$ 
\end{tcolorbox}
\\[-5mm]
\midrule
Orientation of angles is not defined in Hilbert's system of axioms.
&\text{\bf \Prop Axiom \refcolor{MyBlue}\ref{axm:EqAs_EqTs}} (Angles orientation){\bf \Prop .} \em 
Given two flags. If angular measures of corresponding angles are equal, then orientations of these flags should also be equal.
\\[2mm]
& 
\begin{tcolorbox}
[colback=MyBlue!5!white, frame hidden, enhanced, attach boxed title to top right={yshift=-5.2mm,xshift=-0.2mm}, 
boxed title style={size=small, frame hidden, colback=MyBlue},
title=\footnotesize\Axiom{\refcolor{white}\ref{axm:EqAs_EqTs}}]
\footnotesize $\Forall{X \w Y}{\Flag}{\AngCongRs{X_0}{X_1}{Y_0}{Y_1} \rightarrow {\mathbb O^{\rm +}}(X) = {\mathbb O^{\rm +}}(Y)}$ 
\end{tcolorbox}
\\[-5mm]
\bottomrule
\end{tabular*}
  \caption{Group IIIb: Axioms of angle congruence.}
  \label{tab:GIIIb}
\end{table*}

\newpage
\section{Conclusion} 

I would like to thank Andrei Rodin for interesting discussions on Euclid's methods and for providing valuable feedback on this article.

Despite the fact that our approach is completely based on the CoC as it was implemented in Coq, it is worth noting that other logics are also possible, which describe the constructive-deductive method of Euclid equally well. One such approach, called QHC logic, was recently proposed by Sergey Melikhov
\cite{Melikhov2015}.
I would like to thank Lev Beklemishev, Stepan Kuznetsov and Sergey Melikhov, who helped me understand the broader logical framework for the proposed constructive-deductive method in geometry. 

I am also grateful to Andrei Rodin, Julien Narboux and Yves Bertot for drawing my attention to references \cite{Avigad2009, Kellison2019, Beeson2016b, Beeson2019} on recent development in the field.

\bibliographystyle{unsrt}
\bibliography{neo-geo} 

\end{document}